\documentclass[12pt]{article}

\usepackage[T2A]{fontenc}
\usepackage[utf8]{inputenc}

\usepackage{amsfonts}
\usepackage{amssymb}
\usepackage{amsmath}
\usepackage{amsthm}

\def \le {\leqslant}
\def \ge {\geqslant}

\theoremstyle{plain}
\newtheorem{theorem}{Tеорема}
\newtheorem{lem}{Лемма}
\begin{document}

\begin{huge}
\centerline{\bf  Khintchine's singular systems }
\centerline{\bf and their applications}
 \end{huge}
\begin{Large}
 \vskip+1.0cm \centerline{\bf Nikolay G. Moshchevitin }
\end{Large}
 \vskip+2.0cm

 {\it
\hskip11.0cm To my wonderful teacher

\hskip11.0cm Valerii Vasilievich Kozlov

\hskip11.0cm in occasion of his 60-th birthday
}
 \vskip+2.0cm
This paper is a survey
of old and recent results related to Khintichine's singular matrices and their applications in the theory of Diophantine approximations.

In
 1926
A.Khintchine find out that the behaviour of two-dimensional Diophantine approximations differs cardinally from the behaviour
of one-dimensional Diophantine approximations. He proved the existence of two-dimensional real vectors which adimt 
"extremely good" rational approximations (he considered simultaneous  approximations as well as approximations for a linear form in two variables). Later this two-dimensional result was extended to the general systems of $n$ linear forms in $m$ variables.
The coressponding real matrices admitting   "extremely good" rational approximations are known now as 
Khintichine's singular matrices.
Most of the theory of singular matrices was constructed by A.Khintchine and V.Jarnik  in 1920 - 1950.
 
In our survey we discuss the definition of singular matrix and the existence theorems (Section \ref{vokrugo}).
 Section
\ref{ppnp} 
is devoted to the theory of subspaces of the best Diophantine approximations.
The phenomenon of the degeneracy of dimension of the subspaces of the best approximations is also discussed there.
We discuss one-dimensional Diophantine approximations (Section \ref{odp}) and simultaneous Diophantine approximations
 (Section  \ref{sszsdp}).
Section 
   \ref{sdty}
is devoted to the uniform and individual Diophantine characteristics.
We consider ingomogeneous Diophantine approximations in Section
  \ref{neopribl}.
Then in Section
 \ref{prore} 
we mention some problems of geometry of numbers related to the topic.
  A brief survey on classic and modern transference inequalities is given in Section
 \ref{tippo}.
In Section  \ref{HD}
we collect together results about the Hausdorff dimension of  sets
of singular matrices.
 
We would like to point out that Diophantine analysis related to singular matrices was the major tool for the
solution of the problem concerning  oscillating properties of the integral of a conditionally periodic function
proposed by V.V.Kozlov in 1978. This is the  subject of Section
 \ref{Zadako}.

The final Section \ref{dobavl} is devoted to some related problems including the application of a recent method due to Y.Peres and W.Schlag.
 
The  paper is written in Russian. English version should  appear in
"Russian Mathematical Surveys"  in the beginning of 2010.
\newpage

\begin{Large}
\centerline{\bf Сингулярные диофантовы системы А.Я. Хинчина}
\centerline{\bf  и их применение}
 \vskip+0.5cm \centerline{\bf Мощевитин Н.Г.\footnote{ Работа выполнена при поддержке гранта РФФИ № 09-01-00371а
и гранта поддержки ведущих научных школ  НШ-691.2008.1}  }
\end{Large}
 \vskip+1.0cm

 {\it
\hskip12.0cm моему  замечательному учителю

\hskip12.0cm Валерию Васильевичу Козлову

\hskip12.0cm к его шестидесятилетию
}
 \vskip+1.0cm

В 1926 году 
 А.Я.Хинчин обнаружил (см. \cite{HINS}), что в двумерных задачах теории диофантовых приближений
возникают явления, существенно отличающиеся от тех, которые имеют место  в задаче о приближении рациональными числами одного иррационального числа. В частности, он  конструктивно доказал существование двумерных вещественных векторов,  "аномально хорошо" приближаемых рациональными числами (в смысле линейной формы и в смысле совместных приближений). В дальнейшем конструкция работы 
\cite{HINS} была обобщена на случай произвольного количества линейных форм от произвольного количества целочисленных переменных, а  матрицы коэффициентов 
линейных форм, допускающих аномально хорошие приближения получили название  {\it сингулярных}.

В настояшей работе мы делаем обзор результатов из теории диофантовых приближений и ее приложений, тем или иным образом связанных с использованием сингулярных матриц А.Я.Хинчина. 
В основном, эта часть теории диофантовых приближений была  построена в статьях
А.Я.Хинчина 
\cite{H1925} - \cite{dobpere} и В.Ярника  \cite{JTBIL} - \cite{jP59}.
Основополагающая работа \cite{HINS} является одной из лучших (скорее, самой лучшей) работой А.Я.Хинчина в области теории диофантовых приближений.
 Многие статьи  А.Я.Хинчина были недавно переизданы в книге его избранных трудов   \cite{H}.
К сожалению, замечательные работы В.Ярника малодоступны,  на долгое время они были забыты.

Данный обзор начинается с обсуждения 
определения сингулярной матрицы, теорем существо-вания и связи с наилучшими диофантовыми приближениями (пункт \ref{vokrugo}).
В пункте \ref{ppnp} излагаются общие результаты о  подпространствах, порожденных наилучшими приближениями.
Отдельно рассматривается задача о совместных приближениях (пункт \ref{sszsdp}) и восходящая к В.Ярнику задача о связи 
равномерных и индивидуальных диофантовых характеристик (пункт \ref{sdty}). Там же дается усиление одного результата В.Ярника. 
В нескольких местах излагаются результаты автора о вырождении размерности подпространств, порожденных наилучшими приближениями.
Задачам о неоднородных  линейных диофантовых приближении посвящен пункт \ref{neopribl}.
В
пункте \ref{prore} рассматри-ваются задачи теории решеток, связанные с вопросами диофантовых приближений.
В пункте \ref{tippo} перечисляются классические и современные результаты, касающиеся принципа переноса.
В пункте \ref{HD} собраны известные автору  утверждения о размерности Хаусдорфа множеств сингулярных матриц.

Вопросы, связанные с применением теории диофантовых приближений естественным образом возникают в некоторых задачах классической механики, имеющих отношение к так называемой  проблеме "малых знаменателей",
например, в теории  возмущений условно-периодических движений
(см. \cite{kUMN}). В частности оказалось, что сингулярные системы А.Я.Хинчина могу быть применены для построения систем с "быстро разбегающимися траекториями", что было сделано В.В.Козловым и Н.Г.Мощевитиным  \cite{kOMo}.  Более того, многомерный диофантов анализ оказался основным инструментов а решении задачи В.В.Козлова об осцилляции интеграла условнопериодической функции. Этой задаче посвящен пункт \ref{Zadako}
настоящего обзора.

В добавлении (пункт \ref{dobavl}) мы  хотели упомянуть ряд утверждений и теорем, имеющих отношение  к некоторым 
затронутым вопросам
и их развитию.
Это, во-первых, результаты, которые получаются с помощью оригинального метода, предложенного недавно Ю.Пересом и В.Шлагом \cite{PS}.  Во-вторых, это основные понятия и некоторые результаты теории выигрышных множеств,
построенной В.М.Шмидтом \cite{SCH1}.

Некоторые вопросы, затронутые в настоящей статье
(равно как и многие другие результаты теории диофантовых приближений)
нашли свое отражение в замечательном обзоре М.Вальдшмидта \cite{WALL}.
Вопросам теории диофантовых приближений посвящены ставшие классическими монографии Дж.Ф.Коксмы \cite{Koks},
Дж.В.С.Касселса \cite{Cassil}  и В.М.Шимдта \cite{SCH}.

Кроме перечисленных выше, в настоящем обзоре затронуты и некоторые другие задачи. Это, например, задача о диофантовых приближениях с положительными числами (пункт \ref{neotri}) и проблемы существования
матриц с чрезвычайно синуглярными диофантовыми свойствами (пункт \ref{specivid}).

Автор считает своей приятной обязанностью поблагодарть всех участников теоретико-числовых семинаров и спецкурсов механико-математического факультета МГУ и Независимого московского университета за многочисленные обсуждения
задач и теорем из  настоящего обзора. В особенности автор благодарен И.П.Рочеву и О.Н.Герману.

 \section{Вокруг определения сингулярной системы.}\label{vokrugo}

Всюду ниже 
$
{\bf x} = (x_1,...,x_m)
$
обозначает целочисленный вектор. Через $|\cdot ||$  обозначено расстояние до ближайшего челого числа.
Рассматривается матрица
\begin{equation}
\Theta =\left(
\begin{array}{ccc}
\theta_{1}^1&\cdots&\theta_{1}^m\cr
\cdots&\cdots&\cdots\cr
\theta_{n}^1&\cdots&\theta_{n}^m
\end{array}
\right),
\label{1f}
\end{equation}
состоящая из вещественных чисел
$\theta_j^i,\,\, 1\le j\le n,\,\, 1\le i \le m$,
и соответствующая ей система линейных форм

\begin{equation}
{\bf L }({\bf x})= {\bf L}_\Theta ({\bf x}) = \{L_j ({\bf x}),\,\,\, 1\le j\le n\},
 \,\,\,\,\,
L_j ({\bf x})= \sum_{i=1}^m \theta_{j}^i x_i.
\label{2f}
\end{equation}
 Через  $^t\!\Theta$  будем обозначать матрицу,
транспонированную к матрице  $\Theta$.

Согласно принципу Дирихле для любой матрицы $\Theta$ и при любом  вещественном $t\ge 1$ система диофантовых неравенств
$$
 \max_{1\le j\le n}||L_j({\bf x})||\le \frac{1}{t},\,\,\,\,\,\,\,\, 0<\max_{1\le i\le m}|x_i|\le t^{\frac{n}{m}}
$$
имеет целочисленное решение ${\bf x}\in \mathbb{Z}^m$.

\subsection{Определение Хинчина}

Сначала мы сформулируем определение сингулярной системы именно в том виде, в каком его давал А.Я.Хинчин (см.
\cite{HRSU}). Матрица  $ \Theta $
(или, как писал сам А.Я.Хинчин, набор вещественных чисел
$\theta_{j}^i, \,\, 1\le i\le m,\,\,1\le j\le n$)
 представляет из себя {\it  сингулярную систему} если при любом вещественном $\varepsilon >0$
найдется некоторое $ t_0= t_0(\varepsilon)$ такое, что при всяком $t\ge t_0$ система диофантовых неравенств
$$
 \max_{1\le j\le n}||L_j({\bf x})||\le \frac{1}{t},\,\,\,\,\,\,\,\, 0<\max_{1\le i\le m}|x_i|<\varepsilon t^{\frac{n}{m}}
$$
имеет целочисленное решение ${\bf x}\in \mathbb{Z}^m$.

Отметим, что матрицу $ \Theta $, не являющуюся сингулярной, А.Я. Хинчин называл {\it регулярной}
(он использовал	терминологию
{\it регулярная система} чисел
$\theta_{j}^i, 1\le i\le m,\, 1\le j\le n$).

Таким образом, матрица  $ \Theta $ является регулярной  если найдется положительное $\mu (\Theta )$ такое, что
для
  некоторой последовательности  положительных чисел $t_\nu$, стремящейся к бесконечности, в областях
пространства $\mathbb{R}^n$, задаваемых неравенствами
\begin{equation}
 \max_{1\le j\le n}||L_j({\bf x})||\le \frac{\mu (\Theta) }{t_\nu},\,\,\,\,\,\,\,\, 0<\max_{1\le i\le m}|x_i|\le \mu (\Theta )\cdot t_\nu^{\frac{n}{m}}
\label{regul}
\end{equation}
отсутствуют ненулевые целые точки.

В этой статье мы будем использовать, в основном, матричную терминологию. Тем не менее, иногда удобно употреблять понятие сингулярной системы вещественных чисел.

\subsection{ Теоремы существования.}\label{tees}

При $n=m=1$, как легко видеть, сингулярными системами в смысле определения Хинчина являются только рациональные числа $\theta_{1}^1$ (это
утверждение мы прокомментируем ниже в пункте \ref{cd}). Впервые существование сингулярных систем было доказано А.Я. Хинчиным   в 1924 году  в работе
\cite{HINS} при $m=2, n=1$ и при  $m=1,n=2$. Здесь мы приведем формулировки соответствующих утверждений (Hilfssatz I
и Satz 2 из \cite{HINS}, см. также книгу
\cite{Cassil}, гл. V, теорема XIV).

\begin{theorem}\label{T1}
Пусть $\psi (t)$ положительная, непрерывная, убывающая к нулю при $t\to +\infty$ функция вещественного переменного $t$.
Тогда существуют два линейно независимых вместе с единицей над $\mathbb{Z}$ вещественных числа $\theta^1,\theta^2$, такие что для каждого
достаточно большого  $t$ система диофантовых неравенств
$$
||x_1\theta^1+x_2\theta^2 || <\psi (t),\,\,\,\,\,\,\,\, 0<\max_{j=1,2}|x_j|<t
$$
имеет целочисленное решение $(x_1,x_2)\in \mathbb{Z}^2$.
\end{theorem}

\begin{theorem}\label{T2}
Пусть $\psi (t)$ положительная, непрерывная, убывающая к нулю при $t\to +\infty$ функция вещественного переменного $t$, и пусть  функция $t\mapsto t\psi (t)$ монотонно возрастает к бесконечности при  $t\to +\infty$.
Тогда существуют два линейно независимых вместе с единицей над $\mathbb{Z}$ вещественных числа $\theta_1,\theta_2$, такие что для каждого
достаточно большого  $t$ система диофантовых неравенств
$$
\max_{j=1,2}
||x\theta_j || <\psi (t),\,\,\,\,\,\,\,\, 0<x<t
$$
разрешима в $x\in \mathbb{Z}$.

\end{theorem}
Итак, в случае $m=2, n=1$ теорема \ref{T1} устанавливает существование сингулярных систем с  соответствующей функцией $\psi (t)$, сколь угодно быстро
стремящейся к нулю. 
В случае $m=1,n=2$  теорема \ref{T2}  устанавливает существование сингулярных систем только в случае, когда соответствующая  функция $t\mapsto t\psi (t)$ монотонно возрастает к бесконечности.

В связи с этим уместно следующее определение. Пусть 
непрерывная
функция $\psi (t)$ монотонно убывает к нулю при $t\to +\infty$ и $ \psi
(t) =o(t^{-\frac{m}{n}}),\, t\to +\infty$. Матрицу $ \Theta $ 
(или набор $mn$ вещественных чисел)
мы будем называть {\it $\psi$-сингулярным} если
 при всяком  достаточно большом $t $ система диофантовых неравенств
$$
 \max_{1\le j\le n}||L_j({\bf x})||\le \psi (t) ,\,\,\,\,\,\,\,\, 0<\max_{1\le i\le m}|x_i|<  t
$$
имеет целочисленное решение ${\bf x}\in \mathbb{Z}^m$.

Сформулированная выше теорема \ref{T1} непосредственно обобщается на случай произвольной размерности следующим образом.

\begin{theorem}\label{T3}
Пусть $n$ -- произвольное натуральное число, а $m$ -- натуральное число, не меньшее, чем $2$. Пусть $\psi (t)$
положительная, непрерывная, убывающая к нулю при $t\to +\infty$ функция вещественного переменного $t$.
Рассмотрим множество ${\cal M}\subset \mathbb{R}^{mn}$,
 состоящее
из матриц $\Theta$,  таких что

$\bullet$
\,\,  числа  $\theta^i_j,\,\,\,
 1\le i \le m,\,\,\, 1\le j\le n$ линейно независимы вместе с единицей над $\mathbb{Z}$,

$\bullet$\,\, матрица $ \Theta $ является $\psi$-сингулярной.

Тогда для любого отрытого множества ${\cal G} \subset \mathbb{R}^{mn}$ пересечение
${\cal M}\cap{\cal G}$  имеет мощность континуум.
\end{theorem}

 Отметим, что при  $m=1$  
дело обстоит несколько по другому.
Итак, пусть $m=1$, и $ \Theta =\{ \theta_{1},...,\theta_n\} $ есть некоторый  набор вещественных чисел . Условимся обозначать через ${\rm
dim}_\mathbb{Z}\Theta $  максимально возможное количество линейно независимых над $\mathbb{Z}$ чисел из набора $\theta_{1},...,\theta_n,1$.
Сформулируем обобщение теоремы \ref{T2}.

\begin{theorem}\label{T4}
Пусть $m=1$ и $n\ge 2$.

{\rm (i)} Пусть $\psi (t)$ положительная, непрерывная, убывающая к нулю при $t\to +\infty$ функция вещественного переменного $t$, причем
\begin{equation}
\lim_{t\to\infty}{\psi (t)}\cdot{t} =+\infty .\label{psi1}
\end{equation}

Рассмотрим множество ${\cal M}\subset \mathbb{R}^{n}$   наборов вещественных чисел $ \Theta
=\{ \theta_{j},\,\,\, 1\le j\le n\} $, таких что

$\bullet$\,\, ${\rm dim}_\mathbb{Z}\Theta =n+1$,

$\bullet$\,\, система $ \Theta $ является $\psi$-сингулярной.

Тогда для любого отрытого множества ${\cal G} \subset \mathbb{R}^{n}$ пересечение
${\cal M}\cap{\cal G}$  имеет мощность континуум.

{\rm (ii) }  Пусть  $n\ge 2$, и  для положительнозначной функции
  $\psi (t)$
    выполнено
 \begin{equation}
 \limsup_{t\to+\infty}
{\psi (t)}\cdot{t} <+\infty . \label{psi3}
\end{equation}
Тогда, если набор $ \Theta =\{ \theta_{1},...,\theta_n\} $ представляет из себя
  $\psi$-сингулярную систему, то $ {\rm dim}_\mathbb{Z}\Theta \le 2 $.
  Более того, если
 \begin{equation}
 \limsup_{t\to+\infty}
{\psi (t)}\cdot{t} =0 , \label{psi4}
\end{equation}
то $ {\rm dim}_\mathbb{Z}\Theta =1 $, то есть все числа $\theta_{j},\,\,\, 1\le j\le n$ суть рациональные.
\end{theorem}

Теорема \ref{T2} и часть (i)  теоремы \ref{T4} имеются у В.Ярника \cite{J59}, частный случай  рассмотрен у К.Шаботи и Е.Лютц \cite{LUTZ}.
Отметим, что В.Ярник формулирует и доказывает несколько более сильный результат: существование матриц
$\Theta$,  состоящих из вещественных чисел  $\theta_{j}^i$,
 {\it алгебраически  } независимых в совокупности.  

Утверждение  (ii)  теоремы \ref{T4} очевидным образом вытекает из  следcтвия 5 пункта \ref{nvnp}, которое тоже принадлежит В.Ярнику
(см. \cite{JTBIL} и \cite{J59}), результат в предположении (\ref{psi4}) мгновенно  получается из
 следствия 1 пункта \ref{nvnp}.
Таким образом, в случае $m =1, n\ge 2$   при условии (\ref{psi3}) нетривиальных сингулярных систем   не существует.

Теорема \ref{T4}
 допускает некоторое уточнение, имеющееся у Ж.Леска \cite{les}.

\begin{theorem}\label{TTT4} {\rm (Ж.Леска \cite{les})}
Пусть $m=1$ и $n\ge 2$.
Пусть $\psi (t)$ положительная, непрерывная, убывающая к нулю при $t\to +\infty$ функция вещественного переменного $t$, 
и
$$
\limsup_{t\to\infty}{\psi (t)}\cdot{t} =+\infty .
$$
Тогда множество $\psi$-сингулярных  наборов  $ \Theta
=\{ \theta_{j},\,\,\, 1\le j\le n\} $, состоящее из алгебраически независимых 
  вещестенных чисел  в пересечении с произвольным 
 отрытым множеством ${\cal G} \subset \mathbb{R}^{n}$  
  имеет мощность континуум.
\end{theorem}

Теорему \ref{TTT4}  уместно сопоставить с формулой В.Ярника 
(\ref{formuja}) и с
теоремой \ref{T7} из пункта \ref{nvnp}.
Отметим, что
Ж.Леска в \cite{less} получил $p$-адический вариант теоремы \ref{T3}.

 Приведем теорему существования, доказанную  А.Апфельбеком \cite{apf}:

\begin{theorem}\label{Ta}{\rm (А.Апфельбек \cite{apf})}
Пусть  $n,m\ge 2$.  Пусть монотонно убывающая к нулю непрерывная функция  $\psi (t)$ такова, что
 функция  $t\mapsto t\psi (t)$ тоже  убывает.  Тогда найдется $\psi$-сингулярная матрица  $\Theta$,  состоящая из линейно независимых вместе с единицей  над  $\mathbb{Z}$
чисел  $\theta^i_j$, такая что матрица  $^t\!\Theta$  тоже $\psi$-сингулярна.
\end{theorem}

На самом деле у А.Апфельбека  условие линейной независимости несколько слабее, чем  в сформулированной теореме \ref{Ta}.   Условие монотонности функции $t\mapsto t\psi (t)$ не является принципиальным.

Использование в определении сингулярной системы прилагательного "сингулярный" объясняется тем, что такие системы при каждом значении размерностей $m,n$, как отмечал сам А.Я. Хинчин (см. \cite{HRSU}), образуют  множество нулевой меры Лебега. Этот факт связан с простым применением леммы Бореля-Кантелли. Подробное доказательство можно найти в гл.V, \S 7 книги \cite{Cassil}. 
Более сильный результат был получен Г.Давенпортом и В.М.Шмидтом  \cite{DSdir1}.  Мы обсудим его ниже в пункте 
  \ref{daves}.

\subsection{Сингулярные системы и наилучшие приближения.}\label{ssinp}

Для целочисленного вектора $ {\bf x} = (x_1,...,x_m) \in \mathbb{Z}^m$ рассмотрим величины
$$
M({\bf x}) = \max_{1\le i\le m} |x_i|,\,\,\,\,\,\, \zeta ({\bf x}) =  \max_{1\le j\le n}||L_j({\bf x})||.
$$
Мы будем называть целочисленную точку $ {\bf x} = (x_1,...,x_m) $ {\it наилучшим приближением} для матрицы $\Theta$  если
\begin{equation}
\zeta ({\bf x})=\min_{{\bf x}'} \zeta ({\bf x}'), \label{min}
\end{equation}
 где минимум берется по всем ненулевым целым точкам  $ {\bf x}' = (x_1',...,x_m')
\in \mathbb{Z}^m $, подчиненным условию
\begin{equation}
0<  M({\bf x}')\le M({\bf x}). \label{min1}
\end{equation}
При таком определении оказывается, что наилучшие приближения будут встречаться парами: вместе с каждой точкой $ {\bf x}$, являющейся наилучшим
приближением, наилучшим приближением будет также точка $ -{\bf x}$. Однако, если рассмотреть такую пару наилучших приближений $ \pm{\bf x}$, то
для нее величины $M({\bf x})$ и $\zeta ({\bf x})$ будут определяться однозначно и вне зависимости от выбора знака $\pm$.

Отметим следующее. Вообще говоря, может оказаться что для двух  целочисленных точек $ {\bf x}_1\neq\pm {\bf x}_2$ с одинаковым значением $
M({\bf x}_1)=M({\bf x}_2) \neq 0$ выполняется
\begin{equation}
 \max_{1\le j\le n}||L_j({\bf x}_1)||=\max_{1\le j\le n}||L_j({\bf x}_2)||.
\label{two}
\end{equation}
  Но если, например,  предположить,
что все числа $\theta_{j}^i,\,\,\, 1\le i \le m,\,\,\, 1\le j\le n$ линейно независимы вместе с единицей над $\mathbb{Z}$, то равенство
(\ref{two}) невозможно. Таким образом, в случае, когда все числа $\theta_{j}^i,\,\,\, 1\le i \le m,\,\,\, 1\le j\le n$ линейно независимы вместе
с единицей над $\mathbb{Z}$, наилучшие приближения можно перенумеровать и расположить в виде  бесконечной последовательности
$$
\pm {\bf x}_1,\pm {\bf x}_2,..., \pm{\bf x}_\nu,  \pm {\bf x}_{\nu+1}, ...
$$
так, чтобы
\begin{equation}
M({\bf x}_1)< M({\bf x}_2)<...<M({\bf x}_\nu)<M({\bf x}_{\nu+1})<... , \label{1}
\end{equation}
\begin{equation}
\zeta({\bf x}_1)> \zeta({\bf x}_2)>...>\zeta({\bf x}_\nu)>\zeta({\bf x}_{\nu+1})>... . \label{2}
\end{equation}

Для краткости мы будем использовать обозначения
$$
M_\nu = M({\bf x}_\nu),\,\,\,\,\,\zeta_\nu = \zeta({\bf x}_\nu).
$$

Иногда в дальнейшем нам придется иметь дело с наилучшими приближениями для матрицы  $\Theta$, не удовлетворяющих условию линейной независимости.  В этом случае мы не
можем определить последовательность точек ${\bf x}_\nu$ однозначно. Тем не менее, последовательности (\ref{1},\ref{2})  величин $M({\bf x}_\nu),
\zeta({\bf x}_\nu)$ однозначно определяются (только каждому из значений $M({\bf x}_\nu), \zeta({\bf x}_\nu)$, вообще-говоря, может
соответствовать несколько пар целочисленных точек $\pm {\bf x}_\nu$ и последовательности, вообще-говоря, могут быть конечными).

В связи со сказанным выше мы даем следующее определение.
Матрицу $\Theta$ мы будем называть {\it правильной}, если последовательности (\ref{1},\ref{2}) бесконечны, и для каждого
достаточно большого  значения  $\nu$  векторы  ${\bf x}_\nu$  определены однозначно с точностью до знака.

Надо 
сделать еще одно замечание. Определим единичные векторы
$$
{\bf e}^k=
\left(
\begin{array}{c}
0\cr
\vdots\cr
1\cr
\vdots
\cr
0
\end{array}
\right),\,\,\,
1\le k \le n
$$
(  у $k$-того вектора единица стоит на $k$-том месте, остальные компоненты равны нулю)
 и рассмотрим в $\mathbb{R}^n$ набор, состоящий из   $m+n$
векторов
\begin{equation}
\theta^1,\theta^2,...,\theta^m
,{\bf e}^1,{\bf e}^2,...,{\bf e}^n
\label{veee}
\end{equation}
(Здесь $\theta^j$ обозначаает $j$-тый столбец матррицы $\Theta$). Легко видеть, что последовательности 
наилучших приближений (\ref{1},\ref{2})  являются бесконечными в том и только  том случае, когда
набор векторов (\ref{veee}) 
является линейно независимым над  $\mathbb{Z}$.  В этом случае В.Ярник предлагал говорить, что матрица  $\Theta$ является {\it невырожденной} (в  настоящей статье мы тоже будем использовать эту терминологию).
В противном случае
последовательности (\ref{1},\ref{2}) будут конечными, и для последнего значения индекса $\nu$  будет выполнено $\zeta_\nu = 0$.

Если  хотя бы одно из чисел  $m$  или  $n$  равно единице, то матрица 
 $\Theta$  представляет из себя один столбец  или одну строку и может быть отождествлена просто с набором из $m$  
или $n$
вещественных чисел, и можно говорить о величине  ${\rm dim }_\mathbb{Z}\Theta$, как о характеристике этого набора вещественных чисел.  Именно в этом смысле обозначение  ${\rm dim}_\mathbb{Z}\Theta$
 употреблялось в теореме 
\ref{T4}. В этом же смысле это обозначение нам понадобится и в дальнейшем.
В частности, получается, что  для матрицы $\Theta$
величина ${\rm dim}_\mathbb{Z}\Theta$
определена только при  $m=1$  или  $n=1$.
Например,
при $m = 1$  невырожденность матрицы  $\Theta$ означает,
 что хотя бы одно из чисел  $\theta_j =\theta^1_j$  является иррациональным, то есть
 ${\rm dim }_\mathbb{Z} \Theta \ge 2$.  При $ n = 1$  невырожденность матрицы  $\Theta$  означает, что числа  $1,\theta_1^1,...,\theta_1^m$  линенйно независимы над  $\mathbb{Z}$,  то есть
 ${\rm dim}_\mathbb{Z}\Theta = m+1$.

Нам понадобится еще одна характеристика независимости.
 Для матрицы $\Theta$  через  ${\rm DIM}_\mathbb{Z}\Theta$
обозначим максимальное количество векторов из набора
(\ref{veee})  таких, что они линейно независимы над  $\mathbb{Z}$.
  Видим, что при  $n=1$  всегда выполнено
\begin{equation}
{\rm DIM}_\mathbb{Z}\Theta ={\rm dim}_\mathbb{Z}\Theta,
\label{didid}
\end{equation}
 а при  $m=1$  величина  ${\rm DIM}_\mathbb{Z}\Theta$
 может принимать только два значения
$n$  и $n+1$.  
В последнем случае вместо равенства  (\ref{didid}), естественно, имеет место равенство
\begin{equation}
{\rm DIM}_\mathbb{Z}^t\!\Theta ={\rm dim}_\mathbb{Z}\,^t\!\Theta,
\label{dididi}
\end{equation}

Ясно, что для произвольной матрицу  $\Theta$ имеют место неравенства

\begin{equation}
n \le {\rm DIM}_\mathbb{Z}\Theta \le n+m.
\label{nnm}
\end{equation}
Невырожденнось матрицы $\Theta$  эквивалентна равенству
$${\rm DIM}_\mathbb{Z}\Theta = m+n.
$$

Поясним взаимосвязь между правильностью и невырожденностью матрицы $\Theta$. Если $m=1$ или $n=1$, то из невырожденности следует правильность. В других случаях это неверно. Приведем пример.
 
 Рассмотрим иррациональное число  $\xi$.  Пусть  $m=n=2$  и матрица $\Theta$ имеет вид
$$
\Theta = 
 \left(
\begin{array}{cc}
\theta^1_1& \theta_1^2\cr\theta_2^1 & \theta^2_2
\end{array}
\right) =
\left(
\begin{array}{cc}
\xi& 0\cr 0& \xi
\end{array}
\right)
.$$
Тогда
$$
{\rm DIM}_\mathbb{Z} \Theta =
{\rm DIM}_\mathbb{Z}\, ^t\!\Theta = 4,
$$
и матрица является невырожденной.
Но, как легко видеть,
каждому значению $\nu$ соответсвтуют векторы
 $$
(\pm q_\nu,\pm q_\nu ),\,\,\, (\pm q_\nu,0),\,\,\, (0,\pm q_\nu)
$$
 (где $ q_\nu$ есть знаменатель  некоторой подходящей к $\xi$ дроби), с одинаковыми значениями величин $M_\nu,\zeta_\nu$.
Однако, набор вещественных чисел
$1, \theta^i_j,\,\,1\le i,j \le 2$  линейно зависим над  $\mathbb{Z}$,
 поскольку среди чисел  $\theta^i_j$ имеются нули.

 Продолжим анализ понятия наилучшего приближения. 
Обозначим через $ {\bf y}_\nu =(y_{1,\nu},...,y_{n,\nu})\in \mathbb{Z}^n$ целочисленный вектор, состоящий из чисел $y_{j,\nu}$ таких, что
$$
||L_j({\bf x}_\nu)||= |L_j({\bf x}_\nu)+y_{j,\nu}|.
$$
Нам понадобится обозначение
$$
{{\bf z}_\nu} = (x_{1,\nu},...,x_{m,\nu}, y_{1,\nu},..., y_{n,\nu} ) \in \mathbb{Z}^{d},\,\, d = m+n
$$
для "расширенного" вектора наилучших приближений.

То, что "расширенный" вектор
$$
{\bf z} = (x_1,...,x_m,y_1,...,y_n)
$$
 является наилучшим приближением, в частности,
означает, что внутри параллелепипеда
\begin{equation}
\{
{\bf z}'= (x'_1,...,x'_m , y_1',...,y_n'):\,\,\,\,
M({\bf x}')\le M({\bf x}),\, \max_{1\le j\le n}|L_j({\bf x}') +y_j'|
\le  \max_{1\le j\le n}|L_j({\bf x})+y_{j}|\}
\label{dopara}
\end{equation}
нет других целых точек кроме ${\bf 0}$.
Более того, внутри  параллелепипедов вида
\begin{equation}\{
{\bf z}'= (x'_1,...,x'_m , y_1',...,y_n'):\,\,\,\,
M({\bf x}')\le M({\bf x}_{\nu+1}),\, \max_{1\le j\le n}|L_j({\bf x}) + y_j'|
\le  \zeta_\nu = \max_{1\le j\le n}|L_j({\bf x}_\nu)+y_{j,\nu}|\}
\label{dopara1}
\end{equation}
нет других целых точек кроме ${\bf 0}$. 
Отсюда по теореме Минковского о выпуклом теле получаем, что
\begin{equation}\label{minbody}
\zeta_\nu^n M_{\nu+1}^m\le 1.
\end{equation}
Кроме того, отметим два простых и важных свойства.

1. Каждый расширенный вектор  ${\bf z}_\nu$
 является {\it примитивным}, то есть
$$
\text{н.о.д.}
(x_{1,\nu},...,x_{m,\nu}, y_{1,\nu},..., y_{n,\nu} )
=1.
$$

2. Любые два последовательных вектора наилучших приближений
 ${\bf z}_\nu, {\bf z}_{\nu+1}$
 дополнимы до базиса решетки  $\mathbb{Z}^d$.

Свойство матрицы $\Theta$ быть сингулярной легко может быть переформулировано на языке наилучших приближений.

{\bf Предложение 1.}\,\,\,{\it 
Пусть 
для непрерывной и монотонной функции  $\psi(t)$  выполнено
 $\psi (t) = o(t^{-m/n})$  при  $t\to +\infty$.
Невырожденная  матрица $\Theta$ 
является  $\psi$-сингулярной тогда и только тогда, когда для всех достаточно больших значений $\nu$
выполняется
\begin{equation}
\zeta_\nu \le \psi (M_{\nu+1}).
 \label{3}
\end{equation}
 }

Связь свойства системы вещественных чисел $\Theta$ быть сингулярной с понятием наилучших приближений в работах А.Я.Хинчина явно не отмечалась.
Эта связь неявно используется в   работах В.Ярника и в книге Дж.В.С.Касселса \cite{Cassil}. 
В частности,
у В.Ярника (см. \cite{JTBIL},\cite{J41},\cite{JRUS},\cite{J59})
используется  кусочно-постоянная функция
\begin{equation}
\psi_\Theta ( t)  = \,\,\,\,\,\min_{{\bf x}\in\mathbb{Z}^m: \, 0<M({\bf x})\le t
}\,\,\,\,\,\max_{1\le j\le n}||L_j({\bf x})||
 ,
\label{fj}
\end{equation}
а в книге \cite{Cassil}  в ряде доказательств (см. \cite{Cassil}, гл.  V, \S\S 6,7)
используется функция
$$
\eta(\rho ) = \,\,\,\,\,\min_{{\bf x}\in \mathbb{Z}^m: \, 0<
 |x_1|^2+...+ |x_m|^2\le\rho^2
}\,\,\,\,\,\max_{1\le j\le n}||L_j({\bf x})||.
$$
Точки разрыва этих функций фактически определяют наилучшие приближения (у В.Ярника  в $\sup$-норме, как это сделано в настоящей статье, а у Дж.В.C.Касселса  в евклидовой
норме).

Заметим, наконец, что матрица  $\Theta$ невырождена тогда и только тогда, когда
 функция  $\psi_\Theta(t)$
 никогда не обращается в ноль при $ t \ge 1$.
Невырожденная матрица  $\Theta$  является $\psi$-сингулярной если и только если при всех достаточно больших значениях $t$
для функции Ярника (\ref{fj}) выполнено
$$
\psi_\Theta(t) \le \psi (t).
$$
\section{Подпространства, порожденные наилучшими \\ приближениями}
\label{ppnp}

В этом пункте мы обсудим свойства линейных подпространств в 
$\mathbb{R}^d$,
в которых  содержатся  векторы  наилучших приближений  ${\bf z}_\nu$.

\subsection{ Оценки размерностей}\label{or}

Под размерностью ${\rm dim}\,\Lambda$ решетки $\Lambda$ 
 мы  будем подразумевать размерность минимального линейного подпространства
${\rm span}\, \Lambda \subset \mathbb{R}^d$, содержашего решетку $\Lambda$.
Pассмотрим в пространстве $\mathbb{R}^{d}$ линейное подпространство $\pi $. 
Пересечение $\pi \cap \mathbb{Z}^{d}$ представляет из себя некоторую решетку  $\Lambda =\Lambda (\pi )$ (возможно  состоящую из одной точки ${\bf 0}$).   
Подпространство $\pi \subseteq\mathbb{R}^{d}$ 
мы будем называть {\it вполне рациональным}, если 
$${\rm dim }\, \pi = {\rm dim }
\Lambda (\pi ).
$$
Далее, для данного линейного подпространства  $\pi$
определим  $\hbox{\got H}(\pi )$  как максимальное вполне рациональное подпространство,
{\it содержащееся}
в  $\pi$.  Определим также  $\hbox{\got R} (\pi )$  как минимальное вполне рациональное подпространство,
 {\it  содержащее} подпространство  $\pi$.  Имеем
$$
\hbox{\got H}(\pi )\subseteq \pi \subseteq\hbox{\got R}(\pi ) .
$$

Определим теперь векторы
$$
\overline{\theta}_i =(\theta_i^1,...,\theta_i^m,0,...,0,1,0,...,1),\,\,\, 1\le i \le n
$$
(единица
 здесь стоит на $m+i$-том месте).

 Рассмотрим в  $\mathbb{R}^d,\, d = m+n$
линейное подпространство
 ${\cal N} (\Theta)$,  натянутое на векторы  $\overline{\theta}_1,...,\overline{\theta}_m$,
и ортогональное дополнение к нему   ${\cal L} (\Theta )$. 
Видим, что
  $$
{\rm dim }
\,{\cal N}(\Theta ) = n,\,\,\, {\rm dim }\,{\cal L} (\Theta ) = m.
 $$
Рассмотрим теперь подпространства
$$
\hbox{\got H}_\Theta = \hbox{\got H} ({\cal L}(\Theta)), \,\,\, 
 \hbox{\got R}_\Theta = \hbox{\got R} ({\cal L}(\Theta))
.
$$
 Видим, что
\begin{equation}\label{razmeri}
{\rm DIM}_\mathbb{Z}\Theta  +{\rm dim} \,\hbox{\got H}_\Theta =d ,\,\,\,
{\rm DIM}_\mathbb{Z}\,^t\!\Theta  = {\rm dim} \,\hbox{\got R}_\Theta .
\end{equation}
 В частности, невырожденность матрицы $\Theta$, согласно   первому из неравенств (\ref{razmeri}),
равносильна условию 
$\hbox{\got H}_\Theta=\{{\bf 0}\}.$

Следующее утверждение общеизвестно. Оно фактически доказано
 В.Ярником в \cite{JRUS}.   Случай $m=1$  есть у  Дж.Лагариаса \cite{LL}. 
Случай $m=2, n = 1$  имеется у Г.Давенпорта и В.М.Шмидта в \cite{DS}.
(см. также работы 
автора \cite{MUMN},\cite{ME}).

 Для правильной матрицы $\Theta$ рассмотрим величину
$$
R(\Theta )=\min \{ r:\text{найдется подрешетка} \,\Lambda \subseteq \mathbb{Z}^{n+m}, {\rm dim } \,\Lambda  = r\,\text{и}
\,\nu _0 \in \mathbb{N}\,\,\text{такие, что}\,\,{\bf z}_ \nu  \in \Lambda   \hspace{2mm} \forall \nu \ge \nu_0\}. $$
Пусть $\Lambda_\Theta$  есть именно та решетка, которая фигурирует в определении величины  $R(\Theta)$
и $\pi_\Theta ={\rm span}\,\Lambda_\Theta$.
Положим
\begin{equation}\label{ones}
K(\Theta) ={\rm dim}(\pi_\Theta \cap {\cal L}(\Theta))\ge 1.
\end{equation}
 Последнее
неравенство имеет место постольку, поскольку векторы наилучших приближений, лежащие в подпространстве  $\pi_\Theta$
 приближаются сколь угодно близко к подпространству  ${\cal L}(\Theta )$.

\begin{theorem}\label{NT1} 
Для  правильной  матрицы  $\Theta$   выполняется  

{\rm (i)} 
\,\ $2\le R(\Theta) \le{\rm DIM}_\mathbb{Z}\,^t\!\Theta 
$;

{\rm (ii)}
 если  $m=1$, то
$
R(\Theta) ={\rm DIM}_\mathbb{Z}\,^t\!\Theta = {\rm dim}_\mathbb{Z}\Theta  
$;

{\rm (iii)}\, если  $ m >n$, то
$
R(\Theta)\ge 3
$;

 {\rm (iv)}
если $K(\Theta ) =1$  
и  $R(\Theta) > 2$,
то 
$m<n$.
\end{theorem}

Поясним доказательства неравенств теоремы \ref{NT1}.

{\bf 1.}
Нижняя оценка утверждения (i)  
следует из линейной независимости векторов  ${\bf z}_\nu,{\bf z}_{\nu+1}$.

 Для того чтобы доказать
 верхнее неравенство из  (i) и утверждение (ii)  введем несколько обозначений, которые понадобятся не только сейчас, но и в дальнейшем.

Расстояние между множествами ${\cal A},{\cal B}\subset \mathbb{R}^{d}$
 будем обозначать
 ${\rm dist}({\cal A},{\cal B})$
(нам удобно считать, что в данном случае имеется ввиду расстояние  в sup-норме).
 Если у нас имеется некоторая решетка $\Lambda\subset \mathbb{R}^{d}$, причем
${\rm dim}\,\Lambda <d$,
то расстояние 
 от точек множества  $\mathbb{Z}^{d}\setminus \Lambda$ до линейного подпространства
${\rm span}\, \Lambda$ отделено от нуля. Его мы будем обозначать
 \begin{equation}
\rho (\Lambda ) = {\rm dist }\, (\mathbb{Z}^{d}\setminus \Lambda, {\rm span}\, \Lambda) >0.
\label{rho}
\end{equation}
Если $\pi \subset\mathbb{R}^{d}$  есть вполне рациональное линейное подпространство, то будем обозначать
\begin{equation}\label{ropi}
\rho (\pi ) 
=\rho (\pi\cap \mathbb{Z}^{d}).
\end{equation}

{\bf 2.}
Итак,  теперь мы можем прокомментировать доказательство верхнего неравенства утверждения   (i).

Поскольку $\rho ( \hbox{\got R}_\Theta)>0$, но $ 
{\rm dist} ({\bf z}_\nu,  \hbox{\got R}_\Theta)\le
{\rm dist} ({\bf z}_\nu, {\cal L} (\theta))\to 0, \nu\to+\infty$,
то при всех достаточно больших значениях $\nu$ получается
   ${\bf z}_\nu \in  \hbox{\got R}_\Theta\cap\mathbb{Z}^{n+1}, $  и все доказано.

{\bf 3.}  Прокомментриуем доказательство утверждения  (ii).
Заметим, что при условии  $m=1$  пространство  ${\cal L}(\Theta )$   одномерно, $d=m+1$.

Если ${\rm dim}_\mathbb{Z}\Theta= r$,
то это означает, что вектор $\Theta$  лежит в некотором вполне рациональном линейном подпространстве $\pi\subset\mathbb{R}^{n+1}$  размерности $r$ 
( которое и есть $\hbox{\got R}_\Theta$)
 и не лежит ни в каком вполне рациональном подпространстве меньшей размерности. 

 Оценка  $ R(\Theta ) \le r$  следует из  верхнего неравенства пункта  (i).
 
Если $R(\Theta ) <r$, то, поскольку
$ {\rm dist} ({\bf z}_\nu, {\cal L}(\theta))\to 0, \nu\to+\infty$, получается, что $\Theta$
лежит в некотором вполне рациональном подпространстве размерности $<r$, что невозможно.

Итак, $R(\Theta) = r$, и  утверждение  (ii) доказано.

{\bf 4.}   Для доказательства утверждения  (iii)  надо воспользоваться неравенством (\ref{minbody}),
из которого при  $ m>n$  получаем
\begin{equation}\label{dmdmd}
\zeta_\nu M_{\nu+1}\to 0,\,\,\,\,\nu \to +\infty.
\end{equation}
Предположив, что  $R(\Theta ) = 2$, получаем, что  в двумерном подпространстве  $\pi_\Theta={\rm span}\,\Lambda_\Theta$ 
(здесь   $\Lambda_\Theta$  есть именно та решетка, которая фигурирует в  определении величины  $R(\Theta )$)
находится одномерное подпространство $\pi_\Theta= \cap {\cal L}(\Theta )$,
не содержащее ненулевых точек решетки   $\Lambda_\Theta $ (ибо $\hbox{\got H}_\Theta=\{{\bf 0}\}$).
Теперь  (\ref{dmdmd}) противоречит предложению 2 , которое будет доказано в пункте \ref{dr}
(норма $|\cdot |_\bullet$  будет индуцирована $\sup$-нормой в $\mathbb{R}^d$).

{\bf 5.}  
Для доказательства утверждения  (iv)  надо снова воспользоваться неравенством (\ref{minbody}),
из которого при  $ m\ge n$  будет следовать, что
$$
\zeta_\nu M_{\nu+1} \le 1,\,\,\,\,\forall \, \nu
.
$$
Далее, рассмотрим вполне рациональное подпространство  $\pi_\Theta$
 и одномерное подпространство  $\ell_\Theta = \pi_\Theta \cap {\cal L}(\Theta)$.
Подпространство  $\ell_\Theta$ не лежит ни в каком собственном вполне рациональном подпространстве подпространства  $\pi_\Theta$.
Применяя предложение 3 из пункта  \ref{nvnp},
получаем
$$
\zeta_\nu M_{\nu+1} \to +\infty ,\,\,\,\, \nu
\to+\infty,
$$
что есть противоречие.

 Хочется заметить, что 
 в пункте {\bf 4}  выше мы фактически доказали несколько более общее утверждение:
в дополнение к пункту (iii) теоремы  \ref{NT1}  мы формулиируем следующий результат.
\begin{theorem}
\label{ngem}

Пусть  $2\le m\le n$.
 Пусть  правильная  матрица  $\Theta$  является  $\psi$-сингулярной с некоторой функцией  $\psi (t) = o(t^{-1}), t\to +\infty$.
 Тогда   $R(\Theta ) \ge 3$.

\end{theorem}

  В условиях теоремы \ref{ngem}  при  $m=n\ge 2$ требуется,
чтобы матрица  $\Theta$
была сингулярна в смысле оригинального определения А.Я.Хинчина,
при $n>m$  требуется нечто большее.

{\bf Замечание.}\,\,
 Во всех случаях когда  $R(\Theta) \ge 3$ найдется бесконечно много  значений индекса  $\nu$  таких, что
векторы наилучших приближений
$$
{\bf z}_{\nu -1},\,\,\, {\bf z}_\nu,\,\,\,
{\bf z}_{\nu+1}
$$
 линейно независимы.

Как мы увидим в следующем пункте,
подпространство $\pi_\Theta$ может иметь малую размерность  $R(\Theta)$. Тем не менее,
такое вырождение размерности накладывает иногда довольно жесткие условия на элементы матрицы  $\Theta$:

\begin{theorem}\label{qwqw}
Пусть матрица  $\Theta$   правильная.
Пусть
$
R(\Theta) \le  n +K(\Theta)-1
$.
Тогда  матрица  $\Theta$ состоит из чисел   $\theta^i_j$, связанных алгебраическим  соотношением степени $\le \min (m,R(\Theta) -K(\Theta) +1)$.
\end{theorem}

Теорема \ref{qwqw}  практически очевидна. Подпространство  ${\cal L}(\Theta)$ порождено векторами
$$
\underline{\theta}^1=\left(
\begin{array}{c}
-1\cr
0\cr
\vdots
\cr
0
\cr
\theta_1^1\cr
\vdots
\cr
\theta^1_n
\end{array}
\right),
\underline{\theta}^2=
\left(
\begin{array}{c}
0\cr
-1\cr
\vdots
\cr
0
\cr
\theta_1^2\cr
\vdots
\cr
\theta^2_n
\end{array}
\right),\cdots
,
\underline{\theta}^m=
\left(
\begin{array}{c}
0\cr
0\cr
\vdots
\cr
-1
\cr
\theta_1^m\cr
\vdots
\cr
\theta^m_n
\end{array}
\right).
$$
 В подпространстве
$\pi_\Theta$  имеется  базис из 
$R(\Theta)
$ 
независимых
целочисленых векторов  
$$
{\bf u}^j=\left(
\begin{array}{c}
u^j_1\cr
\vdots\cr
u^j_{m+n}
\end{array}
\right)
,\,\,\,\,
1\le j \le R(\Theta).
$$
Поскольку 
$$
{\rm dim }\,( {\rm span } ( 
\pi_\Theta\cup {\cal L}(\Theta)
))
= R(\Theta) + m - K(\Theta)
$$
и выполнено  (\ref{ones}), то векторы
$
\underline{\theta}^1,...,\underline{\theta}^m,
{\bf u}^1,...,{\bf u}^{R(\Theta)}$
линейно зависимы над $\mathbb{R}$.
Более того, поднабор 
\begin{equation}\label{podnabo}
\underline{\theta}^1,...,\underline{\theta}^m,
{\bf u}^1,...,{\bf u}^r,\,\,\,\,\ r={R(\Theta)-K(\Theta)+1},\,\,\,\, 1\le r\le n
 \end{equation}
тоже состоит из линейно зависимых над $\mathbb{R}$ векторов.
По условию  теоремы   $m+R(\Theta) -K(\Theta) +1\le m+n$.
Таким образом, у матрицы размера $(m+n)\times (m+R(\Theta) -K(\Theta) +1)$,  составленной из координат векторов (\ref{podnabo})
все миноры максимального порядка являются нулевыми.

Так как векторы ${\bf u}^1,...,{\bf u}^r$  независимы, то найдется набор индексов
$$
1\le i_1<...<i_r\le d$$
 такой, что
$$
\left|
\begin{array}{ccc}
u^1_{i_1}&\cdots &u^r_{i_1}\cr
\cdots&\cdots&\cdots\cr
u^1_{i_r}&\cdots &u^r_{i_r}
\end{array}
\right|\neq 0.
$$
  Набор, состоящий из $r$ индексов $
1\le i_1<...<i_r\le d$ дополним до набора, состоящего из $r+m$ индексов, вида
$$
1,2,...,m,j_1,...,j_r;\,\,\,\,\,
m<j_1<...<j_r \le d.
$$

Имеем
 $$
\left|
\begin{array}{cccccc}
-1&\cdots &0&u^1_1&\cdots& u_1^r
\cr
\cdots&\cdots&\cdots&\cdots&\cdots&\cdots \cr
0&\cdots&-1&u^1_m&\cdots&u^r_m\cr
\theta^1_{j_1}&\cdots &\theta^m_{j_1}& u^1_{{j_1}}&\cdots & u^r_{{j_1}}\cr
\cdots&\cdots&\cdots&\cdots&\cdots&\cdots \cr
\theta^1_{j_r}&\cdots &\theta^m_{j_r}& u^1_{{j_r}}&\cdots & u^r_{{j_r}}
\end{array}
\right|
=0.
$$
Ясно, что из столбцов  $u^j$  этого определителя выбирается ненулевой минор порядка  $r$,
причем найдется таковой ненулевой минор, отличный от правого нижнего.
Таким образом,  равенство нулю определителя задает на элементы матрицы  $\Theta$
нетривиальное алгебраическое соотношение степени  $\le  \min (m,r)$.

Из пункта  (i) теоремы \ref{NT1} и второго из равенств (\ref{razmeri}) сразу
получаем

{\bf  Следствие 1.}\,\,{\it
Пусть матрица  $\Theta$   правильна.
Тогда
если  $K(\Theta) = m$  то $\pi_\Theta = \hbox{\got R}_\Theta$ и
 $R(\Theta) = {\rm DIM}_\mathbb{Z}\, ^t\!\Theta$. }

Из теорем \ref{NT1},\ref{qwqw} и неравенства (\ref{ones})
получаем

{\bf  Следствие 2.}\,\,{\it
Пусть матрица  $\Theta$   правильна.
Тогда

{\rm (i)} если элементы матрицы $\Theta$  алгебраически независимы в совокупности, то
$R(\Theta )\ge n+1$;

{\rm (ii)}если $(m,n)\neq (1,1)$  и элементы матрицы $\Theta$  алгебраически независимы в совокупности, то
$R(\Theta ) \ge  3 $.
}

Для  случая  $m=n$   из теоремы \ref{NT1}  (пункт (iv))  и теоремы \ref{qwqw}  имеем

{\bf Следствие 3.}
{\it

{\rm (i)}
Если для  правильной матрицы  $\Theta$ известно, что
$m=n>1$  и $ R(\Theta)\le n+1$, то
 элементы матрицы  $\Theta$  алгебраически зависимы.

{\rm (ii)}
Если для правильной матрицы  $\Theta$ известно, что
$m=n>1$  и $ K(\Theta)=1$, то
$R(\Theta) = 2$  и
 элементы матрицы  $\Theta$  алгебраически зависимы.}

{\bf Следствие 4.}
{\it  Если  $m=n=2$  и матрица  $\Theta$  правильна, то

{\rm (i)}  если   $K(\Theta) = 1$  то  $R(\Theta) = 2$;

{\rm (ii)}  если   $K(\Theta) = 2$  то 
$\pi_\Theta \supseteq {\cal L}(\Theta) $  и 
 $R(\Theta) = {\rm DIM}_\mathbb{Z}\, ^t\!\Theta$.

В частности, если ${\rm DIM}_\mathbb{Z}\, ^t\!\Theta=4$,
 то либо  $R(\Theta) = 2$,  либо  $R(\Theta) = 4$.}

 \subsection{Вырождение размерности, $R(\Theta ) = 2$}\label{vr2}

Пусть  $\xi \in (0,1)$  есть иррациональное число
и $p_\nu/q_\nu, \nu = 1,2,3,... $ -  суть все подходящие дроби к числу  $\xi$. 
Пусть $m\ge 2,n\ge 3$. Предположим, что элементы
$\theta^i_j , 1\le i\le 2, 1\le j \le n$  связаны соотношением
\begin{equation}\label{zavi}
\theta^2_j = -\xi \theta^1_j, \,\,\,\, 1\le j \le n
\end{equation}
 и рассмотрим матрицу
\begin{equation}\label{zavi1}
 \Theta = 
 \left(
\begin{array}{ccccc}
\theta^1_1& \theta_1^2&\theta_1^3
&\cdots &\theta_1^m
\cr
\cdots &\cdots &\cdots &\cdots&\cdots \cr
 \theta_n^1 & \theta^2_n&\theta_3^n&\cdots &\theta^m_n
\end{array}
\right)=
\left(
\begin{array}{ccccc}
\theta^1_1& -\xi\theta_1^1&\theta_1^3
&\cdots &\theta_1^m\cr
\cdots &\cdots &\cdots &\cdots&\cdots \cr
\theta_n^1 & -\xi \theta_n^1 & \theta^3_n&\cdots &\theta^m_n
\end{array}
\right)
\end{equation}
(если  $m=2$  в этой матрице два столбца).
Отметим, что все элементы матрицы (\ref{zavi1})   при надлежащем выборе величин
$\theta^1_j, j =1,2,3$ могут быть линейно независимы вместе с единицей над  $\mathbb{Z}$
и тогда
$$
{\rm DIM}_\mathbb{Z} \Theta =
{\rm DIM}_\mathbb{Z}\,^t\!\Theta = m+n,
$$
и матрица правильна.
В то же время условия  (\ref{zavi})  показывают, что элементы матрицы (\ref{zavi1})
 всегда представляют из себя набор алгебраически зависимых чисел (так что, теорема \ref{yy1}, которую мы формулируем ниже, согласуется с теоремой \ref{qwqw}).

\begin{theorem}\label{yy1} Пусть  $m <n$.  Тогда для
 почти всех наборов из  $n(m-1)$ чисел  $(\theta^i_1,\theta^i_2,...,\theta^i_n)$, $ i =1,3,4,...,m$  таких, что
$1/3<\theta_1^1,\theta_2^1,...,\theta_{n-1}^{1}<2/3<\theta^1_n<1$,
$\theta^i_j \in(0,1),\, 3\le i\le m,\, 1\le j \le n$
 последовательность векторов наилучших приближений для матрицы
(\ref{zavi1})  отличается от последовательности векторов
$$
(p_\nu,q_\nu,\underbrace{0,...,0}_{m+n-2}),\,\,\,\,\, \nu = 1,2,3,...
$$
 не более чем на конечное число точек, и, следовательно
 $R(\Theta) = 2$.

\end{theorem}

{\bf Следствие.}
\,\,{\it  При  $m\ge 2,n>m$  существуют матрицы  $\Theta$, состоящие из элементов линейно независимых вместе с единицей над  
$\mathbb{Z}$  и такие, что   $R(\Theta ) = 2$.
}

Доказательство теоремы  \ref{yy1}.

 Видим, что
$$
\max_{1\le j\le n}
||p_\nu \theta^1_j+q_\nu \theta^2_j||=
\max_{1\le j\le n}
||(p_\nu  -q_\nu\xi) \theta^1_j|| =
 |(p_\nu  -q_\nu\xi) \theta^1_n| 
 .
$$
С другой стороны
из свойств непрерывных дробей (см. пункт \ref{cd}, формулы (\ref{approx},\ref{qt}))
 получаем
\begin{equation}\label{dobava}
\max_{1\le j\le n}
||p_\nu \theta^1_j+q_\nu \theta^2_j||=
  |(p_\nu  -q_\nu\xi) \theta^1_n| 
<\frac{1}{q_{\nu+1}}.
\end{equation}

Чтобы вектор  $(p_\nu,q_\nu )$ 
 был бы вектором наилучшего приближения для матрицы  (\ref{zavi1}),
а вектор
$(p_{\nu+1},q_{\nu+1} )$ 
 был бы  в точности следующим вектором наилучшего приближения для матрицы  (\ref{zavi1})
достаточно
чтобы
\begin{equation}\label{maxx}
\min_{(x_1,x_2,...,x_m;y_1,y_2,...,y_n)}
\max_{1\le j \le n}
|(x_1-\xi x_2)\theta^1_j + x_3\theta^3_j + \cdots +x_m\theta^m_j +y_j|
>
|(p_\nu  -q_\nu\xi) \theta^1_n| ,
\end{equation}
где минимум берется по  целым точкам
$(x_1,...,x_m;y_1,...,y_n)$ таким, что
\begin{equation}\label{ppo}
0\neq
\max_{1\le i \le m} |x_i| \le q_{\nu+1},\,\,\,\,
\max_{1\le j\le n} |y_j| \le 
  1+  |x_1-\xi x_2|+ |x_3|+...+|x_m|.
\end{equation}
и
$$
(x_1,...,x_m;y_1,...,y_n) \neq \pm
(p_\nu, q_\nu,\underbrace{0,...,0}_{m+n-2}),\pm  (p_{\nu +1}, q_{\nu+1} ,\underbrace{0,...,0}_{m+n-2}).
$$

Если  $x_3=\cdots = x_m =y_1=\cdots =y_n =0$,  то от свойств чисел $\theta_j^i$  инчего не зависит и
$$
\max_{1\le j \le n}
|(x_1-\xi x_2)\theta^1_j + x_3\theta^3_j + \cdots +x_m\theta^m_j +y_j|
=
|
(x_1-\xi x_2)\theta^1_n |
>|(p_\nu  -q_\nu\xi) \theta^1_n|,
$$
 поскольку  $p_\nu < q_\nu$  и подходящие дроби задают наилучшие одномерные приближения
(см. пункт \ref{cd}).

Условие (\ref{maxx})  будет выполнено, если для каждого рассматриваемого целочисленного набора
$(x_1,...,x_m,y_1,...,y_n)$ 
с условием (\ref{ppo}) и такого что
\begin{equation}\label{oop}
|x_3|+...+|x_m|+|y_1|+...+|y_n|\neq 0
\end{equation}
при каждом $j, 1\le j\le n$ 
будем иметь
$$
 (\theta_j^1,\theta_j^3,...,\theta_j^m)\not\in
{\cal J}_\nu(x_1,x_2,...x_m;y_j),
$$
где 
$$
{\cal J}_\nu(x_1,x_2,...x_m;y_j)=\{(\theta_j^1,\theta_j^3,...,\theta_j^m)\in [0,1]^{m-1}:\,
|(x_1-\xi x_2)\theta^1_j + x_3\theta^3_j + \cdots +x_m\theta^m_j +y_j|\le 1/q_{\nu+1}\}.
$$
Заметим, что для меры Лебега последнего множества выполнено
\begin{equation}\label{mere}
\mu ({\cal J}_\nu(x_1,x_2,...x_m;y_j))\le 
\frac{m^{3/2}}{q_{\nu+1}(|x_1-\xi x_2|+ |x_3|+...+|x_m|)}.
\end{equation}
Кроме того, если  $x_3=...=x_m=0$  то из (\ref{oop})   получаем, что некоторое  $y_j\neq 0$.
Тогда  если 
$
|x_1-\xi x_2|< \frac{1}{2}
$,
то (\ref{maxx}) выполнено автоматически, и, стало быть, можно считать, что
$$
|x_1-\xi x_2|\ge \frac{1}{2}.
$$
Итак, для выполнения условия (\ref{maxx})  достаточно чтобы
$$
\left(
\begin{array}{cccc}
\theta_1^1&\theta_1^3&\cdots&\theta_1^m
\cr
\cdots&\cdots&\cdots&\cdots
\cr
\theta_n^1&\theta_n^3&\cdots&\theta_n^m
\end{array}
\right)
\not\in {\cal J}_\nu =
J_\nu^{(1)}\cup J_\nu^{(2)},
$$
где
$$
J_\nu^{(1)}
=
\bigcup_{
\begin{array}{c}
(x_1,x_2):
\cr
0<\max\{|x_1|,|x_2|\}\le q_{\nu+1}
\cr
|x_1-\xi x_2|\ge 1/2
\end{array}
}
\bigcup_{
\begin{array}{c}
(y_1,...,y_n):\cr
\max_{1\le j \le n}
|y_j|\le 4|x_1-\xi x_2|
\end{array}
}
$$
$$
J_\nu(x_1,x_2,\underbrace{0,...,0}_{m-2};y_1)\times\cdots\times
J_\nu(x_1,x_2,\underbrace{0,...,0}_{m-2};y_n),
$$

$$
J_\nu^{(2)}
=
\bigcup_{
\begin{array}{c}
(x_1,x_2):
\cr
0<\max\{|x_1|,|x_2|\}\le q_{\nu+1}
 \end{array}
}
\bigcup_{
\begin{array}{c}
(x_3,...,x_m):
\cr
0<\max_{3\le i\le m}|x_i|\le q_{\nu+1}
 \end{array}
}
\bigcup_{
\begin{array}{c}
(y_1,...,y_n):\cr
\max_{1\le j \le n}
|y_j|\le 2(|x_1-\xi x_2|+\sum_{i=3}^m|x_i|)
\end{array}
}
$$
$$
J_\nu(x_1,x_2,x_3,...,x_m;y_1)\times\cdots\times
J_\nu(x_1,x_2,x_3,...,x_m;y_n).
$$
Из (\ref{mere})  получаем, что
$$
\mu (J_\nu^{(1)})\ll q_{\nu+1}^{2-n},\,\,\,
\mu (J_\nu^{(2)})\ll q_{\nu+1}^{m-n},\,\,\,
\mu ({\cal J}_\nu) \ll q_{\nu+1}^{m-n}\le q_{\nu+1}^{-1}.
$$
 
 Ряд  $\sum_{\nu=1}^\infty \frac{1}{q_\nu}$
  сходится и теорема \ref{yy1}, таким образом, вытекает из леммы Бореля-Кантелли.

В заключении этого пункта отметим, что при $m=n=2$ пример матрицы $\Theta$, состоящей из чисел, линейно независимых вместе с единицей над  $\mathbb{Z}$  и такой, что
 $R(\Theta) =2$, автору неизвестен.

\subsection{Вырождение размерности, $R(\Theta ) = 3$}\label{vr3}

Приведенные  в предыдущем пункте результаты  показывают, насколько точна
  нижняя оценка пункта (i)  теоремы \ref{NT1}   при $m\ge 2$.
  Ниже
мы установим, что нижняя оценка из утверждения (iii) теоремы  \ref{NT1} тоже точна.
Для этого мы будем использовать сингулярные матрицы.

Сначала рассмотрим весьма специальный случай.

Пусть числа  $\xi^1 , \xi^2 \in (0,1/2)$  линейно независимы вместе с единицей над  $\mathbb{Z}$.
Пусть  $(x_{1,\nu},x_{2,\nu},x_{3,\nu})$  есть "расширенный" вектор наилучших приближений для набора
$(\xi^1,\xi^2)$
(тут мы используем нестандартные для настоящей статьи обозначения, но это будет удобно).
Тогда
$$
||x_{1,\nu}\xi^1+x_{2,\nu}\xi^2||=
|x_{1,\nu}\xi^1+x_{2,\nu}\xi^2+x_{3,\nu}|< M_{\nu+1}^{-2},
$$
где можно считать, что
$$
M_\nu  = \max_{i=1,2,3} |x_i|
.
$$

Положим
$$
\theta^1_j = \xi^1\theta^3_j,\,\,\,
\theta^2_j = \xi^2\theta^3_j,\,\,\,
1\le j \le n
$$
 и рассмотрим матрицу
\begin{equation}\label{tritri}
\left(
\begin{array}{ccc}
\theta^1_1&\theta^2_1&\theta^3_1
\cr
\cdots&\cdots&\cdots\cr
\theta^1_n&\theta^2_n&\theta^3_n
\end{array}
\right)
=
\left(
\begin{array}{ccc}
\xi^1\theta^3_1&\xi^2\theta^3_1&\theta^3_1
\cr
\cdots&\cdots&\cdots\cr
\xi^1\theta^3_n&\xi^2\theta^3_n&\theta^3_n
\end{array}
\right)
\end{equation}
\begin{theorem}\label{nnn}
 Пусть ряд
\begin{equation}\label{sese}
\sum_{\nu=1}^\infty M_{\nu+1} ^3||x_{1,\nu}\xi^1+x_{2,\nu}\xi^2||^n
\end{equation}
 сходится.
Тогда для почти всех наборов  $(\theta^3_1,...,\theta^3_n) \in \mathbb{R}^n$
таких, что  $ 0<\theta^3_2,...,\theta^3_n <1/2<\theta^3_1<1$
 последовательность наилучших приближений для матрицы
$\Theta$, определенной в  (\ref{tritri}) отличается от последовательности 
$$
(x_{1,\nu},x_{2,\nu},x_{3,\nu},\underbrace{0,...,0}_{ n})
$$
 не более чем на конечное число точек.
\end{theorem}

Доказательство теоремы \ref{nnn}
подобно доказательству теоремы \ref{yy1}:
 для почти всех  $(\theta^3_1,...,\theta^3_n) \in \mathbb{R}^n$
 надо добиться того, чтобы
$$
\min_{x_1,x_2,x_3;y_1,...,y_n}\,\,
\max_{1\le j\le n} |\theta_3^j(x_1\xi^1+x_2\xi^2+x_3)+ y_j|>
|x_{1,\nu}\xi^1+x_{2,\nu}\xi^2+x_{3,\nu }|,
$$
где минимум берется по целочисленным векторам, удовлетворяющим соотношениям
$$
0\neq\max_{i=1,2,3} |x_i|\le M_{\nu+1},
$$
$$
(x_1,x_2,x_3,y_1,...,y_n)\neq
\pm
(x_{1,\nu},x_{2,\nu},x_{3,\nu},\underbrace{0,...,0}_{n}),
\pm
(x_{1,\nu+1},x_{2,\nu+1},x_{3,\nu+1},\underbrace{0,...,0}_{n})
.
$$
Естественно, следует различать случаи ${\bf y}={\bf 0}$   и ${\bf y}\neq {\bf 0}$.

Из сходимости ряда  $\sum_{\nu}M_\nu
$
(эта сходимось обеспечивается экспоненциальным ростом величины  $M_\nu$, см., например, лемму
из работы  \cite{sado}) получаем

{\bf Следствие.}\,\,\,{\it
 Если  $n \ge 2$, то для почти всех наборов 
$(\theta^3_1,...,\theta^3_n) \in \mathbb{R}^n$
для матрицы (\ref{tritri})  выполнено  $ R(\Theta) = 3$.

}

Теперь мы сформулируем и докажем утверждение, которое  обобщает теорему Н.Г.Мощевитина и О.Н.Германа \cite{sado},
 которая, в свою очередь, обобщает результата из   автора \cite{MDAN}.
 Отметим, что в работе \cite{sado} авторы анонсировали  еще несколько результатов, в частности, результаты о существовании сингулярных систем специального вида  (при  $n=1$),
подробное доказательство которых пока не опубликовано. Формулировки мы приведем в пункте
\ref{specivid}.

 Пусть  $ m^*> m$.
 Наряду с матрицей  $\Theta$ вида (\ref{1f})  рассмотрим "расширенную" матрицу
$$
\Theta^* =
\left(
\begin{array}{cccccc}
\theta_{1}^1&\cdots&\theta_{1}^m&\theta_1^{m+1}&\cdots &\theta_1^{m^*}\cr
\cdots&\cdots&\cdots&\cdots&\cdots&\cdots\cr
\theta_{n}^1&\cdots&\theta_{n}^m&\theta_n^{m+1}&\cdots &\theta_n^{m^*} 
\end{array}
\right)
$$
 с вещественными элементами. Таким образом, мы  увеличили матрицу
$\Theta $, добавив  $n(m^*-m) $  элементов.
Набор добавленных элементов можно отождествить с некоторым вектором
$\underline{\Theta} =(\theta_1^{m+1},...,\theta_{n}^{m^*})$
 в пространстве
$\mathbb{R}^{n(m^*-m) }$. 
Справедливо следующее утверждение,
непосрдственно обобщающее результат Н.Г.Мощевитина и О.Н.Германа  из \cite{sado} и
 и развивающее подход автора  \cite{MDAN},\cite{ME}.

\begin{theorem} \label{t:zeta_convergence_implies_degeneracy}
Пусть точки
$$
{\bf z}_\nu
=
(x_{1,\nu},...,x_{m,\nu}, y_{1,\nu},...,y_{n.\nu})
 \in \mathbb{Z}^{m+n}
,\quad\nu=1,2,3,\ldots 
$$
образуют последовательность наилучших приближениий для матрицы  $\Theta$.
 Пусть величины  $M_\nu$  и $\zeta_\nu$  суть элементы последовательностей  (\ref{1},\ref{2}).
Рассмотрим последовательность целых точек
\begin{equation} \label{sequ}
{\bf z}_\nu^*=
(
x_{1,\nu},....,x_{m,\nu}, \underbrace{0,\dots,0}_{m^*-m} ,
y_{1,\nu},...,y_{n,\nu})\in \mathbb{Z}^{m^*+n}.
\end{equation}
Предположим, что ряд
\begin{equation} \label{er}
\sum_{\nu=1}^\infty M_{\nu+1}^{\max (m+n,m^*)}(\log M_{\nu+1})^{{\delta (m^*,m+n)}}\zeta_\nu^n,
\quad\delta (a,b)=
\begin{cases}
1,\quad  a=b \cr
0,\quad a\neq b
\end{cases},
\end{equation} сходится. Тогда для почти всех
"добавляемых" векторов 
$\underline{\Theta}\in\mathbb{R}^{n(m^*-m)}$ (в смысле меры Лебега) последовательность наилучших приближений
для матрицы $\Theta^*$ отличается от последовательности
${\bf z}_\nu^*$ на не более, чем конечное число точек.
\end{theorem}

Из теоремы \ref{t:zeta_convergence_implies_degeneracy} в свете предложения 1 непосредственно  получаем

{\bf Следствие.}

{\it
Пусть $m\geq 2$ и пусть  матрица $\Theta$  является   $\psi$-сингулярной. Пусть
также ряд
\begin{equation} \label{er1}
 \sum_{\nu=1}^{\infty}M_\nu^{\max (m+n,m^*)}(\log M_\nu)^{{ \delta (m^*,m+n )}}(\psi(M_\nu))^n
\end{equation}
сходится. Тогда для почти всех 
"добавляемых" векторов $\underline{\Theta}\in\mathbb{R}^{n(m^*-m) }$
последовательность всех наилучших
приближений для матрицы $\Theta^*$ отличается от
последовательности (\ref{sequ}) на не более, чем конечное число точек.
}

Естественно, если положить  $n=1, m=2, \, \psi (t) = t^{-m^*-2}$  и взять числа
 $\theta_j^1,\theta_j^2$,
 существование которых обеспечивается теоремой \ref{T1},  то  наше следcтвие  превращается в
следующий результат, имевшийся в \cite{MDAN},\cite{ME}.

\begin{theorem} \label{dime2}
 
 При $n=1$  и каждом  $m\ge 2$
 существует  набор  $\Theta$,  состоящий из
 алгебраически независимых
чисел  $\theta^i, \,1\le i \le m$,
такой, что  вектора наилучших приближений  ${\bf z}_\nu$
 при всех достаточно больших значениях  $\nu$  лежат в некотором
трехмерном подпространстве, то есть
 $R(\Theta) = 3$.

\end{theorem}

При произвольном $n$  и  $m\ge 3$  из теорем  \ref{nnn}  и \ref{t:zeta_convergence_implies_degeneracy} получаем следующий результат.

\begin{theorem} \label{dime22}
 
 При $m\ge 3$  каждом  $n$
 существует  матрица  $\Theta$,  состоящая из
 линейно  независимых вместе с единицей над $\mathbb{Z}$ 
чисел  $\theta^i_j$,
такая, что  вектора наилучших приближений  ${\bf z}_\nu$
 при всех достаточно больших значениях  $\nu$  лежат в некотором
трехмерном подпространстве, то есть
 $R(\Theta) = 3$.

\end{theorem}

{Доказательство теоремы \ref{t:zeta_convergence_implies_degeneracy}.}

Можно считать,  что  $ \theta_j^i\in [0,1]$  при всех значениях  $i,j$.

Нам надо доказать, что
 при фиксированной матрице $\Theta$,  удовлетовряющей условию теоремы, 
   для почти каждого "дополняющего"  набора  $\underline{\Theta}$
 существует
такое $\nu_0$, что для всех $\nu\geq \nu_0$
будет выполняться
\begin{equation} \label{eq:min_geq_zeta}
  \min\,\,
\max_{1\le j \le n}\left|
\sum_{1\le i \le m^*} \theta_j^i + y_j \right|>\zeta_\nu,
\end{equation}
где минимум берется по всем целым точкам 
$$
{\bf z}
=
(x_{1},...,x_{m^*}, y_{1,\nu},...,y_{n})
 \in \mathbb{Z}^{m^*+n}\setminus \{{\bf 0}\}
$$
таким, что
$$
\max_{1\leq j\le m^*}|x_j|\leq M_{\nu+1},\quad 
{\bf z} \neq {\bf z}_{\nu}^*.
$$

Условие (\ref{eq:min_geq_zeta})  будет выполнено, если 
 для каждого вектора  
$$(x_{1},...,x_{m^*}, y_1,...,y_{n})\in \mathbb{Z}^{m^*+n}\setminus \{{\bf 0}, \pm {\bf z}_\nu, \pm {\bf z}_{\nu+1}\}
$$
такого, что
$\max_{1\le i\le m^*}|x_i|\le M_{\nu+1}
$ имеет место следующее:
 для некоторого  $j$  из промежутка
$1\le j \le n$
выполняется 
 
\begin{equation} \label{momo}
  x_{m+1} \theta_j^{m+1}+ +\cdots+x_{m^*}\theta^{m^*}_j  \not\in  J_\nu( y_j, x_1,...,x_m),
\end{equation}
где 
$$
 J_\nu(y_j,x_1,...,x_m) = 
(-y_j-x_1\theta_j^1+\cdots + x_m\theta_j^m -\zeta_\nu,
-y_j-x_1\theta_j^1+\cdots + x_m\theta_j^m +\zeta_\nu
).$$

 Заметим, что условие (\ref{momo}) означает, что расстояние
от точки $(\theta_j^{m+1},...,\theta_j^{m^*} )\in [0,1]^{m^*-m}$ до подпространства
$$ 
\Big\{ (u_{m+1},\ldots,u_{m^*})\in\mathbb{R}^{m^*-m} \,:\
   x_{m+1}u_{m+1}+\cdots+x_{m^*}u_{m^*}=-y_j-x_1\theta_j^1-\cdots-x_m\theta_j^m \Big\} $$
не меньше, чем $\zeta_\nu\cdot (x_{m+1}^2+\cdots+x_{m^*}^2)^{-1/2}$.

 Положим
$$
\Omega_\nu (
{\bf z} ) =
\Omega_\nu ({\bf x},{\bf y}) =
\Big\{ \underline{\bf \Theta}\in[0,1]^{n(m^*-m) } \, :\
$$
$$
  x_{m+1} \theta_j^{m+1}+ \cdots+x_{m^*}\theta^{m^*}_j  \not\in  J_\nu( y_j, x_1,...,x_m),
\,\,\,\, 1\le j \le n,
 \Big\}
$$
и
$$
\Omega_\nu =\bigcup_{\bf y}\bigcup_{\bf x} \Omega_\nu ({\bf x},{\bf y})
,$$
причем в последней формуле объединения берутся по областям целых точек
$$
\{{\bf y}\, :\,\,
\max_{1\le j \le n } |y_j|\le (m^*+1)M_{\nu+1}\}
,\,\,\,
\{{\bf x}\,:\,\,
0<\max_{1\le i\le m^*}|x_i|\le M_{\nu+1}
\}.$$
Согласно лемме Бореля-Кантелли, для доказательства теоремы достаточно убедиться в том, что
\begin{equation}
\sum_{\nu\ge \nu_0}
\mu \left(  \Omega_\nu \right) \to 0,\,\,\, \nu_0 \to +\infty.
\label{ene}
\end{equation}
 Но
$$
\mu \left(  \Omega_\nu \right)
\ll\zeta_\nu^{n }M_{\nu+1}^{n }\,\sum_{ x_1,...,x_m}
\,\sum_{x_{m+1},...,x_{m^*}}
 \frac{1}{(\max_{ m\le i \le m^*}|x_i|)^n}
\ll
$$
$$
\ll \zeta_\nu^{n}M_{\nu+1}^{n+m }\,
\sum_{1\le t\le M_{\nu+1}}\,
t^{m^*-m-n-1}\ll \zeta_\nu^n
M_{\nu+1}^{\max (m+n,m^*)}(\log M_{\nu+1})^{ \delta (m^*,m+n)}.
$$
Теперь соотношение (\ref{ene}) следует из сходимости  ряда (\ref{er}), и теорема доказана.

\section{Одномерные диофантовы приближения.}\label{odp}

Обсудим простейшую ситуацию, когда $m=n=1$. В этом случае мы имеем дело с задачей о приближении одного числа $\alpha =\theta_{1}^1$
рациональными дробями.
\subsection{Цепные дроби}\label{cd}

 Хорошо известно, что  вопрос о наилучших приближениях в этом случае полностью решается посредством аппарата цепных дробей
(см. \cite{HZT}). Напомним, что если вещественное $\alpha$ представлено в виде цепной дроби
$$
\alpha= [a_0;a_1,a_2,...,a_t,...]= a_0+\frac{1}{\displaystyle{a_1+\frac{1}{\displaystyle{a_2 + \dots+\frac{1}{\displaystyle{a_t +
\displaystyle{\dots} }}}}}},\,\,\,\,\, a_0 \in \mathbb{Z},\,\,\,\,\, a_t \in \mathbb{N},\,\,\,\,\, t =1,2,3,...
 $$
(конечной или бесконечной, в зависимости от того, является ли $\alpha$ рациональным числом или нет), то подходящими дробями к $\alpha$
называются рациональные дроби вида
$$
\frac{p_\nu}{q_\nu} = [a_0;a_1,a_2,...,a_\nu]= a_0+\frac{1}{\displaystyle{a_1+\frac{1}{\displaystyle{a_2 + \dots+\frac{1}{\displaystyle{a_\nu}}}}}}.
 $$
 Наилучшие приближения, определенные выше в пункте \ref{ssinp} (точнее "расширенные" вектора ${\bf z}_\nu\in \mathbb{Z}^2$)
  будут (в терминологии книги \cite{HZT}) наилучшими приближениями {\it второго рода}.
Для них справедливо следующее утверждение (см. \cite{HZT}, теорема 16).

\begin{theorem}\label{T6}
Всякое наилучшее приближение ${\bf z}_\nu$ имеет вид $ {\bf z}_\nu=(q_\nu, p_\nu)$, где $q_\nu$ и $p_\nu$
  являются знаменателем и числителем некоторой подходящей к $\alpha$ дроби. 
\end{theorem}

Отметим, что если $\alpha \neq \frac{p_\nu}{q_\nu}$ (в частности, если $\alpha$ есть число иррациональное), то для его приближения подходящей дробью
имеют место неравенства (cм. \cite{HZT}, теоремы 9,13)
\begin{equation}
\frac{1}{q_\nu(q_\nu+q_{\nu+1})}< \left| \alpha -\frac{p_\nu}{q_\nu} \right| < \frac{1}{q_\nu q_{\nu+1}} ,
\label{approx}
\end{equation}
в частности,
\begin{equation}
||q_\nu \alpha || > \frac{1}{2q_{\nu+1}},
\label{qt}
\end{equation}
или, на языке обозначений пункта \ref{ssinp}  (и при $\alpha \in (0,1)$)
\begin{equation}
\zeta_\nu \ge \frac{1}{2M_{\nu+1}}, \,\,\,\ \forall \nu \in \mathbb{N}. \label{dim1}
\end{equation}
Пользуясь случаем, отметим, что для разности из (\ref{approx}) имеется простое и изящное равенство
\begin{equation}
\left| \alpha -\frac{p_\nu}{q_\nu} \right| =
\frac{1}{q_\nu^2 (\alpha_{\nu+1}+\alpha^*_\nu)},
\label{approx1}
\end{equation}
где
$$
\alpha_{\nu+1} =[a_{\nu+1};a_{\nu+2},a_{\nu+3},...],\,\,\,\,\,
\alpha^*_\nu =[0;a_\nu,a_{\nu-1},a_{\nu-2},...,a_1].
 $$
Это равенство позволяет провести весьма детальное исследование одномерных диофантовых приближений,
например, свести основные задачи изучения  {\it спектра Лагранжа }
$$
\mathbb{L}=
\{\lambda \in \mathbb{R}:\,\, \exists \alpha\in \mathbb{R}\,\text{такое что}\,
\lambda =(\liminf_{q \to +\infty } q||q \alpha|| )^{-1}\}
$$
к задачам о бесконечных последовательностях
\cite{Cusick}. 
О спектре Лагранжа  $\mathbb{L}$  имеется много работ --
см., например, библиографию из книги  \cite{Cusick} и замечательный обзор А.В.Малышева \cite{Malyshev}.
Никакой многомерный аналог равенства (\ref{approx1}), по-видимому, неизвестен.

Если теперь предположить, что $\alpha$ иррационально и представляет из себя есть сингулярную систему с $ m=n=1$ то при все достаточно больших
$t$  согласно предложению 1  будет выполняться $\zeta_\nu < \frac{1}{2M_{\nu+1}}$, что противоречит (\ref{dim1}). Значит иррациональных $\theta$,
являющихся сингулярными системами не бывает.

 Равенство (\ref{approx1})  естественно переписать в виде
$$
q_\nu ||q_\nu \alpha || = 
\frac{1}{\alpha_{\nu+1}+\alpha^*_\nu}.
$$
Следует сказать, что имеет место (см., например \cite{DSdir0})  похожее равенство
\begin{equation}
\label{approx11}
q_{\nu+1} ||q_\nu \alpha || = 
\frac{1}{1+\displaystyle{\frac{1}{\alpha_{\nu+2}\cdot\alpha_{\nu+1}^{**}}}},
\end{equation}
где
$$
\alpha_{\nu+1}^{**}=
a_{\nu+1} + \alpha^*_\nu=
[a_{\nu+1};a_\nu,...,a_1].
$$
 
Это равенство, в частности,  позволяет исследовать
{\it спектр Дирихле}
$$
\mathbb{D}=
\{\lambda \in \mathbb{R}:\,\, \exists \alpha\in \mathbb{R}\,\text{такое что}\,
\lambda =(\limsup_{q \to +\infty } q||q \alpha|| )^{-1}\}
.$$
 В отличие от спектра Лагранжа, изучению спектра Дирихле посвящено
мало работ (см. библиографию в статье \cite{john}).

  Естественно, определения спектров Лагранжа и Дирихле можно давать в терминах функции  $\psi_\alpha (t)$.

\subsection{Функция  $\psi_\alpha(t)$}

Для случая  $m=n=1$ при $\theta^1_1 =\alpha$ функция (\ref{fj}) 
принимает вид
$$
\psi_\alpha (t) = \min_{1\le x\le t} ||\alpha x||.
$$

Недавно И.Д.Кан и Н.Г.Мощевитин \cite{KAMOarxiv} получили следующий результат.

\begin{theorem}
\label{alphabet}
Для произвольных двух иррациональных чисел 
  $\alpha , \beta$ таких, что
 $\alpha\pm \beta \not\in \mathbb{Z}$  функция разности 
$$
\psi_\alpha (t) - \psi_\beta (t)
$$
 бесконечно много раз меняет знак при 
 $t\to +\infty$.
\end{theorem}

 Доказательство этой теоремы мы здесь не приводим.
Отметим лишь, что оно связано с применением формул
(\ref{approx1},\ref{approx11}).

Отметим также, что теоремы \ref{T1},\ref{T2} А.Я.Хинчина, процитированные нами в пункте \ref{tees},
показывают, что результат теоремы \ref{alphabet} не может быть перенесен на большие размерности.

 Возьмем в теореме \ref{T1}
  $\psi (t) = o(t^{-2}),\, t \to +\infty$.
 Пусть  $\theta^1,\theta^2$ - те числа, существование которых утверждает теорема \ref{T1}.
Возьмем другие числа  
 $(\beta^1,\beta^2)$  плохо приближаемыми (в смысле линейной формы):
$$
\inf_{(x_1,x_2) \in  \mathbb{Z}^2\setminus \{(0,0)\}} \left(
||x_1\beta^1+x_2\beta^2|| \cdot \max(|x_1|,|x_2|)^2 \right) >0.
$$
Видим, что при достаточно больших
  $t$ получается
$$
\psi_{(\theta^1,\theta^2)}
(t)<
\psi_{(\beta^1,\beta^2)}
(t).
$$

Ситуация для совместных приближений аналогична.
Надо взять в теореме \ref{T2}  функцию
   $\psi_1 (t) = o(t^{-1/2}),\, t \to +\infty$, числа
$\theta_1,\theta_2$  из теоремы  \ref{T2}, а также совместно плохо приближаемые числа
$\left(\begin{array}{c}\beta_1\cr \beta_2\end{array}\right)$:
$$
\inf_{x \in  \mathbb{Z}\setminus \{0\}}
\,\left(\max_{j=1,2}
||x\beta_j|| \cdot |x|^{1/2} \right) >0.
$$
Видим, что

$$
\psi_{
\left(
\begin{array}{c}
\theta_1\cr \theta_2
\end{array}\right)}
(t) <
\psi_{
\left(
\begin{array}{c}
\beta_1\cr \beta_2
\end{array}\right)}
(t)
$$
для всех достаточно больших $t$.

\subsection{Двумерные решетки}
\label{dr}

  Здесь мы хотим обобщить неравенство  (\ref{approx}) (а, значит, и оценку снизу (\ref{dim1}))  удобным для приложений образом.

Рассмотрим двумерное вполне рациональное линейное подпространство  $\pi\subset \mathbb{R}^d$
и двумерную решетку  $\Lambda =\pi \cap \mathbb{Z}^d$, в нем содержащуюся. 
Будем считать, что множество
$$
{\cal B}=
\pi \cap \{{\bf z}=({\bf x},{\bf y}):\,\,\, \max_{1\le i\le m }|x_i|\le 1\}
$$
компактно. Тогда sup-норма на пространстве векторов  ${\bf x}$  задает на двумерном подпространстве 
$\pi$   индуцированную норму
   $|\cdot|_\bullet$    (таким образом, "единичный шар"  ${\cal B}=\{ {\bf z} \in \pi:\,\, |{\bf z}|_\bullet = 1\}$
 есть ограниченного выпуклое центрально симметричное множество).  Рассмотрим одномерное подпространство  $\ell\subset \pi$.
Пусть  ${\bf 0}$  есть единственная точка решетки  $\Lambda$,
 которая лежит в  $\ell$.
 Пусть, более того,  
$$
\ell \cap \{{\bf z}=({\bf x}', {\bf y}'):\, x_1'=...=x_m' =0\}
=\{
{\bf 0}\}.
$$
Пусть имеется некоторое  линейное подпространство ${\cal L}$,
  такое что
$$
{\cal L}\cap \pi = \ell,\,\,\,\,\,
{\rm dim }\,{\cal L} \ge 1.
$$
Поскольку  $\ell$  не является вполне рациональным, то ${\cal L}$  тоже
  не является вполне рациональным. 
 
Определим наилучшее приближение точками решетки  $\Lambda$  подпространства  $\ell$  относительно нормы  $|\cdot|_\bullet$
как такую точку  ${\bf z}= ({\bf x}, {\bf y})  \in \Lambda$,
 что 
$${\rm dist }  ({\bf z},{\cal L} ) = \min 
{\rm dist } ({\bf z}',{\cal L} )$$
где минимум берется по всем точкам ${\bf z}'\in \Lambda$ таким,  что
$$
0<|{\bf z}'|_\bullet\le |{\bf z}|_\bullet
$$
 и 
$  
{\rm dist }  ({\bf z},{\cal L} )  
 $ 
обозначает расстояние в sup-норме (в $\mathbb{R}^d$) от точки ${\bf z}$  до подпространства ${\cal L}$.

  Естественно, наилучшие приближения должны располагаться в виде бесконечных последова-тельностей
$$
\pm {\bf z}_1,
\pm {\bf z}_2,..., \pm {\bf z}_\nu,\pm {\bf z}_{\nu+1},... ,
$$
$$
{\rm dist } ({\bf z}_1,{\cal L} )>
{\rm dist }({\bf z}_2,{\cal L} )>
\cdots
>
{\rm dist }({\bf z}_\nu,{\cal L} )>
{\rm dist }({\bf z}_{\nu+1},{\cal L} )
>\cdots ,
$$
$$
|{\bf z}_1|_\bullet <
|{\bf z}_2|_\bullet<\cdots <
|{\bf z}_\nu|_\bullet
<
|{\bf z}_{\nu+1}|_\bullet
<\cdots ,
$$
аналогичных последовательностям  (\ref{1},\ref{2}).
Строгая монотонность  величин ${\rm dist }({\bf z}_\nu,{\cal L} ), |{\bf z}_\nu|_\bullet$
определя-ется постановкой задачи. Будем говорить, что набор
$(\pi,\ell,{\cal L})$  является {\it правильным}, если для каждого достаточно большого значения индекса  $\nu$
целочисленный вектор  ${\bf z}_\nu \in \Lambda$   определен однозначно с точностью до знака.

Очевидным обобщением неравенства (\ref{approx})
 является следующее 

{\bf  Предложение 2.}\,\,{\it 
Если набор
$(\pi,\ell,{\cal L})$  является   правильным, то с некоторыми положительными постоянными  $C_1,C_2$,
 зависящими от  $\pi,\ell$     и ${\cal L}$  при всех  $\nu$ выполняется
$$
C_1  \le
{\rm dist }({\bf z}_\nu,{\cal L} )\cdot
|{\bf z }_{\nu+1}|_\bullet
\le
C_2  
.$$
}

Чтобы убедиться в справедливости предложения 2 достаточно рассмотреть двумерную выпуклую центрально-симметричную область
$$
\{ {\bf z}\in \pi:\,\,
|{\bf z}|_\bullet \le  |{\bf z}_{\nu+1}|_\bullet ,\,\,
{\rm dist }({\bf z} ,  ({\cal L} )
\le {\rm dist } ({\bf z}_\nu,{\cal L}  )
\}
,
$$
на границе которой есть две линейно  независимые точки  ${\bf z}_\nu, {\bf z}_{\nu+1} \in \Lambda$.

\section{Сингулярные системы в задаче о совместных  диофантовых приближенияx.}
\label{sszsdp}

Случай $m = 1$  называется задачей о  совместных диофантовых приближениях.
Всюду ниже в настоящем пункте мы полагаем $m=1$ и используем обозначение
$\theta_{j}^1=\theta_j,\,1\le j \le n$.
В этом
 случае 
$$
{\bf x} = x_1 = x\in \mathbb{R}^1,\,\,\,
M({\bf x})=|x|,\,\,\,
L_j(x) = \theta_j x
$$
и
речь идет о малости величины
$$
\max_{1\le j \le n} || x\theta_j||.
$$
Пустые параллелепипеды (\ref{dopara},\ref{dopara1}) теперь имеют вид
 
$$
\Omega_\nu^{(0)} =\{
{\bf z}'= (x', y_1',...,y_n'):\,\,\,\,
|x'|\le |x_\nu|,\, \max_{1\le j\le n}|x'\theta_j + y_j'|
\le  \max_{1\le j\le n}||x_\nu\theta_j||\}
$$
и
\begin{equation}
\Omega_\nu=
\{
{\bf z}'= (x', y_1',...,y_n'):\,\,\,\,
|x'|\le |x_{\nu+1}|,\, \max_{1\le j\le n}|x'\theta_j +y_j'|
\le  \max_{1\le j\le n}||x_\nu\theta_j||\}
\label{para}
\end{equation}
В первом из них нет других целых точек, кроме точек
${\bf 0},\pm{\bf z}_\nu$. Во втором
нет других целых точек кроме ${\bf 0}$ и $\pm {\bf z}_\nu, \pm {\bf z}_{\nu+1}$.

Отметим, что для вещественнозначной функции
$\psi (t) =o(t^{-1/n})$
набор  вещественных чисел $ \Theta\in \mathbb{R}^n $ будет $\psi$-сингулярным если
 при всяком  достаточно большом $t $ система диофантовых неравенств
$$
 \max_{1\le j\le n}||x\theta_j||\le \psi (t) ,\,\,\,\,\,\,\,\, 0<\max_{1\le i\le m}|x_i|<  t
$$
имеет целочисленное решение ${ x}\in \mathbb{Z}$, или, что то же самое, для всех достаточно больших $\nu$
 для знаменателей $x_\nu$  наилучших совместных приближений выполнено
выполнено
$$
\max_{1\le j\le n}||x_\nu\theta_j||\le\psi (x_{\nu+1}).
$$

Напомним, что 
для вектора  $\Theta = (\theta_1,...,\theta_n)\in \mathbb{R}^n$
через  ${\rm dim}_\mathbb{Z}\Theta$  мы обобзнчаем максимальное количество линейно независимых чисел из набора  $1,\theta_1,...,\theta_n$

\subsection{Независимость векторов наилучших приближений}\label{nvnp}

Сейчас мы сформулируем 
в виде следствий утверждения, которые  мгновенно выводятся из описанных выше результатов. 
 
Из предложения 2 пункта \ref{dr} получаем

{\bf  Следствие 1. }\,\,\,{\it   Пусть ${\rm dim }_\mathbb{Z} \Theta = 2$. Тогда с некоторой постоянной $\gamma = \gamma (\Theta )>0$  выполнено
$$
\zeta_\nu \ge \frac{\gamma}{M_{\nu+1}}, \,\,\,\ \forall \nu \in \mathbb{N}.
$$
}

  Из  утверждения  (ii)  теоремы \ref{NT1} и замечания в конце пункта \ref{or} получаем

{\bf Следствие 2.}
\,\,\,{\it В случае $n\ge 2, {\rm dim}_\mathbb{Z}\Theta\ge 3$ найдется бесконечно много значений $\nu$
таких, что векторы  наилучших приближений ${\bf z}_\nu$, ${\bf z}_{\nu+1}$  ${\bf z}_{\nu+2}$  линейно независимы.
}

{\bf Следствие 3.}
\,\,\,{\it  Если
$n=2,{\rm dim}_\mathbb{Z} \Theta = 3$,
то для бесконечно многих значений $\nu$  
для векторов наилучших приближений
$$
{\bf z}_{\nu+i} =(x_{\nu+i}, y_{1,\nu+i}, y_{2.\nu+i}),\,\,\,
i =0,1,2$$
выполнено
$$
{\rm det}\,
\left(
\begin{array}{ccc}
x_\nu & y_{1,\nu}&y_{2,\nu}
\cr
x_{\nu+1} & y_{1,\nu+1}&y_{2,\nu+1}
\cr
x_{\nu+2} & y_{1,\nu+2}&y_{2,\nu+2}
\end{array}
\right)
\neq 0.
$$}
 Следствие 3 имеется в явном виде у Дж. Лагариаса \cite{LL}. В неявном виде у 
В.Ярника (\cite{JRUS}, c. 333 - 337) фактически имеется более общее утверждение:

{\bf Следствие 4.}
\,\,\,
{\it Пусть $n\ge 2$ и ${\rm dim}_\mathbb{Z} \Theta\ge 3$.  
Тогда для 
некоторых двух  индексов $j_1\neq j_2$ для бесконечно многих значений $\nu$ 
  для векторов наилучших приближений
$$
{\bf z}_{\nu+i} =(x_{\nu+i}, y_{1,\nu+i},..., y_{n.\nu+i}),\,\,\,
i =0,1,2$$
выполняется
\begin{equation}
{\rm det}\,
\left(
\begin{array}{ccc}
x_\nu & y_{j_1,\nu}&y_{j_2,\nu}
\cr
x_{\nu+1} & y_{j_1,\nu+1}&y_{j_2,\nu+1}
\cr
x_{\nu+2} & y_{j_1,\nu+2}&y_{j_2,\nu+2}
\end{array}
\right)
\neq 0.
\label{opr}
\end{equation}

}

Формально у Ярника \cite{JRUS}  присутствует чуть более слабое утверждение.
  Поясним, как   получается следствие 4 (ср. c леммой на стр. 333  из \cite{JRUS}).

Согласно следствию 2 для бесконечного  множества значений индекса $\nu$
 матрица
\begin{equation}
\left(
\begin{array}{cccc}
x_\nu & y_{1,\nu}&\cdots &y_{n,\nu}
\cr
x_{\nu+1} & y_{1,\nu+1}&\cdots&y_{n,\nu+1}
\cr
x_{\nu+2} & y_{1,\nu+2}&\cdots &y_{n,\nu+2}
\end{array}
\right)
\label{matt}
\end{equation}
 имеет ранг 3.
  Если все определители 
(\ref{opr})
 обращаются в ноль, то
для некоторых $j_1,j_2,j_3$
имеем
$$
{\rm det }
\left(
\begin{array}{ccc}
 y_{j_1,\nu}&y_{j_2,\nu}&y_{j_3,\nu}
\cr
 y_{j_1,\nu+1}&y_{j_2,\nu+1}&y_{j_3,\nu+1}
\cr
y_{j_1,\nu+2}&y_{j_2,\nu+2}
&y_{j_3,\nu+1}
\end{array}
\right)
\neq 0,
$$
то есть
трехмерные векторы (=столбцы матрицы (\ref{matt})) 
\begin{equation}
\left(
\begin{array}{c}
 y_{j_1,\nu}\cr
 y_{j_1,\nu+1}
\cr
 y_{j_1,\nu+2}\end{array}
\right),
\left(
\begin{array}{c}
 y_{j_2,\nu}\cr
 y_{j_2,\nu+1}
\cr
 y_{j_2,\nu+2}\end{array}
\right),
\left(
\begin{array}{c}
 y_{j_3,\nu}\cr
 y_{j_3,\nu+1}
\cr
 y_{j_3,\nu+2}\end{array}
\right)
\label{veve}
\end{equation}
 линейно независимы.
 Но тогда получается, что вектор
$$
\left(
\begin{array}{c}
 x_\nu\cr
 x_{\nu+1}
\cr
 x_{\nu+2}\end{array}
\right)
$$
 линейно зависим с любыми двумя векторами из  (\ref{veve}), чего не может быть.

Прежде чем формулировать следствие 5 проведем некоторые рассуждения.
Пусть  ${\rm dim}_\mathbb{Z}\Theta \ge 3$.
 Рассмотрим определитель из следствия 4, не обращающийся в ноль при некотором значении параметра  $\nu$. Тогда

$$
1\le
|
{\rm det}\,
\left(
\begin{array}{ccc}
x_\nu & y_{j_1,\nu}&y_{j_2,\nu}
\cr
x_{\nu+1} & y_{j_1,\nu+1}&y_{j_2,\nu+1}
\cr
x_{\nu+2} & y_{j_1,\nu+2}&y_{j_2,\nu+2}
\end{array}
\right)
|=|
{\rm det}\,
\left(
\begin{array}{ccc}
x_\nu & y_{j_1,\nu}-\theta_{j_1}x_\nu&y_{j_2,\nu}-\theta_{j_2}x_\nu
\cr
x_{\nu+1} & y_{j_1,\nu+1}-\theta_{j_1}x_{\nu+1}&y_{j_2,\nu+1}-\theta_{j_2}x_{\nu+1}
\cr
x_{\nu+2} & y_{j_1,\nu+2}-\theta_{j_1}x_{\nu+2}&y_{j_2,\nu+2}-\theta_{j_2}x_{\nu+2}
\end{array}
\right)|\le
$$
\begin{equation}\label{det33}
\le
6x_{\nu+2}\times\max_{1\le j\le n }||x_\nu \theta_j||
\times\max_{1\le j\le n }||x_{\nu+1} \theta_j||=
6M_{\nu+2}\zeta_\nu\zeta_{\nu+1}.
\end{equation}
  Таким образом,
 для бесконечно многих значений $\nu$  имеем
$$
M_{\nu+2}\zeta_{\nu+1}\ge \frac{1}{6\zeta_\nu}
.
$$
Очевидным становится

{\bf Следствие 5.}
\,\,\,
{\it
Пусть  ${\rm dim}_\mathbb{Z}\Theta \ge 3$.
Тогда
\begin{equation}
\limsup_{\nu\to +\infty} M_{\nu+1}\zeta_{\nu} =+\infty.
\label{tbilli}
\end{equation}
}
Утверждение (\ref{tbilli})  имеется у В.Ярника
в \cite{JTBIL}  при $n=2$ (см. Satz 9 из \cite{JTBIL}).  В общем случае оно сфомулировано в  \cite{J59}.
Отметим, что в \cite{JTBIL}  соотношение (\ref{tbilli})  доказывается
с помощью соображений переноса.  Точнее, В.Ярник везде доказывает  равенство 
\begin{equation}\label{formuja}
\limsup_{t\to+\infty} t\cdot\psi_\Theta (t) =+\infty,
\end{equation}
которое является следствием
соотношения
(\ref{tbilli}).
Здесь мы приводим более точный результат из работы  автора
\cite{msig} (в работе \cite{msig}  автор не ссылается на работы В.Ярника, поскольку тогда он не знал об их существовании).

\begin{theorem}\label{T7}
Пусть $ n \ge 2$ и ${\rm dim} _\mathbb{Z}\Theta \ge 3$.
Тогда
 
\begin{equation}
\zeta_\nu M_{\nu +1} \to +\infty,\hspace{2mm} \nu  \to  +\infty . \label{u1}
\end{equation}

\end{theorem}

Доказательство.

{\bf 1.}
Во-первых, отметим, следующее. Пусть $ \Lambda^2 \in \mathbb{ Z}^{n+1} $  --- произвольная двумерная подрешетка и $ {\rm det}_2 \Lambda^2$  --- ее
двумерный фундаментальный объем. Тогда множество решеток
$$
\{ \Lambda^2 \subset \mathbb{ Z}^{n+1}:\hspace{2mm} {\rm det}_2 \Lambda^2 \le \gamma \}
$$
конечно для любого положительного $\gamma $.

{\bf 2.} Во-вторых, если рассмотреть  двумерную решетку $ \Lambda^2_\nu  = \langle {\bf z}_\nu ,  {\bf z}_{\nu +1} \rangle_\mathbb{ Z} $, порожденную двумя
последовательными наилучшими приближениями, то из  условия $ {\rm conv} ({\bf 0}, {\bf z}_\nu ,  {\bf z}_{\nu +1}) \subset \Omega_\nu $ 
(область  $\Omega_\nu$ определена в  (\ref{para}))
следует, что  с некоторойположительной постоянной  $C(\Theta ) $  выполнено

$$
\frac{1}{2} {\rm det}_2 \Lambda^2_\nu = {\rm vol}_2 ({\rm  conv }(0,{\bf z} _\nu , {\bf z} _{\nu +1} )) \le C(\Theta )  \zeta_\nu M_{\nu +1}.
$$

{\bf 3.}
В третьих,     из условия ${\rm dim}_\mathbb{ Z} \Theta \ge 3$ следует  что для каждого фиксированного $\mu$ при
достаточно большом $\nu$ будет выполнено ${\bf z}_\nu \not \in
  \Lambda^2_{\mu} $. Значит, при каждом $\mu$ и
при
 достаточно
большом $\nu$ двумерная решетка  $\Lambda^2_{\nu} $ не будет совпадать ни с одной из решеток $
 \Lambda^2_{1}, \Lambda^2_{2},..., \Lambda^2_{\mu} $.

Теорема  \ref{T7} вытекает из сказанного выше.

 Далее мы формулируем следствие теоремы \ref{T7}, в некотором смысле аналогично предложению 2 из пункта
\ref{dr}.

 По аналогии с тем, что проделано в  пункте   \ref{dr}
рассмотрим вполне рациональное линейное подпространство  $\pi\subset \mathbb{R}^d$
такое, что  ${\rm dim }\,\pi \ge 3$  и пересечение 
$\pi \cap \{{\bf z}=({\bf x},{\bf y}):\,\,\, \max_{1\le i\le m }|x_i|\le 1\}
$
ограничено,
и  решетку  $\Lambda =\pi \cap \mathbb{Z}^d$.
Пусть снова
   $\ell\subset \pi$ есть некоторое одномерное линейнное подпространство в $\pi$.
Но сейчас (в отличие от того, что было в пункте   \ref{dr})
мы будем предполагать, что
 подпространство  $\ell$ 
не лежит ни в каком собственном рациональном подпространстве подпро-странства  $\pi$
 (это нечто большее, чем то, что ${\bf 0}$  есть единственная точка решетки  $\Lambda$,
 которая лежит в $\ell$).

Для индуцированной нормы
$|\cdot|_\bullet$ на  $\pi$ 
рассмотрим последовательность наилучших приближений 
${\bf z}_\nu \in \Lambda $
 к подпространству  $\ell$  относительно нормы  $|\cdot|_\bullet$, соответствующую подпространству ${\cal L}\supset\ell ={\cal L}\cap\pi$
(тут определение наилучших приближений и правильного набора  $(\pi,\ell,{\cal L})$ совершенно аналогично
 определеню, данному в пункте \ref{dr}).

{\bf  Предложение 3.}\,\,{\it  
В описанных выше условиях 
для
правильного набора  $(\pi,\ell,{\cal L})$
выполняется
 $$ 
{\rm dist } ({\bf z}_\nu,{\cal L})\cdot
|{\bf z }_{\nu+1}|_\bullet
\to +\infty,\,\,\,\,
\nu \to +\infty.
$$
}

\subsection{Вырождение размерности наилучших совместных приближений.}\label{vyro}

Сформулируем результат, полученный автором в \cite{MUMN} (см. также \cite{ME}),
представляющий контрпример  к  предположению Дж. Лагариаса \cite{LL}.

\begin{theorem}\label{antilaga}
Пусть 
$n\ge 3$.  Тогда найдутся числа  $\theta_1,...,\theta_n$,
линейно независимые вместе с единицей над $\mathbb{Z}$
 и такие, что среди любых  $n+1$  идущих подряд векторов наилучших приближений
${\bf z}_\nu
, {\bf z}_{\nu+1},...,{\bf z}_{\nu+n}$
 можно выбрать не более трех линейно независимых.
\end{theorem}

Прокомментируем вкратце доказательство этой теоремы.
Набор чисел
$\theta_1,...,\theta_n$
задается с помощью рациональных приближений к нему. Сами рациональные приближения строются индуктивно.
Основным инструментом для осуществления индуктивного шага конструкции является следующее утверждение.

\begin{lem}\label{vivide}
Пусть
для рациональной  точки  $\alpha = \left( \frac{a_1}{q},...,\frac{a_n}{q}\right)$, $(a_1,...,a_n,q) = 1$
выполнено следующее:

{\rm (i)}
конечная последовательность совместных наилучших приближений
$$
{\bf z}_\nu =
(x_\nu,y_{1,\nu},...,y_{n,\nu}) ,\,\,\,\, 1\le \nu \le k,
$$
$$
(x_k,y_{1,k},...,y_{n,k})= (q,a_{1},...,a_{n}) 
$$
определена однозначно;

{\rm (ii)}  найдется двумерное вполне рациональное подпространство  $\pi$ и
номер  $\mu <k$ 
 такие, что
$$
{\bf z}_\nu \in \pi ,\,\,\, \mu\le \nu \le k;
$$

{\rm (iii)} 
для наилучшего приближения с номером  $k$  выполнено
$$
\zeta_k < \rho (\pi )
$$
(величина $ \rho (\pi )$ определена в (\ref{ropi})).

Тогда найдется  положительное  $\varepsilon$ такое, что для всех точек $\beta$,
лежащих в  $\varepsilon$-окрестности точки  $\alpha$,  будет выполнено следующее:

{\rm (i$^*$)}  все векторы наилучших приближений к $\alpha$  будут  наилучшими приближениями и  к $\beta$ тоже;

{\rm (ii$^*$)}
все расширенные векторы наилучших приближений 
к $\beta$,
знаменатель $x$ которых находится  в промежутке 
$x_\mu \le x\le x_k$,
лежат в подпространстве  $\pi$.

\end{lem}

Для осуществления индуктивного шага построения следует выбрать двумерное вполне рацио-нальное подпространство  $\pi_1\ni \alpha$ , а затем  новую рациональную точку $\alpha_1 \in \pi_1$,
лежащую близко к   $\alpha$  так, чтобы она, в свою очередь,
снова удовлетворяла условиям леммы \ref{vivide}  с заменой  $\pi $ на $\pi_1$.
Причем  $\alpha_1$  следует выбирать так, чтобы для нее  имелось  $>n+1$
векторов
 наилучших приближений, лежащих в  подпространстве  $\pi_1$.
Независимость вместе с единицей над  $\mathbb{Z}$
 компонент предельного вектора надо обеспечить посредством выбора  двумерных подпространств так, чтобы
$R(\Theta) = n+1$.

\section{Сингулярность и диофантов тип.}\label{sdty}

 Пусть функция  $\varphi (t)$  монотоно убывает к нулю при  $t\to +\infty$.
 Принято говорить, что набор $\Theta  $  имеет диофантов тип  $\le \varphi (t)$ 
  если для бесконечно многих   ${\bf x}\in\mathbb{Z}^m$  выполняется
$$
\max_{1\le j \le n}||L_j({\bf x})||\le \varphi \left(\max_{1\le i\le m}|x_i|\right).
$$
В работе \cite{JRUS}
 В.Ярник установил простые и замечательные оценки диофантового типа для сингулярных матриц.
Ниже мы их приводим. Как показал М.Лоран  \cite{lora}, неравенства В.Ярника в случаях
$m=1,n=2$ и $m=2,n=1$  неулучшаемы. Точную формулировку теоремы М.Лорана мы приводим в пункте \ref{LORR} нашей статьи.
Ниже в пунктах \ref{mpavno1}, \ref{mravno2}, \ref{borovo}  мы приводим оригинальные результаты В.Ярника.  В пунктах
\ref{vkonetz}, \ref{novoslu}, \ref{echshe}
 имеется усиление теоремы В.Ярника, принадлежащее автору и установленное в препринте \cite{moJJ}.
 
\subsection{ Случай $m=1$.}\label{mpavno1}
 
Сначала рассмотрим задачу о совместных приближениях. В.Ярник
(см. \cite{JRUS}, Теорема 3) получил следующий результат:

\begin{theorem}\label{T8}
Пусть $\psi (t)$   непрерывная убывающая к нулю  при $t\to +\infty$ функция вещественного переменного $t$, причем
 функция $t \psi (t)$ монотонно возрастает к бесконечности при  $t\to +\infty$.
 Обозначим через $\omega (t)$ функцию, обратную к  $t\psi (t)$.
Положим
$$
 \varphi^{[\psi ]}( t) = \psi \left(\omega\left(\frac{1}{6\psi(t)}\right)\right).
$$
  Пусть 
$n \ge 2, {\rm dim }_\mathbb{Z}\Theta\ge 3 $ и
набор  $\Theta$
является $\psi$-сингулярным.
Тогда для бесконечного
  количества натуральных чисел выполнено
$$
  \max_{1\le j \le n}||x\theta_j||\le \varphi^{[\psi ]} (x).
$$
\end{theorem}

Сейчас
  мы докажем более точный факт, а именно, если
наилучшие приближения
 ${\bf z}_\nu$, ${\bf z}_{\nu+1}$  ${\bf z}_{\nu+2}$  линейно независимы, то
\begin{equation}
\zeta_{\nu+1} \le  \varphi^{[\psi ]} (x_{\nu+1}),
\label{uu}
\end{equation}
из которого теорема \ref{T8} следует   в свете следствия 2 из пункта \ref{nvnp}.
 
Действительно,  рассмотрим определитель из следствия 4   пункта \ref{nvnp},  не обращающийся в ноль. Тогда
из неравенства (\ref{det33}), 
 используя предложение 1 из пункта \ref{ssinp}, получаем
$$
1\le  6 M_{\nu+2}\psi (M_{\nu+2})\psi (M_{\nu+1}) =
6 x_{\nu+2}\psi (x_{\nu+2})\psi (x_{\nu+1}), 
$$
откуда  (\ref{uu})  сразу следует.

 Чтобы пояснить результат теоремы \ref{T8}, сформулируем следствие о диофантовых экспонентах, также принадлежащее В. Ярнику
(см. теорему 1, часть I из \cite{JRUS}).

{\bf Следствие.}\,\,\,{\it
 Пусть $n \ge 2$ и ${\rm dim}_\mathbb{Z}\Theta\ge 3$.
Обозначим через  $\alpha (\Theta)$ и $\beta (\Theta)$
 супремумы тех  $\gamma$,
 для которых  для функции Ярника (\ref{fj}), соответственно, выполнено 
\begin{equation}\label{alef}
\limsup_{t\to +\infty} t^\gamma
\psi_\Theta (t) <+\infty,
\end{equation}
\begin{equation}\label{bet}
\liminf_{t\to +\infty} t^\gamma
\psi_\Theta (t) <+\infty
.
\end{equation}
(Ясно, что $
\frac{1}{n}\le\alpha (\Theta)\le \beta (\Theta)\le +\infty.
$)
 Тогда если $\alpha (\theta ) < 1$,   то
\begin{equation}\label{betatet}
 \beta (\Theta ) \ge \frac{\alpha^2(\Theta )}{1-\alpha (\Theta)};
\end{equation}
если же  $\alpha (\Theta ) = 1$,  то $\beta (\Theta) = +\infty. $ 
}

\subsection{ Случай $m=1, n = 3$.}
\label{vkonetz}

Для
$\alpha \in [1/3,1]$ положим
$$
g_1(\alpha ) =
\frac{(1-\alpha)\alpha+\sqrt{(1-\alpha)^2\alpha^2+4\alpha(2\alpha^2-2\alpha+1)}}{4\alpha^2-4\alpha+2}
.$$
Величина  $g_1(\alpha)$  есть корень уравнения 
$$
(2\alpha^2-2\alpha +1)x^2+\alpha(\alpha -1)x -\alpha = 0.$$
Легко видеть, что
\begin{equation}\label{doradoo}
g_1(\alpha ) =
\max_{\delta \ge1,\,\gamma (1-\alpha)-\alpha>0}\,\,
\min \left\{
\delta,\frac{\alpha}{\gamma (1-\alpha)-\alpha},\frac{1}{\alpha}-\left(\frac{1}{\alpha}-1\right)
\frac{\delta}{\gamma}\right\},
\end{equation}
$$
g_1(1/3) = g_1(1) = 1
.$$
Более того, для 
 $ \alpha \in (1/3,1)$ выполнено $ g_1(\alpha) >1$.
 
Пусть
 $\alpha_0 \in (1/2,1)$ есть корень уравнения
$$
x^3-x^2+2x-x = 1.
$$
Видим, что в интервале
 $1/3<\alpha<\alpha_0$
 выполняется
$$
g_1(\alpha) > \max \left\{ 1, \frac{\alpha}{1-\alpha }\right\}.
$$

\begin{theorem}\label{Doradoo}
Пусть
  $m=1,n=3$  и матрица $
\Theta=\left(\begin{array}{c}
\theta_1\cr \theta_2 \cr \theta_3\end{array}\right)$ состоит из чисел, линейно независимых вместе с единицей над $\mathbb{Z}$. Тогда
\begin{equation}\label{novoe0JA}
\beta (\Theta) \ge \alpha (\Theta) g_1 (\alpha (\Theta)).
\end{equation}
 
\end{theorem}

Видим, что неравенство 
  (\ref{novoe0JA})  сильнее, чем результат теоремы   \ref{T8} (и ее следствия) в интервале
  $1/3 < \alpha(\Theta )<\alpha_0$.

Доказательство.
Будем предполагать, что
$
\psi_{\Theta}(t) \le \psi (t)$
 выполнено с некоторой непрерывной убывающей к нулю функцией $\psi(t)$.
Также следует предположить, что функция 
$t\mapsto t\cdot \psi (t)$
монотонно возрастает к бесконечности при
  $t\to +\infty$.

Рассмотрим последовательность векторов наилучших приближений 
  ${\bf z}_\nu = (x_\nu,y_{1,\nu},y_{2,\nu},y_{3,\nu})$.
Из условия линейной независимости видим, что для бесконечной последовательности пар индексов  
    $\nu<k, \nu\to +\infty$   будет выполнено 

$\bullet $  обе тройки 
$$
{\bf z}_{\nu-1},{\bf z}_\nu,{\bf z}_{\nu+1};\,\,\,\,\,\, 
{\bf z}_{k-1},{\bf z}_k,{\bf z}_{k+1}
$$ 
состоят из линейно независимых векторов;  

$\bullet$ 
имеется некоторое двумерное подпространство $ \pi $, такое что  
$$
{\bf z}_l\in \pi,\,\,\, \nu\le l\le k;\,\,\,\,\, {\bf z}_{\nu-1} \not\in \pi,\,\,\, {\bf z}_{k+1} \not\in \pi;
$$

$\bullet$  векторы 
$$
{\bf z}_{\nu-1},{\bf z}_\nu,{\bf z}_{k},{\bf z}_{k+1}
$$
линейно независимы.

Итак
 $$
1\le |
{\rm det}
\left(
\begin{array}{cccc}
y_{1,\nu-1}&y_{2,\nu-1}&y_{3,\nu-1}& x_{\nu-1}\cr
y_{1,\nu}&y_{2,\nu}&y_{3,\nu}& x_{\nu}\cr
y_{1,k}&y_{2,k}&y_{3, k}& x_{k}\cr
y_{1,k+1}&y_{2,k+1}&y_{3,k+1}& x_{k+1}
\end{array}
\right)|
\le
$$
\begin{equation}\label{promejutok0JA}
\le
24 \zeta_{\nu-1}\zeta_\nu\zeta_{k} M_{k+1} \le 
24 \psi (M_\nu) \psi (M_{\nu +1})\psi (M_{k+1})M_{k+1}<
24 \psi (M_\nu) (\psi (M_{\nu +1})^2M_{k+1}
.
\end{equation}
 Далее надо рассмотрет три случая.

{\bf 1$^0$.} Для некоторого $\gamma >1$  имеется бесконечно много пар   $(\nu,k)$
рассматриваемого вида, и таких что
$$
M_{k+1}\le M_{\nu+1}^{\gamma}.
$$
Из неравенства   (\ref{promejutok0JA})  выводим
 $$
   \frac{1}{24\psi(M_\nu)}\le  M_{\nu+1}^\gamma\cdot \psi (M_{\nu+1})\cdot \psi (M_{\nu+1}^\gamma )
$$
Дополнительно предполагая возрастание функции
$t\mapsto t^\gamma \cdot \psi (t)\cdot \psi (t^\gamma )$  
и обозначая обратную к ней через
$\rho (t)$, видим, что
\begin{equation}\label{aJA}
\zeta_\nu\le \psi (M_{\nu+1})\le
\psi\left(
\rho \left(
\frac{1}{24\psi (M_\nu)}
\right)
\right).
 \end{equation}

{\bf 2$^0$.} Для некоторого $\delta \ge 1$ имеется бесконечно много рассматриваемых пар индексов  $(\nu,k)$,
таких что
$$
M_{k+1}\ge  M_k^{\delta}.
$$
Тогда сразу получаем
\begin{equation}\label{bJA}
\zeta_k \le \psi (M_{k+1}) \le \psi ( M_k^{\delta}).
\end{equation}

{\bf 3$^0$.} Для бесконечного количества рассматриваемых пар индексов   $(\nu,k)$
выполнено
$$
 M_{\nu+1}^\gamma \le M_{k+1}\le M_k^\delta.
$$
 Тогда для двумерной 
решетки  $\Lambda =\mathbb{Z}^4\cap \pi$ 
при 
 $\nu\le l \le k-1$  будет выполняться
$$
\zeta_{l}M_{l+1} \asymp_\Theta {\rm det}_2\Lambda.
$$
 Отсюда получаем неравенства
\begin{equation}\label{cJA}
 \zeta_{k-1}\ll_\Theta
\frac{M_{\nu+1}\psi (M_{\nu+1})}{M_k}\ll_\Theta
M_k^{\frac{\delta}{\gamma}-1}\psi (M_k^{\frac{\delta}{\gamma}})\le
M_{k-1}^{\frac{\delta}{\gamma}-1}\psi (M_{k-1}^{\frac{\delta}{\gamma}})
 \end{equation}
 (конечно, надо дополнительно предполагать монотонное убывание   функции
$t\mapsto t^{\frac{\delta}{\gamma}-1}\psi (t^{\frac{\delta}{\gamma}})$).

Для завершения доказательства надо рассмотреть малое
 $\varepsilon >0$ и положить $\psi(t) =t^{-\alpha+\varepsilon}$.
 В свете соотношения  (\ref{doradoo})  неравенство  (\ref{novoe0JA}) вытекает из  (\ref{aJA},\ref{bJA},\ref{cJA}).

\subsection{ Случай $m=2$.}\label{mravno2}

Сформулируем теорему В.Ярника о диофантовом типе из работы \cite{JRUS}, аналогичную теореме \ref{T8}.

\begin{theorem}
\label{m=2} 
Пусть  $ m= 2, n \ge 1$ и невырожденная матрица
$$
\Theta= 
\left(
\begin{array}{cc}
\theta_1^1&\theta_1^2\cr

\theta_2^1&\theta_2^2\cr
\vdots &\vdots\cr
\theta_n^1 &\theta_n^2
\end{array}
\right)
$$
является
  $\psi$-сингулярной с некоторой монотонно убывающей  функцией  $\psi (t) = o(t^{-1}), t\to +\infty$.
 Положим
$$
\varphi^{[\psi ]}_2(t) = \psi\left(\frac{1}{ 6t\psi (t)}\right)
$$
Тогда   
 найдется бесконечно много	  целочисленных векторов
${\bf x} = (x_1,x_2)$ таких, что
$$
\max_{1\le j\le n}
||\theta_j^1x_1+\theta_j^2x_2||\le
\varphi^{[\psi ]}_2\left(\max_{i=1,2}|x_i|\right) .
$$
 
\end{theorem}

Отметим, что в условиях теоремы для правильной матрицы $\Theta$  выполнено  $R(\Theta)\ge 3$.
(При  $n=1$ это следует из утверждения (iii)  теоремы \ref{NT1}. 
При  $n\ge 2$  это  вытекает  из $\psi$-сингулярности и из теоремы \ref{ngem}.)
Теперь
для доказательства надо
установить аналог следствия 4   
 пункта \ref{nvnp}. Видим, что для некоторого индекса 
$k$ и
для бесконечно многих значений  $\nu$
составленный из координат векторов наилучших приближений
  определитель
$$
\left|
\begin{array}{ccc}
x_{1,\nu} & x_{2,\nu} & y_{k,\nu}
\cr
x_{1,\nu+1} & x_{2,\nu+1} & y_{k,\nu+1}
\cr
x_{1,\nu+2} & x_{2,\nu+2} & y_{k,\nu+2}
\end{array}
\right|
$$
не обращается в ноль
(на самом деле это место требует более подробного объяснения).
Затем   для такого значения  $\nu$ надо
провести оценку, аналогичную (\ref{det33}):
$$
1\le 6M_{\nu+2}M_{\nu+1}\zeta_\nu\le
6M_{\nu+2}M_{\nu+1}\psi (M_{\nu+1})
.$$
 Затем надо  воспользоваться монотоностью функции  $\psi $.

Мы
не останавливаемся здесь на доказательстве подробно (например, следует разобраться с ситуацией, когдя матрица $\Theta$  не являеися правильной).
Все детали имеются в оригинальной работе В.Ярника.

Из теоремы  \ref{m=2}  В.Ярник выводит следующее утверждение.

{\bf  Следствие.}\,\,
{\it
Пусть  $ m= 2, n \ge 1$ и   матрица $\Theta$ невырожденная.
Тогда для величин
 $\alpha(\Theta),\beta(\Theta)$,
определяемых как супремум тех  $\gamma$,
для которых выполнено  (\ref{alef}) и  (\ref{bet}) соответственно,
 имеет место соотношение
\begin{equation}\label{lgl}
\beta (\Theta) \ge \alpha (\Theta) (\alpha (\Theta)- 1).
\end{equation}
}

\subsection{ Случай $m=n=2$}\label{novoslu}

Для $\alpha \ge 1$ 
положим 
$$
g_3 (\alpha ) =
\frac{1-\alpha +\sqrt{(1-\alpha)^2+4\alpha (2\alpha^2-2\alpha+1)}}{2\alpha} 
.
$$
Итак $g_3 (\alpha ) $
есть решение уравнения
\begin{equation}\label{eqXXX}
\alpha x^2+(\alpha-1)x - (2\alpha^2-2\alpha+1) =0.
\end{equation}
Видим, что $g_3   (1) = 1$ и при  $\alpha >1$ выполнено $g_3(\alpha ) >1$.
Более того, в интервале
$$1\le\alpha < \left(\frac{1+\sqrt{5}}{2}\right)^2
$$
имеем
$$g_3 (\alpha ) >\max (1, \alpha - 1).
$$ 
\begin{theorem}\label{xoXO}
Пусть
числа $\theta_j^i,\,\, i,j=1,2$ линейно независимы вместе с единицей над $\mathbb{Z}$.
рассмотрим матрицу
$$
\Theta =
\left(
\begin{array}{cc}
\theta^1_1&\theta^2_1\cr
\theta^1_2&\theta^2_2
\end{array}
\right)
$$
Тогда
$$
\beta(\Theta) \ge \alpha(\Theta)g_3(\alpha(\Theta)).
$$
 \end{theorem}
Теорема
\ref{xoXO} сильнее теоремы \ref{m=2}  при   $\alpha (\Theta) \in \left(1,\left(\frac{1+\sqrt{5}}{2}\right)^2\right)$.

Доказательство.

 Если 
 $R(\Theta) = 2$ то из теоремы \ref{ngem} 
получаем $\alpha (\Theta ) = 1$, и доказывать нечего.
Согласно следствию 4 из пункта \ref{or} не бывает того, чтобы
 $R(\Theta) =3$.
Итак  $R(\Theta) = 4$.       Тогда  для бесконечнойпоследовательности пар индексов                                                                                       
 
   $\nu<k, \nu\to +\infty$    имеем следующее:

$\bullet $  обе тройки
$$
{\bf z}_{\nu-1},{\bf z}_\nu,{\bf z}_{\nu+1};\,\,\,\,\,\, 
{\bf z}_{k-1},{\bf z}_k,{\bf z}_{k+1}
$$ состоят из линейно независимых векторов;  

$\bullet$  имеется двумерное подпространство    $\pi$, такое что
$$
{\bf z}_l\in \pi,\,\,\, \nu\le l\le k;\,\,\,\,\, {\bf z}_{\nu-1} \not\in \pi,\,\,\, {\bf z}_{k+1} \not\in \pi;
$$

$\bullet$  векторы
$$
{\bf z}_{\nu-1},{\bf z}_\nu,{\bf z}_{k},{\bf z}_{k+1}
$$
линейно независимы.

Тогда
 $$
1\le |
{\rm det}
\left(
\begin{array}{cccc}
x_{1,\nu-1}&x_{2,\nu-1}&y_{1,\nu-1}& y_{2,\nu-1}\cr
x_{1,\nu}&x_{2,\nu}&y_{1,\nu}& y_{2,\nu}\cr
x_{1,k}&x_{2,k}&y_{1,k}& x_{2,k}\cr
x_{1,k+1}&x_{2,k+1}&y_{1,k+1}& y_{2,k+1}
\end{array}
\right)|
\le
24 \zeta_{\nu-1}\zeta_\nu M_k  M_{k+1}  
.
$$
предполагаем монотонное убывание функции 
$t\mapsto t\cdot \psi (t)$. Пусть снова
$
\psi_\Theta (t) \le \psi (t)$. Тогда
\begin{equation}\label{uroXXX}
1\le 24 M_{k+1}M_k\psi (M_{\nu+1})\psi (M_\nu).
\end{equation}
Рассматриваем два случая.

{\bf 1$^0$.} 
 Для некоторого  $\gamma >1$ для бесконечно многих рассматриваемых пар  $(\nu,k)$
выполнено
$$
M_{k+1}\ge  M_k^{\gamma}.
$$
Тогда немедленно получаем
\begin{equation}\label{bXXX}
\zeta_k \le \psi (M_{k+1}) \le \psi ( M_k^{\gamma}).
\end{equation}

{\bf 2$^0$.} Для бесконечного количества рассматриваемых пар   $(\nu,k)$
выполнено
$$
  M_{k+1}\le M_k^\gamma.
$$
Тогда  из 
 (\ref{uroXXX})  видим, что
\begin{equation}\label{urodXXX}
M_k \ge (\psi (M_\nu ))^{-\frac{2}{1+\gamma}}.
\end{equation}
Рассмотрим двумерную подрешетку
  $\Lambda = \pi \cap \mathbb{Z}^4$ 
с определителем 
${\rm det }\,\Lambda$.
Для  точки ${\bf z}\in \pi$ расстояние от  ${\bf z}$ до двумерного подпротсранства 
$$
{\cal L} =\{ {\bf z}=(x_1,x_2,y_1,y_2):\,\,
\theta_1^1x_1+\theta^2_1x_2+y_1 =
\theta_2^1x_1+\theta^2_2x_2+y_2=0\}
$$
пропорционально расстоянию от нее до одномерного подпространства
${\cal L}\cap \pi$.
Пусть $\delta$ 
есть коэффициент этой пропорциональности.  
Параллелепипед
$$
\{
{\bf z}=(x_1,x_2,y_1,y_2):\,\,
|x|< M_{l+1},\,\, \max_{1\le j\le 3}|\theta_jx-y_j|<\zeta_l\}
$$
не содержит внутри себя ненулевых целых точек.
 Значит
  \begin{equation}\label{urodaXXX}
\gamma_1(\Theta
)  \delta \,{\rm det}\,\Lambda \le 
\zeta_l M_{l+1} \le \gamma_2(\Theta
)  \delta \,{\rm det}\,\Lambda ,\,\,\,\,\,\,
\nu\le l\le k-1
\end{equation}
 с некотороыми постоянными 
$\gamma_i(\Theta),\, i =1,2$.
Из (\ref{urodXXX}) и (\ref{urodaXXX})  выводим
\begin{equation}\label{cXXX}
 \zeta_\nu \ll_\Theta  
\frac{\psi (M_k)M_k}{M_{\nu+1}}\ll_\Theta  
M_\nu^{-1}(\psi (M_\nu))^{-\frac{2}{1+\gamma}}\psi\left( (\psi(M_\nu))^{-\frac{2}{1+\gamma}}\right).
\end{equation}
Далее рассматриваем функцию
$ \psi (t) = t^{-\alpha(\Theta)+\varepsilon}$ с малым положительным $\varepsilon$.
Поскольку величина
$\gamma =g_3(\alpha (\Theta)) $ удовлетворяет  (\ref{eqXXX}),
получаем утверждение теоремы \ref{xoXO}.

\subsection{Случай $m >2$.}\label{borovo}

В случае  $m>2$  В.Ярник  в \cite{JRUS} использует более громоздкие рассуждения Поэтому он формулирует и доказывает не общий результат с произвольной функцией 
$\psi$,  а только лишь утверждение про диофантовы показатели. 

\begin{theorem}\label{mbo2}
Пусть  $m\ge 3$  и про невырожденную матрицу  $\Theta$  известно, что
$$
\alpha (\Theta) > (5m^2)^{m-1}.
$$
Тогда
$$
\beta (\Theta ) \ge 
(\alpha(\Theta))^\frac{m}{m-1} - 3\alpha(\Theta). 
$$

\end{theorem}

Доказательство этой теоремы мы здесь не приводим. Оно основывается на изящном построении последовательности {\it независимыx} векторов
 приближений. 
У автора нет сомнений в том, что результат теоремы \ref{mbo2}  не является оптимальным и может быть улучшен.
В следующем пункте мы как раз и приводим подобного рода улучшение (в частном случае).

\subsection{Случай $m =3,n=1$.}\label{echshe}
 
Для   $\alpha\ge 3$ определим функции 
$$
g_2(\alpha ) =\sqrt{\alpha+\frac{1}{\alpha^2}-\frac{7}{4}}+\frac{1}{\alpha}-\frac{1}{2},\,\,\,\,
h(\alpha ) = \alpha - g(\alpha ) -1.
$$
 Заметим, что функции  $g_2(\alpha) $  и  $h(\alpha)$  монотонно возрастают к бесконечности при 
 $\alpha \to+\infty$  и
$$
g_2(3) = h(3) = 1,\,\,\,\,  g_2(\alpha ) \le \alpha - 2.
$$

\begin{theorem}\label{m3n1}
Рассмотрим матрицу-строку $\Theta = (\theta^1,\theta^2,\theta^3)$,
состоящую из линейно независимых вместе с единицей над  $\mathbb{Z}$  чисел.
Тогда  показатели  $\alpha(\Theta)$  и  $\beta (\Theta)$  связаны соотношением
\begin{equation}\label{novoe}
\beta (\Theta) \ge \alpha (\Theta) g_2(\alpha (\Theta)).
\end{equation}

\end{theorem}

Доказательство.

Мы будем рассматривать последовательности векторов  наилучших приближений ${\bf z}_\nu$.

 Сначала рассмотрим случай, когда  $R(\Theta) = 3$. Тогда все векторы  ${\bf z}_\nu$,  начиная с какого-то номера
лежат в некотором трехмерном вполне рациональном подпространстве  $\pi$, и
 мы фактически
имеем 
дело с последовательностью наилучших приближений точками {\it трехмерной}
решетки  $\mathbb{Z}^4\cap \pi $ 
{\it двумерного }
подпространства  $\pi \cap {\cal L}(\Theta)$ (напомним, что определение подпространства ${\cal L}(\Theta)$
 мы давали в пункте \ref{or}). Тогда применима теорема 
\ref{m=2}, и ее следствие дает оценку
(\ref{lgl}),
что сильнее  чем (\ref{novoe}).

Теперь рассмотрим случай, когда  $R(\Theta) = 4$.
В этом случае найдется бесконечно много пар индексов  $\nu<k \,\,\,(\nu\to+\infty )$  таких, что

$\bullet $ тройки
$$
{\bf z}_{\nu-1},{\bf z}_\nu,{\bf z}_{\nu+1};\,\,\,\,\,\, 
{\bf z}_{k-1},{\bf z}_k,{\bf z}_{k+1}
$$
 состоят из линейно независимых векторов каждая; 

$\bullet$ 
имеется некоторое двумерное вполне рациональное подпространство  $\pi$  такое, что
$$
{\bf z}_l\in \pi,\,\,\, \nu\le l\le k;\,\,\,\,\, {\bf z}_{\nu-1} \not\in \pi,\,\,\, {\bf z}_{k+1} \not\in \pi;
$$

$\bullet$  четыре вектора
$$
{\bf z}_{\nu-1},{\bf z}_\nu,{\bf z}_{\nu+1},{\bf z}_{k+1}
$$
линейно независимы.

 Тогда, предполагая, что  набор  $\Theta$  
является  $\psi$-сингулярным
(естественно, считаем, что функция $\psi (t)$  монотонна,
и что монотонна функция $t\mapsto t\psi (t)$), получаем
\begin{equation}\label{promejutok}
1\le|
{\rm det}
\left(
\begin{array}{cccc}
x_{1,\nu-1}&x_{2,\nu-1}&x_{3,\nu-1}& y_{\nu-1}\cr
x_{1,\nu}&x_{2,\nu}&x_{3,\nu}& y_{\nu}\cr
x_{1,\nu+1}&x_{2,\nu+1}&x_{3,\nu+1}& y_{\nu+1}\cr
x_{1,k+1}&x_{2,k+1}&x_{3,k+1}& y_{k+1}
\end{array}
\right)|
\le
24 \zeta_{\nu-1}M_\nu M_{\nu+1}
M_{k+1} \le 
24 \psi (M_\nu) M_\nu M_{\nu+1}M_{k+1}.
\end{equation}
Далее рассматриваем случаи.

{\bf 1$^0$.}  Если для бесконечного количества рассматриваемых пар индексов  $(\nu,k)$
 выполнено
$$
M_{k+1}\le M_\nu^{h(\alpha(\Theta))},
$$
то для таких пар индексов
из (\ref{promejutok})  получаем
$$
M_{\nu+1} \ge \frac{1}{24\psi(M_\nu)M_\nu^{1+h(\alpha(\Theta))}},
\,\,\,\,
\zeta_\nu \le \psi\left(\frac{1}{24\psi(M_\nu)M_\nu^{1+h(\alpha(\Theta))}}
\right).
$$
Из последнего неравенства требуемое соотношение (\ref{novoe}) сразу вытекает.

{\bf 2$^0$.}  Пусть для бесконечного количества рассматриваемых пар индексов  $(\nu,k)$
 выполнено
$$
M_{k+1}\ge  M_k^{g_2(\alpha(\Theta))},
$$
В этом случае сразу получаем
$$
\zeta_k \le \psi (M_{k+1}) \le \psi ( M_k^{g_2(\alpha(\Theta))}),
$$
откуда тоже сразу получаем соотношение (\ref{novoe}).

{\bf 3$^0$.}  Пусть для бесконечного количества рассматриваемых пар индексов  $(\nu,k)$
 выполнено
$$
M_\nu^{h(\alpha(\Theta))}\le
M_{k+1}\le  M_k^{g_2(\alpha(\Theta))},
$$
В этом случае надо более детально исследовать, что происходит в подпространстве  $\pi$.
 Рассмотрим
проекцию подпространства  $\pi$  на подпространство  ${\cal L}(\Theta)$.
Это есть некоторое двумерное линейное подпространство  $\pi^*$,
 пересекающееся с подпространством  $\pi$  по некоторой прямой  $\ell = \pi\cap \pi^*$.
Для точки ${\bf z}\in \pi$
 расстояния от ${\bf z}$
до подпространства ${\cal L}(\Theta)$ и от ${\bf z}$
до подпространства $\ell$  пропорциональны.
Обозначим через  $\delta$ коэффициент этой пропорциональности ("угол" между  пересекающимися двумерными подпространствами
 $ \pi$  и  $\pi^*$).
Так как векторы  ${\bf z}_l$  являются векторами наилучших приближений, то как векторы 
двумерной
решетки
 $\Lambda =\mathbb{Z}^4 \cap \pi$
они автоматически оказываются наилучшими приближениями  точками решетки  $\Lambda$
 прямой  $\ell$  в индуци-рованной норме.
Пусть  ${\rm det}\,\Lambda$  есть фундаментальный объем решетки  $\Lambda$.
 Ясно, что с некоторыми положительными постоянными
$\gamma_i(\Theta),\, i =1,2$
выполнено
$$
\gamma_1(\Theta
)\cdot \delta \,{\rm det}\,\Lambda \le 
\zeta_lM_{l+1} \le \gamma_2(\Theta
)\cdot \delta \,{\rm det}\,\Lambda ,\,\,\,\,\,\,
\nu\le l\le k-1.
$$
 В частности
$$
\zeta_\nu M_{\nu+1}\le
\frac{\gamma_2(\Theta)}{\gamma_1(\Theta)} \cdot \zeta_{k-1}M_k.
$$
Учтем, что  $\zeta_{k-1}\le \psi (M_k)$. Из условия случая получаем, что
$$
\zeta_\nu \le \frac{\gamma_2(\Theta)}{\gamma_1(\Theta)} \cdot
\frac{\psi (M_k)M_k}{M_{\nu+1}}\le \frac{\gamma_2(\Theta)}{\gamma_1(\Theta)} \cdot
\psi\left(M_\nu^{\frac{h(\alpha(\Theta))}{g_2(\alpha(\Theta))}}\right)
M_\nu^{\frac{h(\alpha(\Theta))}{g_2(\alpha(\Theta))}-1}
.
$$
Поскольку
$$
\alpha (g_2(\alpha))^2 +(\alpha -2)g_2(\alpha) -(\alpha -1)^2 =0,
$$
то снова получаем (\ref{novoe}).

Теорема доказана.

 \section{ Неоднородные приближения. }\label{neopribl}
 
\subsection{Одномерные задачи.}
Для вещественных $\theta$ и $\alpha $
рассмотрим величину
$$
\lambda(\theta,\alpha ) = \liminf_{x\to \infty} |x|\cdot ||x\theta -\alpha ||
$$
(здесь, естественно, подразумевается, что $x$ принимает целые значения).  
Согласно классическому результату Г. Минковского (см. \cite{Cassil}, теорема II главы III)
для любых $\theta, \alpha $ выполняется
$$
\lambda(\theta,\alpha ) \le \frac{1}{4}.
$$
 Определим
$$
\lambda (\theta) =\lambda (\theta,0)
$$
и
$$
\mu(\theta ) = \sup_\alpha \lambda (\theta,\alpha),
$$
где супремум берется по всем $\alpha$, не представимым в виде
$a\theta +b$  с целыми $a,b$.

А.Я.Хинчин в работе \cite{HCH}
 доказал неравенство
\begin{equation}
\mu (\theta ) \le \frac{\sqrt{1-4\lambda^2(\theta)}}{4},
\label{meen}
\end{equation}
которое бывает достаточно точным. (Например, если $\theta$  эквивалентно числу
$[0;k,k,k,...]$, задаваемому  периодической непрерывной  дробью, и неполное частное $k$   либо равно единице, либо четно, 
то
в неравенстве (\ref{meen}) имеет место знак равенства.) 
 Из (\ref{meen}), например, следует, что
равенство
 $\mu (\theta) = 1/4$ возможно только для таких $\theta$,
разложение которых в непрерывную дробь  имеет сколь угодно большие неполные частные.

 В работе А.Я.Хинчина \cite{HCH}
 имеется ряд других результатов. Изучению множеств значений функции $\lambda(\theta,\alpha)$  и ей подобных, например, "односторонней" функции
$$
k(\theta, \alpha ) =
 \liminf_{x\to +\infty}x\cdot ||x\theta -\alpha ||,
$$
посвящены работы  Дж.В.С.Касселса \cite{CC}, Е.Барнса \cite{Bar},  Т.Кузика и А.Поллингтона \cite{CUSP} и другие.

Сейчас мы сформулируем
 простой одномерный (но тем не менее весьма нетривиальный) результат, который содержится все в той же работе А.Я Хинчина
\cite{HINS}  (Satz 4).

\begin{theorem}\label{T12}

 Существует абсолютная константа $\gamma$  со следующим свойством:
для каждого действительного $\theta$ можно найти действительное $\alpha$ такое, что для всех натуральных $x$ выполнено
$$
||x\theta -\alpha ||\ge \frac{\gamma}{x}.
$$

\end{theorem}

Отметим, что в книге Касселса  \cite{Cassil}    результат  теоремы \ref{T12} имеется с постоянной  $\gamma = 1/51 $ (см. теорему XI  главы V  \cite{Cassil}; имеется непонятное замечание ко главе  V, в котором речь идет об оценке
точного значания постоянной $\gamma$;
по-видимому, наиболее точный результат принадлежит Г.Годвину \cite{Maloi}).

 В дальнейшем мы обсудим и многомерные обобщения теоремы \ref{T12}, и   моменты ее доказательства, поскольку
доказательство теоремы \ref{T12}  послужило отправной точкой для всех последующих многомерных результатов, в том числе, для теорем В. Ярника, о которых пойдет речь в следующем пункте.
 В дальнейшем мы поговорим как и о многомерных обобщениях  теоремы ??, так и о  моментах  доказательств.

\subsection{Многомерные теоремы}\label{mnt}

В доказательство теоремы \ref{T12} обнаружилась связь между однородными и неоднородными задачами 
в теории линейных диофантовых приближений. Эта связь получила дальнейшее развитие,
причем оказалось что характеристики синуглярности соотвествующих однородных систем играют решающую роль.
А.Я. Хинчин неоднократно возвращался к исследованию неоднородных линейных приближений.  
 В работе 1936 года  \cite{HAA}
он доказал основополагающий многомерный  результат (точнее,  утверждение было сформулировано для произвольного $m$,  а доказано при $m=2$), который мы здесь приводим.

\begin{theorem}\label{T13}{\rm (А.Я.Хинчин \cite{HAA}) }

Следующие два  условия {\rm  (i)}  и {\rm (ii)} эквивалентны.

{\rm  (i)}\,
Для набора вещественных чисел $\theta_1,...,\theta_m$  выполняется
$$
\limsup_{t\to +\infty}
\,\, t^m\cdot
\min_{{\bf x} \in \mathbb{Z}^m: 0<M({\bf x}) \le t}\,\,
||\sum_{1\le i \le m} \theta_i x_i || >0.
$$

{\rm  (ii)}  \, 
 Набор вещественных чисел  $\theta_1,...,\theta_m$
таков, что для любого  вещественного числа
$\alpha$
выполнено
$$
\liminf_{t\to +\infty}
\,\,
t^m\cdot 
\min_{{\bf x} \in \mathbb{Z}^m: M({\bf x}) \le t}\,\,
||\sum_{1\le i \le m} \theta_i x_i -\alpha ||<+\infty.
$$
\end{theorem}

Видим, что условие (i)  эквивалентно регулярности набора
$
\Theta  
$
(\ref{1f})
  с заданным значением $m$ и с $n=1$.

Xочется отметить, что доказательства всех дальнейших 
результатов, обсуждаемые в настоящем пункте, развивают идеи работ А.Я.Хинчина \cite{HINS},\cite{HAA}.

В 1948 году в работе \cite{HRSU} А.Я.Хинчин  доказал общее утверждение, формулировку которого мы тоже хотим здесь привести.

Система $\Theta$ называется {\it системой  Чебышева} если для любого набора вещественных чисел ${\bf \alpha } = (\alpha_1,...,\alpha_n)\in
\mathbb{R}^n$ найдется постоянная $\Gamma =\Gamma ({\bf \alpha})$ такая, что система диофантовых неравенств
$$
 \max_{1\le j\le n}||L_j({\bf x})-\alpha_j||\le {\Gamma} \cdot \left(\max_{1\le i\le m}|x_i|\right)^{\frac{m}{n}}
$$
имеет решения со сколь угодно большим значением величины $\max_{1\le i\le m}|x_i|$.

 Основной результат работы \cite{HRSU} состоит в следующем.

\begin{theorem}\label{T14}

Система $ \Theta$ является регулярной тогда и только тогда, когда она является
 системой Чебышева.
\end{theorem}

Теорема \ref{T13}, сформулированная выше, таким образом, является частным случаем теоремы \ref{T14} ($n=1$).

 Доказательство этой теоремы помимо оригинальной работы А.Я.Хинчина \cite{HINS} изложено, например, в  главе X книги \cite{Cassil}.

Следует отметить, что несколько ранее (в работах \cite{J39} и \cite{J41} 1939 и 1941 годов, соответсвтенно; см. также \cite{J41bis})
В. Ярник получил более общие утверждения, формулировки которых мы приводим ниже.

Наряду с "однородной" функцией Ярника (\ref{fj}) удобно использовать "неоднородную" функцию 
$$
\psi_{\Theta,\alpha } (t) =
\,\,\,\,\,\min_{{\bf x}\in \mathbb{R}^m: \, M({\bf x})\le t}
\,\,\,\,\,\max_{1\le j\le n}||L_j({\bf x})-\alpha_j||,\,\,\,\,
\alpha = (\alpha_1,...,\alpha_n) \in \mathbb{R}^n
$$
 и функцию
$$
\Psi_\Theta^{[inhom]}
= \sup_{\alpha \in [0,1]^n} \psi_{\Theta,\alpha }(t).
$$
 В определении "неоднородной" функции
Ярника можно брать минимум по всем {\it  ненулевым} векторам ${\bf x}$.
Результат от этого не изменится.

 Ясно, что функция Ярника
(\ref{fj})  удовлетворяет равенству
$$
\psi_\Theta (t) =
\psi_{\Theta, {\bf 0}}(t).
$$
Также наряду с системой  чисел $\Theta$
 будем рассмариваль транспонированную систему  $^t\!\Theta$. Для функции Ярника
$\psi_{\, ^t\!\Theta
}(t)$,  очевидно, имеем равенство
$$
\psi_{ \, ^t\!\Theta
}(t)= \min_{{\bf x} \in \mathbb{R}^n:  M(x) < t}
\max_{1\le i \le m}
||L_i^*({\bf x})||
,$$
 где
$$
L_i^*({\bf x}) =\sum_{j=1}^n\theta_{i,j} x_i
.
$$                                                                                  

 Всюду ниже в этом пункте предполагаем, что
функция $\psi (t)$ 
строго монотонно убывает к нулю. 
Обозначим через $\rho (t) $  функцию обратную к функции
$t\mapsto \frac{1}{\psi (t)}$ .

Сначала мы сформулируем более простые результаты работы \cite{J39} (теорему 1 из \cite{J39}).

\begin{theorem}\label{T15} {\rm  (В.Ярник \cite{J39})}

Пусть функция $\psi (t)$  
такова, что при некотором положительном  $\eta$  функция
$t\mapsto \frac{1}{t^\eta\psi (t)}$  монотонно стремится к бесконечности при $ t\to +\infty. $
Предположим, что при всех достаточно больших значениях $t$ выполняется
\begin{equation}
\psi_{\,^t\!\Theta (t)} > \psi (t).
\label{suff}
\end{equation}
   Тогда  для всех достаточно больших значений $t$  выполнено
\begin{equation}
\Psi_\Theta^{[inhom]}
\le\frac{((m+n)!(m+n))^{\frac{\eta+1}{\eta}}}{\rho (t)}.
\label{suffend}
\end{equation}
\end{theorem}

При  $\psi_{\,^t\!\Theta}= c \cdot t^{n/m}$ c  постоянной $c >0$  из теоремы \ref{T15}
получаем

{\bf Следствие.}\,\,\,{\it
 Если  транспонированная система $^t\!\Theta$ является регулярной,
то система  $\Theta$   является  системой Чебышева.
}

Следующая теорема фактически есть результат работы В. Ярника \cite{J41} (теорема 7 из \cite{J41}, см. также \cite{J41bis} -- они отличаются от формулируемой ниже теоремы \ref{T16} только тем, что в них не выписаны в явном виде констатны).

\begin{theorem}\label{T16}

 Пусть при всех достаточно больших $t$  выполнено
\begin{equation}
\psi_{\,^t\!\Theta} (t) \le \psi (t)
\label{ness}
\end{equation}
 Тогда 
 найдется набор вещественных чисел
 $(\alpha_1,...,\alpha_n)$  такой, что при всех достаточно больших $t$  выполняется
\begin{equation}
\psi_{\Theta,\alpha} (t) \ge \frac{1}{24n^{3/2}\rho (8mt)} .
\label{nessres}
\end{equation}
\end{theorem}

 В рассматриваемом вопросе это утверждение предсталяет основную сложность. Оно восходит к теореме А.Я.Хинчина \ref{T12},
 которая была доказана в 
\cite{HINS}. Доказательство именно этой теоремы мы прокомментируем в следующем пункте.
Здесь же мы отметим два важных следствия.

Полагая 
 $\psi_{\,^t\!\Theta}
 = o(t^{-n/m}),\,\, t \to+\infty $,
 получаем

{\bf Следствие.}\,\,\,{\it
Если  транспонированная система   $^t\!\Theta$  сингулярна, то система  $\Theta$ не есть система  Чебышева. 
}

 Это следствие вместе со следствием теоремы \ref{T15} приводит к утверждению, которое естественно считать доказанным В. Ярником:

\begin{theorem}\label{T17}
Транспонированная система $ ^t\!\Theta$ является регулярной тогда и только тогда, когда она является
 системой Чебышева.
\end{theorem}

  Для того, чтобы убедиться в эквивалентности теоремы \ref{T14} доказанной А.Я.Хинчиным и теоремы \ref{T17} доказанной В. Ярником,
надо вспоминть следующую теорему переноса из работы А.Я.Хинчина 1948 года \cite{DAN1948}
(она излагается также в книге
\cite{Cassil}, гл. V, теорема XII) .

\begin{theorem}\label{T17}
Система $\Theta$ сингулярна тогда и только тогда, когда сингулярна транспонированная система $^t\!\Theta$.
\end{theorem}

Но в работе В.Ярника 1941 года \cite{J41}  имеется следующее утверждение.

\begin{theorem}\label{T18}
  Если с некоторой положительной постоянной  $\gamma_1$  при всех достаточно больших значениях $t$
 выполнено
$$
\psi_\Theta
\ge \gamma_1 \cdot t^{-m/n},
$$
то
с некоторой положительной постоянной  $\gamma_2$ при всех достаточно больших значениях $t$
 выполнено
$$
\psi_{^t\!\Theta}
\ge \gamma_2 \cdot t^{-n/m},
$$
\end{theorem}

Теоремы \ref{T18} и \ref{T18} очевидным образом эквивалентны!

О других теоремах переноса пойдет речь  в пункте \ref{tippo}.

Совокупное доказательство теорем \ref{T15},\ref{T16} В.Ярника мало чем отличается от доказательсвта теоремы \ref{T17}, данного А.Я.Хинчиным.   Хинчин доказывает сначала, что если система  $\Theta$ является регулярной, то она есть система Чебышева. Для этого достаточно применить теорему Минковского к системе $m+n$  линейных форм от  $m+n+1$
 переменных $x_1,...,x_n,y_1,...,y_m, u$ вида
 $$
L_j ({\bf x}) - y_j - \xi_j u,\,\,\,\, j = 1,...,n,\,\,\,
x_j,\,\,\,\, j = 1,...,m, u
$$
(здесь  $\xi_j$ суть новые вещественные параметры).
 Затем,  с помощью применения конструктивных соображений (о них пойдет речь в пункте \ref{ododo}) и соображений переноса,
он доказывает обратное утверждение. Ярник же соображения переноса применяет в теореме \ref{T15}, а конструктивной теоремой является теорема \ref{T16}.

Теперь  о втором следствии.
Согласно теорееме Дирихле, упоминавшейся в самом начале настоящей статьи,
если положить
 $\psi (t) 
 = t^{-n/m} $,
 то условие выполнения неравенства (\ref{ness}) из  теоремы В.Ярника \ref{T16} становится пустым условием.
Теперь в качестве следствия мы получаем утверждение, имеющееся в книге Дж.В.С.Касселса (теорема  X главы V  из \cite{Cassil}):

\begin{theorem}\label{T19}

Для произвольных натуральных
$n,m$
существует положительная постоянная $\Gamma_{m,n} $,
обладающая следующим свойством. Для любой стстемы $\Theta$
 найдется неоднородность
   $(\alpha_1,...,\alpha_n)$ такая, что
 $$\inf_{{\bf x} \in \mathbb{Z}^m\setminus\{{\bf 0}\}}
\left( \max_{1\le j\le n}||L_j({\bf x}) -\alpha_j ||\right)^n \left(\max_{1\le i \le m} |x_i|\right)^m >\Gamma_{m,n}.
$$
\end{theorem}

Эта теорема является непосредственным обобщением теоремы \ref{T12} , доказанной А.Я.Хинчиным в \cite{HINS}.
С ней связана интересная история, о которой мы расскажем в пункте \ref{oteore}.

В заключение  настоящего пункта отметим, что в замечательной работе В.Ярника
  \cite{J41} имеется еще несколько теорем, о которых мы здесь не упомянули (в частности, там имеется результат метрического характера).

\subsection{О доказательстве теоремы \ref{T16}}\label{ododo}

 Доказательство теоремы \ref{T19}, имеющееся в книге \cite{Cassil} непосредственно переделывается в доказательство теоремы \ref{T16}. Оригинальные рассужждения В. Ярника несколько более громоздки.   Мы приведем здесь схему доказательства теоремы \ref{T16}, следуя книге \cite{Cassil}.
Приводя схему доказательства, мы ограничимся случаем, когда все числа  $\theta^{i}_{j}$  линейно независимы вместе с единицей над  $\mathbb{Z}$ (исключительно для удобства изложения, чтобы избежать разбора нескольких случаев).

{\bf 1.}  Последовательность наилучших приближений
 ${\bf w}_\nu =  (u_{1,\nu},...,u_{n,\nu}, v_{1,\nu},...,v_{m,\nu})$  для транспониро-ванной системы  $^t\!\Theta$ 
 надо "прорядить" таким образом, чтобы получилась последовательность
${\bf w}_{\nu_k}$  такая, что
\begin{equation}
\max_{1\le i \le m} ||L^*({\bf u}_{\nu_k})|| = \psi_{^t\!\Theta} \left(\frac{M({\bf u}_{\nu_{k+1}})}{3\sqrt{n}}\right),
\,\,\,\,\,\,\, M({\bf u}_{\nu+1})\ge 3\sqrt{n} M({\bf u}_\nu),\,\,\,\,\,\,\, M ({\bf u}_\nu ) =\max_{1\le j \le n} |u_{j, \nu+1}|
\label{redk}
\end{equation}
(см. лемму 4 \S 6  главы V \cite{Cassil} или рассуждения из  \S 2, пункт 2 из \cite{HRSU}, или конструкцию
леммы 1 из \cite{J41}). 

{\bf 2.} Необходимо следующее утверждение о диофантовых приближениях с лакунарной последова-тельностью векторов.
Приводимая ниже формулировка взята из  \cite{Cassil} (лемма 2 \S 6  глава V).  Она непосредственно обобщает лемму А.Я. Хинчина
(Hilfssatz 3
 из \cite{HINS}). Этому утверждению и его развитию посвящен пункт \ref{lakey} добавления в настоящем обзоре.

\begin{lem}\label{YT2}
  Для последовательности векторов  ${
\bf u}_k =(u_{1,k},....,u_{n,k}) \in \mathbb{R}^n$,
 удовлетворяющей  неравенствам
$$
\max_{1\le j \le n} |u_{j,k+1}| \ge  3\sqrt{n}\cdot
\max_{1\le j \le n} |u_{j,k}|,
$$
 найдутся вещественные числа  $\alpha_1,...\alpha_n$  такие, что
$$
||
\sum_{1\le j \le n } u_{j,k}\alpha_j||\ge  \frac{1}{4}\,\,\,\,\, k  = 1,2,3,... .
$$
\end{lem}

Наличие дополнительного (по сравнению с  соответствующими утверждениями из \cite{Cassil}) множителя
 $\sqrt{n}$  связано с тем, что мы рассматриваем sup-норму, а не евклидову норму.

{\bf 3.} Далее доказывается, что числа $\alpha_j$
и есть те самые неоднородности, существование которых утверждается в теореме \ref{T16}.
Для этого используется тождество
$$
\sum_{1\le j \le n } u_{j}\alpha_j
=
\sum_{1\le j \le n } u_j(\alpha_j -L_j({\bf x})
) +
\sum_{1\le j \le n } u_{j}L_j ({\bf x})=
\sum_{1\le j \le n } u_j(\alpha_j -L_j({\bf x})
) +
\sum_{1\le i \le m } x_i L_i^* ({\bf u}),
$$
где, естественно,
$$
 {\bf x}  = (x_1,...,x_m) \in \mathbb{R}^m\,\,\,\,\,
{\bf u}= (u_1,...,u_n) \in \mathbb{R}^n.
$$
В этом тождестве
(примененном для вектора  ${\bf u}_{\nu_k}$,  определенного в п.1 доказательства,   чисел  $\alpha_j$  из леммы  \ref{YT2} и произвольного вектора  ${\bf x} \in \mathbb{R}^m$)  
 переходим к модулям и используем лемму \ref{YT2}:
$$
\frac{1}{4} \le
n\cdot M({\bf u}_{\nu_k})\cdot
\max_{1\le j\le n} || L_j({\bf x}) -\alpha_j ||+m\cdot
\max_{1\le i \le m } |x_i| \cdot \psi_{^t\!\Theta} 
 \left(\frac{M({\bf u}_{\nu_{k+1}})}{3\sqrt{n}}\right).
$$
 Теперь
  если по заданному вектору  ${\bf x}$
выбрать  $k$  из условия
\begin{equation}
 \psi_{^t\!\Theta} 
 \left(\frac{M({\bf u}_{\nu_k})}{3\sqrt{n}}\right)
\ge
\frac{1}{8m \max_{1\le i \le m } |x_i| }\ge
 \psi_{^t\!\Theta} 
 \left(\frac{M({\bf u}_{\nu_{k+1}})}{3\sqrt{n}}\right)
,
\label{usl}
\end{equation}
 то получается
\begin{equation}
\max_{1\le j\le n} || L_j({\bf x}) -\alpha_j ||
\ge\frac{1}{8nM({\bf u}_{\nu_k})}.
\label{us}
\end{equation}
Но из условия теоремы (\ref{ness}) и соотношения (\ref{usl}) получаем
$$
\psi\left(
\frac{M({\bf u}_{\nu_k})}
{3\sqrt{n}}
\right)
\ge
\psi_{^t\!\Theta }\left(
\frac{M({\bf u}_{\nu_k})}
{3\sqrt{n}}
\right)
\ge \frac{1}{8m\max_{1\le i \le m } |x_i| }.
$$
Отсюда, по определению функции  $\rho (\cdot )$,
получаем
$$
\rho (8m\max_{1\le i \le m } |x_i|) \ge
\frac{M({\bf u}_{\nu_k})}{3\sqrt{n}}
.
$$
 Подстановка последнего неравенства в неравенство  (\ref{us}) дает именно (\ref{nessres})
 с $ t = \max_{1\le i \le m } |x_i|$.

\subsection{  О теореме \ref{T19}}
\label{oteore}

В последнее время к линейным неоднородным диофантовым приближениям вновь пробудился интерес.
В частности, теорема \ref{T19} уточнялась несколькими авторами. Здесь мы изложим историю вопроса и 
изложим схему доказательства наиболее сильного на настоящее время результата.

 Следующую теорему доказал Д.Клейнбок \cite{KL}.

\begin{theorem}\label{T21}

Через
  ${\cal B}$
обозначим множество вещественных матриц размера
 $(m+1)\times n$ 
вида
$$(\theta_{i,j},\eta_j ),\,\,\, 1\le i \le m,\,\,\, 1\le j\le n;\,\,\,\,
$$
для которых выполняется неравенство
$$\inf_{{\bf x} \in \mathbb{Z}^m\setminus\{{\bf 0}\}}
\left( \max_{1\le j\le n}||L_j({\bf x}) -\eta_j ||\right)^n \left(\max_{1\le i \le m} |x_i|\right)^m >0.
$$ 

Тогда  ${\cal B}$
есть множество
полной хаусдорфовой размерности в 
$\mathbb{R}^{mn+n}$.
\end{theorem}

Доказательство Д. Клейнбока было связано с рассмотрением специальных потоков на некоторых однородных пространствах.
Я.Бюжо, С.Харрап, 
С.Кристенсен и С.Велани \cite{BU} дали простое доказательство результата Д.Клейнбока. Более того, они получили более сильную теорему, которую мы сейчас и сформулируем.

\begin{theorem}\label{T22}
 \,\, 
Для произвольного набора вещественных чисел
 $$
\Theta =\{
\theta_{i,j},\,\,\, 1\le i \le m,\,\,\, 1\le j\le n
\}$$
рассмотрим множесттво ${\cal B}(\Theta )$, состоящее из векторов 
$$
(\eta_1,...,\eta_n) $$
таких, что 
$$\inf_{{\bf x} \in \mathbb{Z}^m\setminus\{{\bf 0}\}}
\left( \max_{1\le j\le n}||L_j({\bf x}) -\eta_j ||\right)^n \left(\max_{1\le i \le m} |x_i|\right)^m >0.
$$
Тогда
${\cal B}(\Theta)$
есть множество полной размерности Хаусдорфа в
$\mathbb{R}^{n}$.
\end{theorem}

Для формулировки и обсуждения дальнейших результатов необходимо знание
теории $(\alpha,\beta)$-игр В.М.Шмидта.
Основные понятия и факты из этой теории, необходимые для понимания  формулировок трех последующих теорем и дальнейших доказательств из текущего пункта,  мы изложили в пункте \ref{uinni} добавления.

Сейчас мы формулируем примечательный результат, который  был недавно получен Дж. Тсенгом \cite{Ts}.

\begin{theorem}\label{T23}
 \,\,\,
Для любого вещественного числа
$\theta$ множество  ${\cal B}$, состоящее из вещественных чисел $\eta$  таких, что
$$
\inf_{x\in \mathbb{Z}\setminus\{{0}\}} |x|\cdot ||\theta x+\eta||>0
$$
является
 $\alpha$-выигрышным множеством для любого $\alpha \in (0,1/8)$.
\end{theorem}

Согласно замечанию 2.3 из работы Дж.Тсенга \cite{Ts},
Дж.Тсенг и М. Айсиедлер получили обобщение теоремы \ref{T23} на случай систем линейных форм произвольной размерности
(Препринт \cite{tsengA} появился совсем недавно).
 
В работе \cite{MaR} автор, модифицируя доказательство теоремы  \ref{T19}
получил следующее утверждение.

\begin{theorem}\label{T25}
 Пусть $\alpha \in (0,1/2)$. Тогда для любого набора вещественных чисел
$$
\Theta =\{
\theta_{i,j},\,\,\, 1\le i \le m,\,\,\, 1\le j\le n
\}$$
множество ${\cal B}(\Theta )$, состоящее из векторов
$$
(\eta_1,...,\eta_n) $$
таких, что
$$\inf_{{\bf x} \in \mathbb{Z}^m\setminus\{{\bf 0}\}}
\left( \max_{1\le j\le n}||L_j({\bf x}) -\eta_j ||\right)^n\left(\max_{1\le i \le m} |x_i|\right)^m >0
$$
является  $\alpha$-выигрышным множеством в $\mathbb{R}^n$.
\end{theorem}

Справедлив более общий результат:

\begin{theorem}\label{T26}
Пусть  $\alpha \in (0,1/2)$. 
Путсь функция $\psi (t) $
строго монотонно убывает к нулю при $t\to+\infty$  и  пусть  $\rho (t)$ есть функция, обратная к функции
 $t\mapsto 1/\psi (t)$.

Для любого набора вещественных чисел
$
\Theta $
 рассмотрим функцию Ярника  $\psi_{^t\!\Theta}$.  Пусть
при всех достаточно больших $t$  выполнено
 $$
\psi_{^t\!\Theta}\le \psi (t).$$
 
Тогда 
множество ${\cal B}(\Theta )$, состоящее из векторов
$$
(\eta_1,...,\eta_n) $$
таких, что
$$\inf_{{\bf x} \in \mathbb{Z}^m\setminus\{{\bf 0}\}}  
\left( \max_{1\le j\le n}||L_j({\bf x}) -\eta_j ||\right)\cdot \rho \left(\max_{1\le i\le m}||x_i|\right) >0
$$
является  $\alpha$-выигрышным множеством в $\mathbb{R}^n$.

\end{theorem}

Легко видеть 
(т.к.   всегда $\psi_{^t\!\Theta (t)} \le t^{-n/m}$), что сформулированная теорема \ref{T25}
является частным случаем теоремы  \ref{T26}.

Ниже в этом пункте мы  приводим схему доказательства последней теоремы.

Чтобы получить доказательство теоремы \ref{T26} надо всего лишь переформулировать лемму \ref{YT2} о лакунарной последовательности векторов из пункта \ref{ododo}.

\begin{lem}\label{HHT2}
 Пусть последовательность вещественных чисел $t_r, r =1,2,3,...$
удовлетворяет условию лакунарности
\begin{equation}
\frac{t_{r+1}}{t_r} \ge M
 ,\,\,\,\,\,
r=1,2,3,...
 \label{lac}
\end{equation}
для некоторого $ M>1$. 
пусть последовательность
$\Lambda\subset \mathbb{Z}^n$ 
целочисленных векторов ${\bf u}^{(r)}= (u^{(r)}_1, ....,u^{(r)}_n  ) \in \mathbb{Z}^n$
такова, что
\begin{equation}
t_r^2= (u^{(r)}_1)^2+ ....+(u^{(r)}_n)^2. \label{vec}
\end{equation}
Тогда множество
$$
N(\Lambda ) =\{ \eta =(\eta_1,...,\eta_n)\in\mathbb{R}^n:\,\,\, \exists \, c(\eta)>0 \,\,\,\text{такое, что}\,\,\, || u^{(r)}_1\eta_1+ ...+u^{(r)}_n
\eta_n||\ge c(\eta )\,\,\ \forall r \in \mathbb{N} \}
$$
является
$\alpha$-выигрышными при каждом  $\alpha \in (0,1/2)$.
\end{lem}

То, что в условиях теоремы\ref{T26}  всякий набор   $\eta \in N(\Lambda )$  
(при надлежащем выботе векторов ${\bf u}^{(r)}$)
будет принадлежать множеству ${\cal B} (\Theta)$,
сразу 
следует из классических аргуметнов доказательства теоремы \ref{T16}, которые мы напомнили читателям в пункте  \label{ododo}. Доказательства частного случая (теоремы  \ref{T25})  и общего случая (теоремы \ref{T26})   не отличаются.
 
  Ниже приводим схему доказательства 
 леммы \ref{HHT2}. 

Как обычно, для 
$\alpha, \beta \in (0,1)$ полагаем  $\gamma = 1+\alpha\beta -2\alpha >0$.

Рассмотрим  шар  (в евклидовой норме; в этом доказательстве всюду удобнее использовать  именно евклидову норму)
$B\subset \mathbb{R}^n$с центром  $O$ и радиусом $\rho$.
Его границу будем обозначать
$S = \partial B$. 
Через  $\mu $ 
обозначаем нормализованную меру Лебега на
$S$ (так что $ \int_S d\mu = \mu S = 1$).

Пусть 
 $ x \in S$.
Определим $n-1$-мерное аффинное подпространство
 $\pi (x)\subset \mathbb{R}^n$, как подпространство, проходящее через точку
$O$ и ортогональное одномерному подпространству, проходящему через две точки  $O$  и $x$. 
Определим "полупространство"
$\Pi (x)$ с границей  $\pi (x)$ и такое, что $x \in \Pi (x)$.

Для данных $\alpha, \beta \in (0,1)$ 
рассмортим полупространство
$ \Pi_{\alpha, \beta,\rho }(x)$ такое, что
$ \Pi_{\alpha, \beta,\rho }(x)\subset \Pi (x)$ и расстояние от   $ \Pi_{\alpha, \beta,\rho }(x)$
 до 
 $O$  равно $\frac{\gamma \rho}{2}$.  Положим
 $$\Omega (x) =S\cap \Pi_{\alpha, \beta,\rho }(x),\,\,\,
 \Omega^*(x) = \bigcup_{y\in S:\,\, \Pi (y) \supset \Omega (x)} \{y\}.
 $$
Ясно, что мера
 $\mu \Omega^* (x)$ не зависит от  $x\in S$.
 Положим
\begin{equation}
 \omega = \omega (\alpha, \beta ) = \mu \Omega^* (x) \in (0,1).
 \label{om}
\end{equation}

С помощью соображений усреднения сразу получаем следующий вспомогательный результат.

\begin{lem}\label{HHHT2}
  Рассмотрим произвольные аффинные подпространства
 $\pi_1, ... , \pi_k$ размерности $n-1$.
Тогда найдется точка $x\in S$
 такая, что
 $$\Omega (x) \cap \pi_j = \varnothing,\,\,\,
 $$
для по крайней мере 
$  \lceil\omega k \rceil$  штук индексов $j$.
\end{lem}
 
 Следующая лемма принадлежит В.М.Шмидту (лемма 1B главы 3 из книги \cite{SCH}).

\begin{lem}\label{HHHHT2}
 Пусть  $t$ таково, что
 $$
 (\alpha \beta )^t <\frac{\gamma}{2}.
 $$
Предположим, что  в игре встетился шар черных $ B_j$.
Предположим, что $n-1$-мерное аффинное подпространство  $\pi$ 
проходит через центр шара  $B_j$. 
Тогда белые могут играть таким образом, что в результате шар $B_{k+t}$ 
будет целиком содержаться в "полупространстве"  $ \Pi_{\alpha, \beta,\rho_j }(x)$,
граница которого параллельна подпространству $\pi$.
 
\end{lem}

Выбрав параметры
\begin{equation}
t = t(\alpha,\beta) =\left\lceil  \frac{\log(\gamma/2)}{\log (\alpha\beta )}\right\rceil,\,\,\,\,\, \tau_k =
 t\cdot\left\lceil \frac{\log k}{\log\left(\frac{1}{1-\omega}\right)} \right\rceil
, \label{tau}
\end{equation}
где  $\omega$ определено в  (\ref{om}), получаем

{\bf Следствие.}\,\,\,{\it
Пусть  шар $B_j$ с радиусом $\rho_j$
всретился в игре черных. Пусть имеется набор аффинных подпространств
 $\pi_i,\, 1\le i \le k$. 
Тогда белые могут играть таким образом, что для каждой точки $x \in B_{j+\tau_k}$ 
расстояние от точки  $x$  до каждого из подпространств $\pi_i,\, 1\le i \le k$  будет больше, чем $\frac{\rho_{j+\tau_k}\gamma}{2}$.}

Далее выбираем параметр 
 $k=k (\alpha ,\beta,  M)$ такой, что
\begin{equation}
\tau_k\times \frac{\log (1/(\alpha\beta ) )}{\log M} +2 < k.
  \label{ka}
\end{equation}
Не ограничая общности, мы можем предполагать, что в дополнение к условию лакунарности
(\ref{lac}) выполняется
$
\frac{t_{r+1}}{t_r} \le M^2, \,\,\, r=1,2,3,... .
 $
Выбор параметров (\ref{tau},\ref{ka}) позволяет белым вести игру таким образом что за время игры 
$\tau_k$  удавалось бы "убегать" от $q=r_{j-1}-r_j<k$
 семейств "опасных" подпространств вида
 $$
 \{ {\bf y}=(y_1,... , y_n)\in \mathbb{R}^n:\,\, u_1^{(r)}y_1+...+ u_n^{(r)}y_n =a\},\,\,\, a\in \mathbb{Z}
  ,\,\,\,\, r_j\le r<r_{j+1}.
$$

\section{Пространства решеток}\label{prore}
Задачи о диофантовых приближениях естественным образом связаны с вопросами поведения семейств (орбит)
некоторых решеток. На эту тему имеется много работ (см., например, работы Г.А.Маргулиса и Д.Клейнбока \cite{kma},
\cite{kma1},\cite{KL},\cite{KLEI},\cite{kdir}, обзор А.Городника \cite{oogoro} и литературу, цитируемую в этих работах).  Ниже мы остановимся на двух 
задачах, имеющих отношение к сингулярным матрицам.

\subsection{ Метрическая теорема Даверпорта-Шмидта.}\label{daves}

Как отмечалось в конце пункта \ref{tees},  сингулярные системы $\Theta$ образуют в $nm$-мерном пространстве множество нулевой меры Лебега. То же самое можно сказать так.
Рассмотрим множество
$\hbox{\got T}_\mu \subset \mathbb{R}^{mn}$,
состоящее из тех наборов $\Theta$, для которых
при всех достаточно больших значениях параметра $t$
  в каждой области  
\begin{equation}
\{ {\bf x} \in \mathbb{R}^n:\,\,\,
 \max_{1\le j\le n}||L_j({\bf x})||\le \frac{\mu  }{t},\,\,\,\,\,\,\,\, 0<\max_{1\le i\le m}|x_i|\le \mu \cdot t^{\frac{n}{m}}
\}
\label{regul1}
\end{equation}
имется  ненулевая целая точка.
(Ясно, что $\hbox{\got T}_{\mu_1} \subset \hbox{\got T}_{\mu_2}$ при $\mu_1 <\mu_2$.)
Тогда пересечение
$$
\bigcap_{\mu >0}\hbox{\got T}_\mu
$$
имеет лебегову меру нуль.

 Замечательное усиление этого факта было получено Г.Давенпортом и В.М.Шмидтом в работе \cite{DSdir0}, \cite{DSdir1}.

\begin{theorem}\label{T90} {\rm 
(Г.Давенпорт, В.М.Шмидт  \cite{DSdir0}, \cite{DSdir1})}

Для любого $\mu$ из интервала $0<\mu <1$ множество
$\hbox{\got  T}_\mu$ имеет меру нуль.
\end{theorem}

Отметим, что согласно теоремы Дирихле,
упоминавшейся в самом начале нашей работы,
при $\mu =1$  имеем $\hbox{\got T}_1 = \mathbb{R}^{mn}$.

Замечание. На самом деле, в работах  рассмотрен только случай $n=1$ (одна линейная форма) или $m=1$
(совместные приближения); тем не менее, доказательства проходят в общем случае. В случае же когда и $m=1$ и $n=1$
Г.Давенпорт и В.М.Шмидт установили, что из условия 
$\theta \in \hbox{\got T}_\mu$ с некоторым $\mu <1$ вытекает, что $\theta$ является плохо приближаемым числом,
то есть
$$
\liminf _{q\to \infty} q||q\theta || >0
$$
 (в последней формуле, естественно, предполагается
что $q$ принимает целые значения),
или, что то  же самое,
неполные частные разложения $\theta$ в цепную дробь ограничены.
 Фактически этот результат связан с применением формулы  (\ref{approx11}). 

Доказательство теоремы \ref{T90} 
в случае $n=1$
 фактически основывается на следующем утверждении.

\begin{theorem}\label{T11} {\rm  (В.М.Шмидт \cite{schlat})}

Рассмотрим последовательность натуральных чисел $N_\nu$, возрастающих к бесконечности.
Для набора вещественных чисел $\Theta =(\theta^1,....,\theta^m)\in \mathbb{R}^m$ рассмотрим решетку
$$
\Lambda (\Theta, N)= {\cal A}(\Theta, N)\mathbb{Z}^{m+1}
,$$
где матрица ${\cal A}(\Theta, N)$ имеет вид
$$
{\cal A}(\Theta, N)=
\left(
\begin{array}{ccccc}
N^{-1} & 0& 0&  \cdots & 0
\cr
 0  & N^{-1} & 0&\cdots & 0
\cr
 0   &0& N^{-1} &  \cdots
& 0
 \cr
 \cdots &\cdots &\cdots &\cdots
 \cr
 \theta^1N^m  & \theta^2N^m& \theta^3N^m&\cdots & N^m
\end{array}\right)
.$$
Тогда для почти всех (в смысле меры Лебега) наборов $\Theta \in \mathbb{R}^m$
последовательность решеток
$$
\Lambda(\Theta, N_\nu),\,\,\,\, \nu = 1,2,3,...
$$
всюду плотна в пространстве решеток в $\mathbb{R}^{m+1}$ с определителем 1.
\end{theorem}

Случай $m=1$,
связанный с рассмотрением решеток вида 
\begin{equation}\label{lati}
\Lambda^* (\Theta, N) = {\cal A}^* (\Theta, N)\mathbb{Z}^{n+1},\,\,\,\,
{\cal A}^* (\Theta, N) = \left(
\begin{array}{ccccc}
N^{-1} & 0& 0&  \cdots &0 \cr N^{\frac{1}{n}} \theta_1 & N^{\frac{1}{n}} & 0&\cdots & 0 \cr  N^{\frac{1}{n}} \theta_2 &0& N^{\frac{1}{n}} &  \cdots
& 0 \cr \cdots &\cdots &\cdots &\cdots \cr N^{\frac{1}{n}} \theta_n &0&0&\cdots & N^{\frac{1}{n}}
\end{array}\right),\,\,\,
\Theta =
\left(
\begin{array}{c}
\theta_1\cr
\theta_2\cr
\vdots
\cr
\theta_n
\end{array}\right),
\end{equation}
 в оригинальной работе \cite{DSdir1} получается из случая $n=1$
с помощью техники переноса.

В настоящей статье мы не будем вдаваться в подробности, касающиеся пространства решеток и вопросов сходимости в нем. Отметим, лишь что классические работы К. Малера \cite{Mah},\cite{Mah1}
на эту тему изложены в главе V монографии Дж.В. С. Касселса \cite{CGEOM}. Следует также сказать, что
ряд результатов теории диофантовых приближений  оказался связанным с
рассмотрением спецальных динамически систем на пространствах решеток и получил свое развитие в работах Д.Клейнбока. В частности,  в работе \cite{kdir} получено обобщение теоремы  Давенпорта-Шмидта, о которой идет речь в  настоящем параграфе.

\subsection{Задача о последовательных минимумах.}

Пусть дана решетка  $\Lambda \subset \mathbb{R}^d$  и
выпуклое ${\bf 0}$-симметричное тело  $\Omega \subset \mathbb{R}^d$.
 Величины
$$
\mu_l (\Omega, \Lambda) =\inf \{ t:\,
t\Omega \,\,
\text{содержит }
l
\text{ линейно независимых точек решетки}\,\,
\Lambda \},\,\,\
1\le l \le d
$$
 называются {\it  последовательными минимумами } решетки  $\Lambda$  относительно тела  $\Omega$. Знаменитая вторая теорема Минковского о выпуклом теле (см., например, главу VIII из \cite{CGEOM} или главу из IV \cite{SCH}) утверждает, что
$$
\frac{2^d}{d!}\cdot {\rm det }\Lambda \le
\mu_1 (\Omega, \Lambda )\cdots \mu_d (\Omega, \Lambda )
\cdot {\rm mes} \Omega \le 2^d{\rm det} \Lambda.
$$
В настоящем пункте мы ограничимся рассмотрением задачи о совместных приближениях, то есть случаем  $m=1$.
 Теорема Дирихле,  о линейных формах, упоминавшаяся в начале нашей работы, может быть переформулирована
в этом случае  следующим образом.
Рассмотрим куб  $E=[-1,1]^{n+1}\subset \mathbb{R}^{n+1}$ и
для набора  $\Theta $ 
 решетку $\Lambda^*(\Theta, N)$,
 определенную в (\ref{lati}).
 Тогда для любого вещественного  $N\ge 1$  выполняется
$$
\mu_1 (E,\Lambda^*(\Theta,N)) \le 1.
$$
Как отмечал В.М.Шмидт \cite{SCHL}, легко видеть, что для каждого  $k, 1\le k \le n$  найдется последовательность вещественных чисел
 $N_\nu$, стремящаяся к  бесконечности такая, что
$$
\mu_k  (E,\Lambda^*(\Theta,N_\nu)) =
\mu_{k+1}  (E,\Lambda^*(\Theta,N_\nu)) 
.
$$
 В частности, при  $n=1$  из теоремы Минковского получаем, что
$$
1\ll\mu_1(E,\Lambda^*(\Theta,N))\cdot \mu_2(E,\Lambda^*(\Theta,N))\ll 1
$$
и, следовательно, не может быть чтобы  $\lim_{N\to +\infty} \mu_1(E,\Lambda^*(\Theta,N)) =0$.
Но если  $n>1$ то даже для линейно независимых вместе с единицей чисел  $\theta_1,...,\theta_n$
 может оказаться, что 
$$
\lim_{N\to +\infty}
\mu_{n-1}(E,\Lambda^*(\Theta,N)) =0.
$$
(Набор $\Theta$ в этом случае следует брать из теоремы  \ref{T4},
с надлежащим выбором функции $\psi$.)

Автор в \cite{arxsuxx}  доказал следующую теорему, ответив тем самым на вопрос, поставленный В.М.Шмидтом \cite{SCHL}.

\begin{theorem}\label{suxx}

Пусть
$ 1\le k \le n-1$. 
Тогда найдется набор   $\Theta = (\theta_1,...,\theta_n)$
 вещественных чисел
 такой, что

$\bullet$\,\,   $ 1,\theta_1,...,\theta_n$   лиинейно независимы над  $\mathbb{Z}$;

$\bullet$\,\, $\mu_k (\xi, N) \to 0$  при $ N\to +\infty$;

$\bullet$\,\, $\mu_{k+2} (\xi, N) \to +\infty$ при  $ N\to +\infty$.

\end{theorem}

Доказательство этой теоремы в общем случае достаточно громоздко. В случае  $k=1$
 доказательство весьма просто. Оно идейно очень похоже на доказательство теоремы 
\ref{antilaga}.

Замечание. Легко видеть, что утверждение теоремы  \ref{suxx}  становится тривиальным, если не требовать выполнения условия
линейной независимости над 
$\mathbb{Z}$.

Отметим в заключение этого пункта, что недавно усиление теоремы \ref{suxx}
анонсировал И.Чеунг.

\section{ Теоремы переноса.}\label{tippo}

 Теоремы переноса -- это утверждения о том, что если система $\Theta$ обладает некоторыми диофантовыми свойствами,
то и транспонированная система $\Theta^*$ тоже обладает некоторыми диофантовыми свойствами.
  Простейшая теорема переноса, касающаяся сингулярных систем 
уже упоминалась нами в пункте \ref{mnt} (теорема \ref{T17}).
 Утверждения о связи однородной задачи для системы  $^t\!\Theta$  и неоднородной задачи для системы 
$\Theta$  (такие, как теоремы \ref{T14} - \ref{T19}) тоже можно интерпретировать как теоремы переноса.
Теории переноса посвящено огромное количество работ. Следующий список никоим образим не претендующий на полноту,
тем не менее дает некоторое представление о том, какие математики и в связи с какими вопросами занимались
переносом
(более подробную библиографию можно почерпнуть как раз из упоминаемых ниже работ).

$\bullet$  Классические результаты А.Я.Хинчина о диофантовых экспонентах для "обычных" диофантовых приближений 
\cite{H1925},\cite{HINS}; общий случай, рассмотренный Ф.Дайсоном \cite{dai};
точность оценок А.Я.Хинчина, доказанная В.Ярником \cite{jP59}.

$\bullet$  Теоремы переноса для сингулярных систем (работы В.Ярника \cite{JTBIL}, А.Апфельбека \cite{apf}, М.Лорана
и Я.Бюжо \cite{lora}, \cite{lora1}, \cite{blora},).

$\bullet$  Теоремы переноса для приближений рациональных чисел (работы  Н.М.Коробова \cite{koro},\cite{corobo}); интерес к такого рода задачам появился в связи с теорией приближенного интегрирования.

$\bullet$   Теормы переноса, так или иначе связанные с примененнием тригонометрических сумм. Они имеются, например,  у А.О.Гельфонда \cite{gelcas},  В.М.Шмидта и Ю.Вана \cite{ScVa},
Ю.В. Каширского \cite{Kash}.  Этот подход основывается на работе К.Л.Зигеля \cite{zz}.

$\bullet$  Более тонкие результаты типа пареноса, связанные с приближениями рациональными подпространствами
(работы В.М.Шмидта \cite{SCHgra}, М.Лорана \cite{lora}, М.Лорана
и Я.Бюжо \cite{blora}).

Классические результаты изложены в книгах Дж.В.С.Касселса \cite{Cassil} (глава V)
и В.М.Шмидта \cite{SCH}(глава IV).

В настоящей работе нас  наибоее интересуют  теоремы переноса для сингулярных систем,
то есть связанные с поведением величин типа
$$
\limsup_{t\to +\infty} \varphi (t) \psi_\Theta (t),
$$
где $\varphi(t)$  есть некоторая функция, а $\psi_\Theta$  есть функция Ярника  (\ref{fj}).
Это  результаты В.Ярника и А.Апфельбека  и теоремы, недавно полученные М.Лораном  и Я.Бюжо.

Тем не менее мы начнем с того, что приведем классические теоремы А.Я.Хинчина и Ф.Дайсона.

Прежде чем формулировать результаты, напомним определение  экспоненты $\alpha(\Theta)$    -- это супремум тех 
$\gamma$,
 для которых     выполнено 
 $$
\limsup_{t\to +\infty} t^\gamma
\psi_\Theta (t) <+\infty
$$
Здесь же 
напомним определение  экспоненты  $\beta(\Theta)$, которая есть супремум тех 
$\gamma$,
 для которых     выполнено 
 $$
\liminf_{t\to +\infty} t^\gamma
\psi_\Theta (t) <+\infty
$$
(эти  определения уже использовались у нас в статье в пункте \ref{sdty}).

\subsection{Теоремы А.Я.Хинчина и Ф.Дайсона}\label{HHDD}
 Теорему переноса для показателей  $\beta(\Theta)$  для системы
$\Theta =(\theta^1,...,\theta^m)\in \mathbb{R}^m$ и транспонированной системы 
$^t\!\Theta$  А.Я.Хинчин доказал все в той же знаменитой работе  \cite{HINS}.

\begin{theorem}\label{perenoshi}
Для набора $\Theta =(\theta^1,...,\theta^m)\in \mathbb{R}^m$    выполняется
\begin{equation}\label{hincini}
\frac{\beta(\Theta)}{(m-1)\beta(\theta) +m}\le \beta(^t\!\Theta) \le 
\frac{\beta(\theta)-m+1}{m}.
\end{equation}
\end{theorem}

Результат  Ф.Дайсона \cite{dai}  
(простое доказательство есть у А.Я.Хинчина в \cite{dobpere})
выглядит так:
 \begin{theorem}\label{perenosdai}
При произвольных значентях размерностей  $m,n$ для матрицы $\Theta$     выполняется
\begin{equation}\label{hincini}
\beta(^t\!\Theta)\ge
\frac{n\beta(\Theta)+n-1}{(m-1)\beta(\Theta) +m} .
\end{equation}
\end{theorem}

Естественно, теорему \ref{perenoshi} можно получить, применив два раза
теорему \ref{perenosdai} (для матрицы-строки и для матрицы-столбца).

\subsection{Результаты В.Ярника  и А.Апфельбека}\label{JAP}

Работа В.Ярника \cite{JTBIL} посвящена теоремам переноса, связывающим случаи  $m=1$  и  $n=1$.
Мы приведем общий результат (теоремы \ref{www1}, \ref{www2};
в оригинальной работе \cite{JTBIL} -- Satz 7 и  Satz 8) и перечислим некоторые его следствия.
Будем рассматривать вектор-строку
$$
\Theta = (\theta^1,...,\theta^m),\,\,m\ge 2
$$
и  функции
$$
\psi_\Theta (t),\,\,\,\psi_{^t\!\Theta}(t).
$$

\begin{theorem}\label{www1}
Пусть  числа $1,\theta^1,...,\theta^m$
 линейно независимы над  $\mathbb{Z}$.
Пусть $K$  есть некоторая положительная постоянная.
Пусть функция $\varphi (t)$ возрастает,
и функция  $t\mapsto  t\cdot (\varphi (t))^{-1}$  монотонно возрастает к бесконечности.
Пусть выполняется
$$
\limsup_{t\to+\infty}
\varphi (t) \psi_{^t\!\Theta}(t) <K.
$$
Тогда

{\rm (i)}  выполняется неравенство
$$
\limsup_{t\to +\infty} t^{m-1}\varphi (t^m) \cdot \psi_\Theta (t) \le m^{2m} K;
$$

{\rm (ii)}  если же  известно, что  $\varphi (t) \ge t^{\frac{m-1}{m}}$
 при всех достаточно больших значениях $t$  и
функция  $\rho (t)$  есть обратная к функции  $t\mapsto  t\cdot (\varphi (t))^{-1}$,
то выполнено
$$
\limsup_{t\to +\infty} t^{m-2}\rho\left(\frac{t}{2K}\right) \cdot \psi_\Theta (t) \le m^{m} (K+1);
$$

\end{theorem}

\begin{theorem}\label{www2}
Пусть  числа $1,\theta^1,...,\theta^m$
 линейно независимы над  $\mathbb{Z}$.
Пусть $K$  есть некоторая положительная постоянная.
Пусть функция  $\varphi (t)$  возрастает к $+\infty$  при  $t \to+\infty$.
Пусть  $\rho (t)$  есть функция обратная к функции
$ t\mapsto t\varphi ^{\frac{m-1}{m}}(t)$.
Предположим, что
$$
\limsup_{t\to+\infty} \varphi (t)\psi_\Theta (t)<K.
$$
Тогда 

{\rm (i)}  выполняется неравенство
$$
\limsup_{t\to+\infty}
\left(
\frac{t}{\rho (tK^{\frac{m-1}{m}})}\right)^{\frac{1}{m-1}}\cdot \psi_{^t\!\Theta}(t) \le 3m^2(K+1);
$$

{\rm (ii)} 
кроме того, 
если выполнено
$\varphi (t) \ge t^{m(2m-3)}$ и
если при достаточно больших значениях $t$
 функция 
$t\mapsto \varphi (t)\cdot t^{-2m+3} $ возрастает, а  функция $\rho_1(t)$ есть обратная к ней, 
  то
$$
\limsup_{t\to+\infty}
\left(
\frac{t}{\rho_1 (tK)}\right)^{\frac{1}{m-1}} \cdot \psi_{^t\!\Theta}(t) \le 3m^2(K+1).
$$

\end{theorem}

Отдельно обсудим случай  $m=2$. В этом случае согласно (\ref{minbody})  имеем 
$$
\psi_{^t\!\Theta}(t)  t^{\frac{1}{2}}\le 1,\,\,\,
\psi_{\Theta}(t) t^{2}\le 1.
$$
 Таким образом, следует использовать утверждение  (ii)  теоремы \ref{www1}
и утверждение (ii) теоремы \ref{www2}.
Получается

{\bf Следствие 1.} (Satz 2 из \cite{JTBIL}) {\it

{\rm (i)}
Пусть  функция $\varphi (t)$ такова, что  функция  $t\mapsto t^{-1}\varphi (t)$  возрастает, и при достаточно больших  $t$ 
 выполняется  $\varphi (t) \ge t^2$.   
Пусть
 $$
\limsup_{t\to+\infty} \varphi (t)\psi_\Theta (t)<K.
$$
Тогда  для функции  $\rho (t)$  обрантой и  $\varphi(t)$  имеем
$$
\limsup_{t\to+\infty} \frac{t}{\rho (tK)}\cdot \psi_{^t\!\Theta }(t)<12(K+1).
$$

{\rm (ii)}
Пусть  функция $\varphi (t)$ такова, что  функция  $t\mapsto t\varphi (t)^{-1}$  возрастает, и при достаточно больших  $t$ 
 выполняется  $\varphi (t) \le t^{\frac{1}{2}}$.   
Пусть
 $$
\limsup_{t\to+\infty} \varphi (t)\psi_{^t\!\Theta }(t)<K.
$$
Тогда  для функции  $\rho (t)$  обратной к функции  $t\mapsto t\varphi (t)^{-1}$  имеем
$$
\limsup_{t\to+\infty} {\rho\left(\frac{ t}{2K}\right)}\cdot \psi_{\Theta} (t)<4(K+1).
$$

}

Последнее следствие наиболее известно в качестве утверждения о диофантовых экспонентах.
 С учетом  определения вличины $\alpha(\Theta)$  следствие 1 превращается в

{\bf Следствие 2.} (Satz 1 из \cite{JTBIL}) {\it
При  $m=2$  выполнено равенство
\begin{equation}\label{ja22}
\alpha(\Theta) =\frac{1}{1-\alpha(^t\!\Theta)}.
\end{equation}
}
Отметим, что А.Я.Хинчин в работе \cite{dobpere}
 приводит простое и короткое  доказательство
равенства В.Ярника (\ref{ja22}).

Приведем следствие теорем \ref{www1},\ref{www2},
касающееся диофантовых экспонент  $\alpha(\Theta),\alpha(^t\!\Theta)$
 в случае произвольного значения размерности  $m$:

{\bf Следствие 3.} (Satz 3 из \cite{JTBIL}) {\it

{\rm (i)}
Всегда выполняются неравенства
$$
\alpha(\Theta) \ge (m-1) +m\alpha (^t\!\Theta ),
\,\,\,\,
\alpha (^t\!\Theta )\ge
\frac{1}{m-1}\left(1-\frac{m}{(m-1)\alpha(\Theta) +m}\right)
=\frac{\alpha (\Theta )}{(m-1)\alpha(\Theta) +m}
;
$$

{\rm (ii)}
если
$\alpha (^t\!\Theta )>\frac{m-1}{m}$, то
$$
\alpha (\Theta )\ge m-2+\frac{1}{1-\alpha (^t\!\Theta )};
$$

 {\rm (ii)}
если
$\alpha (\Theta )> m(2m-3)$, то
$$
\alpha (^t\!\Theta )\ge 
\frac{1}{m-1}\left(1-\frac{1}{\alpha (\Theta ) -2m+4}\right)
.$$}
Естественно, следствие 2 можно получить, положив в следствии 3 $m=2$.

А.Апфельбек \cite{apf}  обобщил теоремы  \ref{www1},\ref{www2}  на случай произвольных  $m,n$:

\begin{theorem}\label{www1apf}{\rm (А.Апфельбек \cite{apf})}
Пусть  матрица $\Theta$  невырождена.
Пусть $K$ -- положительное число и функция  $\varphi(t)$ 
возрастает к $+\infty$  при  $t \to +\infty$.
Пусть выполняется
$$
\limsup_{t\to+\infty}
\varphi (t) \psi_{\Theta}(t) <K.
$$
Тогда

{\rm (i1)}  при $m=1$ выполняется неравенство
$$
\limsup_{t\to +\infty} t^{n-1}\varphi \left(\frac{t^n}{2(n-1)}\right) \cdot \psi_{^t\!\Theta }(t) \le 2(n+1) K;
$$

{\rm (i2)}  при $m>1$ выполняется неравенство
$$
\limsup_{t\to +\infty} 
\left( \frac{t^n}{\rho\left(K^{\frac{m-1}{m+n-1}}t\right)}\right)^{\frac{1}{m-1}}\cdot
\psi_{^t\!\Theta} (t) \le (2(n+m))^{\frac{1}{m-1}},
$$
где $\rho (t)$  обозначает функцию,обратную к функции  $t\mapsto \left(t^m(\varphi(t))^{m-1}\right)^{\frac{1}{m+n-1}}$;
 
{\rm (ii)}  если же  при  $ m>1$ известно, что 
$$
\varphi (t) \ge 
2^{m+n-2}
K t^{\frac{2(m+n-2)^2+m-2}{n}}
 $$
 при всех достаточно больших значениях $t$  и
функция  $t\mapsto t^{-\frac{2m+n-1}{n}}\varphi (t)$  возрастает, то
  выполнено
$$
\limsup_{t\to +\infty}
\left(\frac{t^n}{\rho_1\left(K^{\frac{m-1}{m+n-2}}t\right)}\right)^{\frac{1}{m-1}}\cdot
  \psi_\Theta (t) \le 3(m+n),
$$
где $\rho_1 (t)$  обозначает функцию,обратную к функции  
$$
t\mapsto \left(t^{-\frac{(m-2)(2m+n-3)}{(m-1)n}}\varphi (t)\right)^{\frac{m-1}{m+n-2}}.
 $$
\end{theorem}
Ясно, что из теоремы \ref{www1apf}  получается следствие о  диофантовых экспонентах
 $\alpha(\Theta),\alpha(^t\!\Theta)$.

{\bf Следствие} (Теорема 6 из  \cite{apf}) {\it

{\rm (i)}
При произвольных  $m,n$  выполнено
$$
\alpha(^t\!\Theta)\ge
\frac{n\alpha(\Theta)+n-1}
{(m-1)\alpha(\Theta)+m};
$$

{\rm (ii)}
если  $m>1$  и известно,что  $$
\alpha(\Theta) >\frac{2(m+n-1)(m+n-3)+m}{n}
$$
то
$$
\alpha(^t\!\Theta )\ge
\frac{1}{m}\left( n+
\frac{n(n\alpha(\Theta)-m)-2n(m+n-3)}{(m-1)(n\alpha(\Theta)-m)+m-(m-2)(m+n-3)}\right)
 $$
}

А.Апфельбек в работе \cite{apf}
 доказал также, что   в случае $\alpha(\Theta) = +\infty$ неравенство  $\alpha (^t\!\Theta) \ge \frac{n}{m-1}$
(которое дает следствие  теоремы \ref{www1apf}) является точным (теорема 11 из \cite{apf},
ее доказательство основывается на построении сингулярных матриц   $\Theta$
 специального вида).

\subsection{Теоремы М.Лорана}\label{LORR}

В работе \cite{lora} М.Лоран получает следующий результат.

\begin{theorem}\label{lorra}
Для экспонент двумерных диофантовых приближений имеет место следующее.

{\rm (i)} Для произвольной вектор-строки  $\Theta =(\theta^1,\theta^2)\in \mathbb{R}^2$, такой, что 
 ${\rm dim}_\mathbb{Z}\Theta = 3$  для величин
\begin{equation}\label{expone}
w=\alpha(\Theta),\,\,\, w^*=\alpha (^t\!\Theta),\,\,\, v=\beta(\Theta),\,\,\, v^*=\beta (^t\!\Theta)
\end{equation}
(определения которых мы приводили  и  в пункте \ref{mravno2}, и в пункте \ref{JAP})
выполнены соотношения
\begin{equation}
\label{positiv}
2\le w\le +\infty,\,\,\,
w=\frac{1}{1-w^*},\,\,\,
\frac{v(w-1)}{v+w}\le v^*\le \frac{v-w+1}{w}. 
\end{equation}

{\rm (ii)} Для произвольного набора  из   четырех чисел  $(w,w^*,v,v^*)$,  удовлетворяющего соотношениям  (\ref{positiv})
 найдется вектор-строка  $\Theta =(\theta^1,\theta^2)\in \mathbb{R}^2$,  
 ${\rm dim}_\mathbb{Z}\Theta = 3$  такая, что
выполнено (\ref{expone}).
\end{theorem}

Скажем несколько слов о теореме \ref{lorra}.
 Первое соотношение из  (\ref{positiv})  есть следствие теоремы Минковского  (\ref{minbody});
второе соотношение есть результат В.Ярника - следствие  2 теорем \ref{www1},\ref{www2};
третье соотношение доказывается М.Лораном с помощью рассмотрения последовательностей наилучших приближений и с использованием соотношения В.Ярника (\ref{ja22}).
Принадлежащие В.Ярнику следствие теоремы \ref{T8}  из пункта \ref{mpavno1}
 и следствие  теоремы \ref{m=2}  из пункта \ref{mravno2} (естественно, при $m=1,n=2$ и $m=2,n=1$,
соответственно)
могут быть получены из утверждения  (i)  теоремы \ref{lorra}.
Утверждение (ii)  теоремы \ref{lorra}, в частности, показывает,  что 
следствия теорем  В.Ярника из пунктов \ref{mpavno1},\ref{mravno2} при $m=1,n=2$ и $m=2,n=1$ содержат точные (неулучшаемые)  оценки.
Более того из теоремы \ref{lorra} вытекает, что  
пары величин
$(\alpha(\Theta), \beta(\Theta))$
и $(\alpha (^t\!\Theta), \beta (^t\!\Theta))$
могут принимать любые допустимые (то есть, согласующиеся с неравенствами (\ref{positiv})),
значения.

Обобщение третьего неравенства из (\ref{positiv}) на случай  одной линейной формы от  $m$
 переменных имеется в работе М.Лорана \cite{lora1}.
Оно
выглядти следующим образом.

\begin{theorem}\label{lllo}
Для вектор-строки $\Theta =(\theta^1,...,\theta^m
), m\ge 2$  с условием  ${\rm dim}_\mathbb{Z}\theta = m+1$  выполнено
\begin{equation}\label{newhin}
\frac{(\alpha(\Theta)-1)\beta(\Theta)}{((m-2)\alpha(\Theta)+1)\beta(\Theta)+(m-1)\alpha(\Theta)}\le
\beta(^t\!\Theta)\le
\frac{(1-\alpha(^t\!\Theta))\beta(\Theta)-m+2-\alpha(^t\!\Theta)}{m-1}
\end{equation}
 
\end{theorem}
В свете равенства В.Ярника (\ref{ja22})  при  $m=2$  соотношение (\ref{newhin})  совпадают с третьим
соотношением из (\ref{positiv}).

Из (\ref{minbody}),  очевидно, имеем
\begin{equation}\label{dobo}
\alpha(\Theta)\ge m,\,\,\,\,\,\
\alpha(^t\!\Theta) \ge \frac{1}{m}.
\end{equation}
В случае, когда   в обеих неравенствах  (\ref{dobo})  имеет  место равенство,
неравенство М.Лорана (\ref{newhin})  совпадает с неравенством
А.Я.Хинчина (\ref{hincini}).  Если же  известно, что неравенства в (\ref{dobo})  строгие, то неравенство М.Лорана
(\ref{newhin})  будет сильнее, чем неравенство А.Я.Хинчина (\ref{hincini}).  Таким образом, теорема \ref{lllo}
утверждает, что теорема
\ref{perenoshi} может быть усилена для сингулярных наборов $\Theta$.
Автору неизвестно, имеется ли аналогичное усиление теоремы Ф.Дайсона  (имеется в виду теорема \ref{perenosdai}),
связанное с рассмотрением сингулярных матриц $\Theta$ .
Ясно, что такого рода результат должен иметь место.

Уточнение  приведенных результатов М.Лорана имеются у  Я.Бюжо и М.Лорана   \cite{blora}.

\section{О размерности Хаусдорфа множеств сингулярных систем}\label{HD}

О точных значениях и оценках размерности Хаусдорфа множеств сингулярных матриц известно не так много.
Ниже мы перечислим некоторые результаты. Отметим, что ряд задач, связанных с отысканием
размерностей Хаусдорфа поставлен Я.Бюжо и М.Лораном \cite{BULmo}.
Насколько известно автору, в основном, исследовался только случай  $n=1$, причем получены только оценки размерности Хаусдорфа.  Исключением является работа И.Чеунга, в которой получен следующий замечательный результат.

\begin{theorem}\label{CHU}{\rm (И.Чеунг \cite{CHUU})}
 Хаусдорфова размерность сингулярных систем
(в смысле определения Хинчина при $n=2,m=1$)
 вида  $\left(\begin{array}{c}\theta_1\cr \theta_2\end{array}\right)$  равна $4/3$.
\end{theorem}

  По причине  переноса видим, что в теореме \ref{CHU}  можно
говорить о хаусдорфовой размерности сингулярных систем  $ (\theta^1,\theta^2), n=1,m=2$.

Далее мы будем рассматривать множества $E_m(\alpha )$ сингулярных векторов
$\Theta=(\theta^1,...,\theta^m)$ (то есть случай произволного  $m\ge 2$  и $ n=1$),
определенное следующим образом:
$$
E_m(\alpha ) =
\{\Theta \in \mathbb{R}^m:\,\, \lim_{t\to +\infty} t^\alpha \psi_\Theta (t) =0\}.
$$

Сформулируем наиболее точные известные автору оценки снизу и сверху хаусдорфовой размерности.

Следующий результат принадлежит Р.К.Бейкеру \cite{bak2}; он улучшает предшествующие результаты
Р.К.Бейкера \cite{bak1}, К.Ю.Явида \cite{bssr} (он первый привел контрпример к  
предположению Р.К.Бейкера о том, что хаусдорфова размерность множества сингулярных векторов чрезвычайно мала, сформулированному в \cite{bak1}, однако затем результат К.Ю.Явида был {\it количествено}
улучшен Б.Ринном)
  и Б.Ринна \cite{ry1}.

\begin{theorem}\label{bakke}{\rm (Р.К.Бейкер \cite{bak2})}
При $m\ge 3$  и $\alpha >m$ для размерности Хаусдорфа множества
$
E_m(\alpha )$
выполняется оценка снизу
$$
{\rm HD}\,
E_m(\alpha )
\ge n-2+n/\alpha.
$$
\end{theorem}

Оценку сверху получил Б.Ринн \cite{ry2}, улучшив результат Р.К.Бейкера из \cite{bak1}:

\begin{theorem}\label{rynne}{\rm (Б.Ринн \cite{ry2})}
При $m\ge 3$  и $\alpha >m$ для размерности Хаусдорфа множества
$
E_m(\alpha )$
выполняется оценка  сверху
$$
{\rm HD}\,
E_m(\alpha )\le
  n-2+ (2n+2)/(\alpha+1).
$$
\end{theorem}

В случае $m=2$  обе наилучшие   оценки (верхняя имеется в \cite{bak2},  нижняя
следует из результата работы \cite{bak1}) некоторое время принадлежали Р.К.Бейкеру:

\begin{theorem}\label{bakke2}{\rm (Р.К.Бейкер \cite{bak1},\cite{bak2})}
Для размерности Хаусдорфа множества
$
E_2(\alpha )$
выполняются оценнки
$$
\frac{2}{\alpha} \le
{\rm HD}\,
E_2(\alpha )
\le \frac{6}{\alpha + 1}.
$$
\end{theorem}

Нижняя оценка  теоремы \ref{bakke2}, по-видимому, является наилучшей известной до сих пор.

Применяя равенство В.Ярника (\ref{ja22}), неравенство (\ref{betatet})  из следствия теоремы \ref{T8} (тоже принадлежащее В.Ярнику) и частный случай ($n=2,m=1$) утверждения, имеющегося в работе М.М.Додсона \cite{dods} (теорема \ref{dodd} ниже),
Я.Бюжо и М.Лоран в работе \cite{BULmo}
 выводят неравенство
$$
{\rm HD}\,
E_2(\alpha )
\le \frac{3\alpha}{\alpha^2-\alpha + 1},
$$
более сильное, чем верхняя оценка из теоремы  \ref{bakke2}.
Точнее, в \cite{BULmo} доказано, что
$$
{\rm HD}\,\{
\Theta \in \mathbb{R}^m,\,\alpha (\Theta ) \ge \alpha \} \le\frac{3\alpha}{\alpha^2-\alpha + 1}.
$$
Приведем общую формулировку использованного результата из  \cite{dods}.
 Напомним, что  показатель $\beta (\Theta )$  обозначает
супремумы тех  $\gamma$,
 для которых     выполнено 
 $$
\liminf_{t\to +\infty} t^\gamma
\psi_\Theta (t) <+\infty
.
$$

\begin{theorem}\label{dodd}
 При каждом  $\tau > m/n$  множество
$$
\hbox{\got W}_{m,n} (\tau ) =\{ \Theta:\,\, \beta (\Theta ) \ge \tau\}
$$
 имеет  размерность Хаусдорфа
$$
{\rm HD}\, \hbox{\got W}_{m,n} (\tau )=
(m-1)n+\frac{m+n}{\tau +1}.
$$
\end{theorem}

Приведем еще два простых следствия из сформулированных выше результатов.
Из теорем \ref{bakke},\ref{rynne}  для множества
$$
E_m(\infty) = \bigcap_{\alpha >m} E_m(\alpha)
$$
получается следующее 

{\bf Следсвтие 1.} \,\,{\it
Для размерности Хаусдорфа множества $
E_m(\infty)$
 выполнено
$$
{\rm HD}\, 
E_m(\infty) = n-2.
$$}

Отметим также, что из теоремы \ref{bakke}  получается

{\bf Следствие 2.} \,\,{\it
Размерность Хаусдорфа множества 
всех сингулярных (в смысле определения А.Я.Хинчина) при  $n=1$  и произвольном  $m \ge 2$
 оценивается снизу величиной  $m-1$.
}

В заключение
настоящего пункта хочется сослаться на книги В.И.Берника и Ю.В.Мельничука \cite{ber}
 и В.И.Берника и М.М.Додсона \cite{berber},
посвященные диофантовым приближениям и размерности Хаусдорфа.

\section{Приближения с целыми неотрицательными числами}\label{neotri}

\subsection{Двумерные приближения}

Положим $\tau = \frac{1+\sqrt{5}}{2}$.
 В работе 
\cite{SCHmon}
В.М.Шмидт получил следующий результат.

\begin{theorem}\label{TT}
Пусть вещественные числа $\theta^1,\theta^2$
линейно независимы вместе с единицей над $\mathbb{Z}$.
Тогда найдется последовательность целых двумерных векторов
 $(x_1(i), x_2(i))$
 таких, что

1.\,\, $x_1(i), x_2(i) > 0$;

2.\,\, $||\theta^1x_1(i)+\theta^2x_2(i) ||\cdot (\max \{x_1(i),x_2(i)\})^\tau \to 0$ при $ i\to +\infty$.
\end{theorem}

Известная гипотеза (см. \cite{SCHmon,SCHL})
о  том, что показатель  $\tau$ в теореме \ref{TT}  может быть заменен на  $2-\varepsilon$ 
с произвольным сколь угодно малым $\varepsilon$, до сих пор не доказана.
Тем не менее, ряд математиков получили обобщения и уточнения теоремы В.М.Шмидта (см., например, \cite{T1}, \cite{T2}, \cite{BK}).
 Доказательство В.М.Шмидта использовало утверждение о независимости наилучших приближений линейной формы
(теорема \ref{NT1}, $m=2,n=1$) и фактически было связано с  разбором  двух случаев, различающихся в зависимости от того, насколько сингулярным является набор  $\theta^1,\theta^2  \, (n=1,m=2)$.
(Фактически, надо разбирать два случая. Первый случай, когда для всех  $\nu$  таких, что наилучшие приближения  
${\bf z}_{\nu-1}.{\bf z}_{\nu},{\bf z}_{\nu+1}$
независимы, выполнено
$$
\zeta_\nu \le M_\nu^{\frac{1}{\tau}}\cdot M_{\nu+1}^{-\frac{\tau}{\tau-1}}.
$$
Второй случай, когда для бесконечно многих   $\nu$ выполняется противоположное неравенство.)

В \cite{MoX}  был разобран случай плохо приближаемых чисел $\theta^1,\theta^2$.
 Сформулируем результат из  \cite{MoX}.
 Для вещественного
 $\gamma \ge 2$  рассмотрим функцию
$$
g(\gamma ) = \tau +\frac{2\tau - 2}{\tau^2\gamma -2}.
$$
Сразу видно, что 
$g(\gamma )$  строго убывает и что 
$$
g(2) =2,\,\,\,\,\,\, \lim_{\gamma\to+\infty} g(\gamma ) = \tau.
$$
Также для положительного
 $\Gamma$ определим величину 
$$
C(\Gamma ) = 2^{18} \Gamma^{\frac{\tau -\tau^2}{\tau^2\gamma-2}}.
$$
 
\begin{theorem}\label{TTT}
 Пусть набор вещественнных чисел $\theta^1,\theta^2$
 является плохо приближаемым в следующем смысле.
Для некоторых положительных постоянных
$\Gamma\in (0,1)$ и $\gamma \ge 2$ неравенство
\begin{equation}
||\theta^1m_1+\theta^2m_2||\ge \frac{\Gamma}{(\max \{|m_1|,|m_2|\})^\gamma}
\label{badq}
\end{equation}
выполнено для всех целых векторов $(m_1,m_2)\in \mathbb{Z}^2\setminus \{(0,0)\}$.
Тогда найдется бесконечная последовательность  целых двумерных векторов
 $(x_1(i), x_2(i))$
такая, что

1.\,\, $x_1(i), x_2(i) > 0$;

2.\,\, $||\theta^1x_1(i)+\theta^2x_2(i) ||\cdot (\max \{x_1(i),x_2(i)\})^{g(\gamma )} \le C(\Gamma )$ для всех  $ i$.

\end{theorem}

 Мы не будем останавливаться на доказательстве этой теоремы. Оно развивает оригинальны идеи В.М.Шмидта и изложено в \cite{MoX}  достаточно подробно.

\subsection{О линейных формах с числом переменных  $k>2$}

В случае  $k \ge 3$    линейная форма от $k$  переменных,
вообще говоря, может не принимать малые значения при
 положительных $x_j$.  Приведем результат В.М.Шмидта из  \cite{SCHmon}:

\begin{theorem}
\label{manno}
Найдется набор чисел  $\Theta = (\theta_1,...,\theta_k), k\ge 3$ такой, что

$\bullet$  ${\rm dim}_\mathbb{Z}\Theta = k+1$;

$\bullet$ для любого положительного $\varepsilon$  найдется положительное  $c(\varepsilon)$
такое, что
$$
|| \theta_1x_1+...+\theta_k
x_k||>c(\varepsilon) \left( \max_{1\le i\le k}
|x_i|\right)^{-2-\varepsilon}
$$
для всех  целых  $x_1,...,x_k$  с условием  $x_i >0$.
\end{theorem}
 Для доказательства теоремы \ref{manno}  используется результат Г.Давенпорта и В.М.Шмидта \cite{DS1} о существовании чисел
с аномальными совместными приближениями:

\begin{theorem}
\label{odnochislo}
Пусть  $n \ge 2$ и $\psi (t)$   есть положительнозначная функция вещественного аргумента  $t$.
Тогда найдется набор чисел  $\Theta = (\theta_1,...,\theta_n)$ такой, что

$\bullet$  ${\rm dim}_\mathbb{Z}\Theta = n+1$;

$\bullet$  для любого достаточно большого $t$  найдется натуральное $q\le t$ такое, что
неравенства
\begin{equation}\label{qoqu}
||q\theta_j||\le \psi (t)
\end{equation}
выполнены для всех значений индекса  $j, 1\le j \le n$,
за возможным исключением одного значения $j_0 = j_0(t)$.
\end{theorem}

Отметим, что при   $n\ge 2$ и $\psi (t) = O(t^{-1}), t\to +\infty$
обеспечить выполнения неравенства (\ref{qoqu})  для всех  $j$
  из промежутка  $1\le j \le n$ 
нельзя, - это противоречило бы теореме \ref{T7}.

Следуя В.М.Шмидту \cite{SCHmon}, покажем, что из теоремы \ref{odnochislo}  
с помощью простейших сображений переноса
 вытекает теорема
\ref{manno}. Для этого надо взять $n = k-1 $ и   $\psi (t) = e^{-t}$. 
Для чисел  $\theta_1,...,\theta_{k-1}$,  существование которых обеспечивается теоремой 
\ref{odnochislo},
найдем такое число 
$\theta_k$,
для которого с некоторой положительной постоянной  $c_1 (\varepsilon)$ выполнено
$$
||\theta_i u + \theta_k v ||>c_1(\varepsilon ) (|u|+|v|)^{-2-\varepsilon} ,\,\,\, 1\le i \le k-1
$$
для всех целых  $u,v$  с условием  $v>0$.
 (Стандарные аргуметны, связанные с применением леммы Бореля-Кантелли, показывают, что это неравенство выполнено для {\it почти всех }  чисел $\theta_k$.  Мы проводили подобные рассуждения в пунктах  \ref{vr2},\ref{vr3} настоящей статьи.) 

Теперь для целой точки  $(x_1,...,x_k)$  с достаточно большим значением
 $M=\max_{1\le i \le k}|x_i|$  полагаем  $t = (\log M)^2$  и выбираем  $q$
 из интервала   $1\le q\le t$,  удовлетворяющим  заключению теоремы \ref{odnochislo}.
 Не ограничая общности, считаем, что  $j_0(t) = 1$.  Значит, если $ x_k >0$,  то
$$
||\theta_1x_1+\cdots+\theta_k x_k ||\ge
q^{-1}||\theta_1qx_1+\cdots +\theta_kqx_k||\ge
$$
$$\ge 
q^{-1} (||\theta_1 qx_1 +\theta_kqx_k||-
M (||\theta_2q||+\cdots + ||\theta_{k-1}q||)) >
$$
$$>
q^{-1} (c_1 (\varepsilon/3) (qM)^{-2-(\varepsilon/3)}
- kM e^{-t})>M^{-2-\varepsilon},
$$
и теорема \ref{manno} доказана.

\section{Задача  В.В. Козлова}\label{Zadako}

Рассмотрим набор вещественных чисел
$\Theta = (\theta^1,...,\theta^m)$,  линейно незвисимых  над  $\mathbb{Z}$
и вещественно-значную функцию  $f(x_1,...,x_m)$, периодичную по каждому из аргументов  $x_i$  с периодом 1,
о которой мы будем предполагать, что она достаточно гладкая (например, непрерывная).
Всюду ниже в этом пункте будем предполагать, что выполняется условие нулевого среднего
\begin{equation}
\label{zeromen}
\int_0^1\cdots \int_0^1 f(x_1,...,x_m) dx_1\dots dx_m = 0
.
\end{equation}
 Согласно знаменитой теореме Г.Вейля \cite{GW} 
для непрерывной функции  $f$
для любой начальной фазы
${\bf y} =(y_1,...,y_m)\in \mathbb{R}^m$
выполнено
$$
I (T,{\bf y})=
I^{[f,\Theta]} (T,{\bf y}) =
\int_0^T f(\theta^1t+y_1,...,\theta^mt+y_m ) dt = o(T),\,\,\, T\to +\infty.
$$
 В.В.Козлов поставил задачу об осцилляции и возвращаемости интеграла  $
I (T,{\bf y})$.

Мы будем говорить, что интеграл $
I (T,{\bf y})$
{\it осциллирует}
(при заданной начальной фазе  ${\bf y}$), если
 множества моментов времени
$$
\hbox{\got T}_+ =\{
T\in \mathbb{R}_+|\,\,
I (T,{\bf y}) >0\},\,\,\,
\hbox{\got T}_- =\{
T\in \mathbb{R}_+|\,\,
I (T,{\bf y}) <0\}  
$$
 оба являются неограниченными сверху.

Мы будем говорить, что интеграл $
I (T,{\bf y})$
{\it возвращается}
(при заданной начальной фазе  ${\bf y}$), если
$$
\liminf_{T\to +\infty} |
I (T,{\bf y})
| = 0.
$$

У описанной выше непрерывной постановки задачи имеется дискретный аналог.

Рассмотрим набор вещественных чисел
$\Theta = (\theta^1,...,\theta^m)$,  линейно незвисимых  вместе с единицей над  $\mathbb{Z}$
и (достаточно гладкую) вещественнозначную функцию  $F(x_1,...,x_m)$, периодичную по каждому из аргументов  $x_i$  с периодом 1, и имеющую нулевое среднее:
\begin{equation}
\label{zeromen1}
\int_0^1\cdots \int_0^1 F(x_1,...,x_m) dx_1\dots dx_m = 0
.
\end{equation}
 Упомянутая выше теорема Г.Вейля о равенстве пространственного и временного средних в этом случае утверждает, что
 для любой начальной фазы
${\bf y} =(y_1,...,y_m)\in \mathbb{R}^m$
выполнено
$$
S (Q,{\bf y})=
S^{[F,\Theta]} (Q,{\bf y}) =
\sum_{s=1}^QF(\theta^1s+y_1,...,\theta^ms+y_m ) dt = o(Q),\,\,\, Q\to +\infty.
$$
Осцилляция и возвращаемость суммы $S (Q,{\bf y})$ при  $Q\to +\infty$
 определяется, естественно, аналогично осцилляции и возвращаемости интеграла
 $I(T,{\bf y})$.

Если для функции $m$ переменных $f(x_1,...,x_m)$  рассмотреть функцию  $F_{z}(y_1,...,y_{m-1})$  от
$m-1$ переменных $y_1,...,y_{m-1}$ (переменная  $z$  теперь рассматривается как параметр), определяему соотношением
$$
F_{z}(y_1,...,y_{m-1})=
\int_0^{\theta_m^{-1}}
 f(\theta^1t+y_1,...,\theta^{m-1}k+y_{m-1},\theta^mt+z ) dt 
$$
 и новый набор частот
 $\Theta^* =(\theta^1/\theta^m,...,\theta^{m-1}/\theta^m)$
(состоящий из чисел, линейно независимых {\it вместе с единицей}
над  $\mathbb{Z}$),
то из возвращаемости 
(или осцилляции) 
суммы
$S^{[F_z,\theta]} (Q,{\bf y})$  при
при некотором значении параметра  $z$  и при некоторой  начальной фазе  ${\bf y } = (y_1,...,y_{m_1})$
  будет следовать возвращаемости (или осцилляция) интеграла
$I(T,{\bf x})$ при соответствующей начальной фазе  ${\bf x} $.
Заметим здесь, что обратное утверждение, вообще говоря, может оказаться неверным.

\subsection{Лемма Переса и теорема Халаса}
В этом пункте мы приведем два общих утверждения из эргодической теории и обсудим, какие следствия они дают в рассматриваемой задаче.

Напомним, что под {\it  динамической системой} обычно понимается  вероятностное просттранства  $(\Omega,{\cal A},\mu)$
с заданным эргодическим преобразованием  $T$  множества  $\Omega$  в себя
(об основных понятиях эргодической теории см. \cite{korsifo}).
Эргодическая теорма Биркгофа утверждает, что  для любой интегрируемой функции $g$ 
из  $\Omega$  в
$\mathbb{R}$  для почти каждого  $x \in \Omega$  выполняется
$$
\frac{1}{Q}\sum_{s=1}^Q g(T^sx) \to \int_\Omega g d\mu,\,\,\,\,
Q\to +\infty.
$$

$\bullet $\,\,
Следующее утверждение известно как лемма Ю.Переса.

\begin{lem}\label{perre} {\rm (Ю.Перес \cite{perre})}
Пусть  $\Omega$  компактно, а преобразование  $T$  непрерывно. Тогда для любой непрерывной функции  $g$  из  $\Omega$  в
$\mathbb{R}$  найдется точка  $x\in\Omega$ такая, что
$$
\frac{1}{Q}\sum_{s=1}^Q g(T^sx) \ge \int_\Omega g d\mu.
$$
 
\end{lem}
В рассматриваемой нами задаче лемма \ref{perre}  влечет

{\bf  Следствие.} \,{\it
Для  произвольного набора $\Theta$, состоящего из линейно независимых вместе с единицей над  $\mathbb{Z}$
частот и
для любой непрерывной функции $F(x_1,...,x_m)$  периодичной по каждому аргументу с периодом 1  и имеющей нулевое среднее 
найдется точка  ${\bf y}\in [0,1)^m$  такая, что
$$
S^{[F,\Theta]} (Q,{\bf y}) \ge 0 \,\,\,\,\forall Q\in \mathbb{N}.
$$
}

Аналогичное утверждение для интеграла $I^{[f,\Theta]} (T,{\bf y})$
  имеется у П.Боля \cite{boll}.
В несколько более сильном виде оно имеется у В.В.Козлова \cite{KOZ1},\cite{KOZ2}:

\begin{theorem}\label{koooo}
{\rm (В.В.Козлов \cite{KOZ1},\cite{KOZ2})}

{\rm (i)} Пусть частоты  $\theta^1,...,\theta^m$  линейно независимы над  $\mathbb{Z}$ и периодическая функция $f$
непрерывна и имеет нулевое среднее.
Тогда найдется точка 
 ${\bf y}\in [0,1)^m$
такая, что
$$
I^{[f,\Theta]} (T,{\bf y})\ge 0,\,\,\, \forall T \in  \mathbb{R},
$$
и, кроме того,
$f ({\bf y}) =0$.

{\rm (ii)} Если  частоты  $\theta^1,...,\theta^m$  линейно зависимы над  $\mathbb{Z}$, то существует, по крайней мере,
две различные точки  ${\bf y}_1,{\bf y}_2\in [0,1)^m$,
для которых выполнено заключение пункта  {\rm (i)}.
\end{theorem}

$\bullet $ \,\, Сформулируем одну из теорем из работы Г.Халаса \cite{halas}, касающейся поведения сумм Биркгофа для интергируемой на $\Omega$  функции $g$.

\begin{theorem}\label{Haaaa}
{\rm (Г.Халас \cite{halas})}
Для любой интегрируемой 
$g$ разность
\begin{equation}\label{veli}
\sum_{s=1}^Q g (T^s{ x})
-Q\int_\Omega gd\mu
\end{equation}
 бесконечно много раз меняет знак (в слабом смысле)  
для почти всех  $x\in\Omega$

\end{theorem}

То, что величина меняет знак в  {\it слабом смысле},
в данном случае означает что, когда величина обращается в ноль, мы можем также считать ее поменявшей знак, то есть 
теорема \ref{Haaaa} утверждает, что  величина (\ref{veli})
не может, начиная с какого-то момента, быть всегда  строго положительной или строго отрицательной.

{\bf Следствие.}\,{\it
Для  произвольного набора $\Theta$, состоящего из линейно независимых вместе с единицей над  $\mathbb{Z}$
частот и для произвольной интегрируемой
функции $F(x_1,...,x_m)$  периодичной по каждому аргументу с периодом 1  и имеющей нулевое среднее 
величина $S^{[F,\Theta]} (Q,{\bf y})$  бесконечно много раз меняет знак (в слабом смысле)  для почти всез значений начальных фаз  ${\bf y}$.

}

\subsection{Индивидуальная возвращаемость}

В случае  $m=2$ возвращаемость интеграла
$I^{[f,\Theta]} (T,{\bf y})$ для  периодической функции $f(x_1,x_2) \in C^2([0,1]^2)$  с нулевым средним, произвольного набора частот
$\Theta =(\theta^1,\theta^2).\,\,\, \theta^1/\theta^2 \not\in \mathbb{Q}$
 и для произвольных начальных фаз  $(y_1,y_2)$
 была доказана самим В.В. Козловым в  \cite{KOZ1},\cite{KOZ2}. В этих же работах В.В.Козлов заметил, что
при  $m=2$
имеет место более точное утверждение о возвращаемости
 интеграла $I(T,{\bf y})$, а именно
найдется последовательность 
 положительных {\it  целых}
 чисел  $T_\nu $  такая, что одновременно
\begin{equation}
\label{nadstro}
\max_{1\le i\le 2}||\theta^iT_\nu ||\to 0,\,\,\,
I(T_\nu,{\bf y}) \to 0,\,\,\,
\nu \to +\infty.
\end{equation}
Далее В.В.Козлов
заметил, что из условия $f(y_1,y_2)\neq 0$
 с помощью уже доказанной  "сильной возвращаемости"  (\ref{nadstro})
 следует осцилляция интеграла
 $I(T,{\bf y})$  при рассмотренном
значении начальной фазы  ${\bf y} = (y_1,y_2)$.

Эти результаты В.В.Козлова были усилены Е.А.Сидоровым 
\cite{SYD},
о чем мы подробно поговорим в следующем пункте.
 В случае произвольного значения размерности
$m$
результаты об осцилляции и возвращаемость были получены
С.В. Конягиным (в случае {\it  нечетной}  функции  $f$ (см. \cite{KONY}))
 и Н.Г. Мощевитиным в общем случае (см. \cite{momatza},\cite{mograz}).
Ниже мы приведем формулировки этих результатов и дадим некоторые  комментарии.

Сначала мы сформулируем две теоремы, доказанные Н.Г.Мощевитиным  \cite{momatza},\cite{mograz}.

\begin{theorem} \label{koz1}
  Предположим что  $m \ge 2$, что функция  $f$  принадлежит классу
 $C^w ([0,1]^m)$, где
\begin{equation}
\label{wewe}
w= w(m)= 
[\exp (20m\log m)] 
\end{equation}
и что выполнено условие нулевого среднего  (\ref{zeromen}).
 Тогда 
для любых частот
$\Theta = (\theta^1,...,\theta^m)$,  линейно незвисимых  над  $\mathbb{Z}$
при любой начальной фазе  ${\bf y}$
 интеграл
 $I(T,{\bf y})$  возвращается.
\end{theorem}

\begin{theorem} \label{koz2}
 Если в условия теоремы \ref{koz1}  дополнительно потребовать, что для начальной фазы  ${\bf y} = (y_1,...,y_m)$
выполняется
$$
f(y_1,...,y_m) \neq 0,
$$
 то интеграл  $I(T,{\bf y})$  осциллирует. 
\end{theorem}
Доказательства теоремы  \ref{koz1} можно найти в работах автора 
\cite{momatza},\cite{mograz}.  Однако там 
 она сформулирована с худшим значением   показателя гладкости  $w(m)$.
  Тем не менее, доказательство из \cite{momatza},\cite{mograz} дословно проходит для показателя (\ref{wewe}).
К сожалению, соответствующая выкладка
проведена только в диссертации автора \cite{diss}.  Тем не мене, она легко восстанавливается.
Теорема \ref{koz2}, вообще говоря, не является непосредственным следствием теоремы  \ref{koz1},
но ее доказательство совершенно аналогично доказательству теоремы \ref{koz1}.
Точное рассуждение тоже приведено в диссертации автора  \cite{diss}.

 Из теоремы \ref{koz1}
очевидным образом получаем 

{\bf Следствие.}\,\,{\it
Пусть функция  $F(x_1,...,x_m)$
принадлежит классу $ C^{w(m)}([0,1]^m)$
и имеет нулевое среднее
(\ref{zeromen1}).
Тогда
для любых частот
$\Theta = (\theta^1,...,\theta^m)$,  линейно незвисимых  вместе с единицей над  $\mathbb{Z}$
и для всякой начальной фазы  ${\bf x}$  выполнено
$$
\liminf_{Q\to +\infty} |
S(Q,{\bf x})
| <+\infty.
$$
}
Видим, что следствие не обеспечивает возвращаемости суммы  $S(Q,{\bf y})$.
По-видимому, вопрос о возвращаемостии суммы $S(Q,{\bf y})$
остается открытым.

Отметим также, что доказательства теорем \ref{koz1},\ref{koz2}  связаны и исследованием распределений наилучших диофантовых приближений линейной формы (случай  $n=1$) и с тонким анализом осцилляции гармоник в разложении
функции $f$ 
в кратный ряд Фурье.

Сформулируем теперь результат С.В.Конягина
из \cite{KONY}.

\begin{theorem} \label{koz3}
  Предположим что  $m \ge 3$, что функция  $f$  принадлежит классу
 $C^{w_1} ([0,1]^m)$, где
\begin{equation}
\label{wewe1}
w_1= w_1(m)= 3m\cdot 2^{m-1}
\end{equation}
и что выполнено условие нулевого среднего  (\ref{zeromen}).
 Дополнительно предположим, что
 для всех  ${\bf x}$  выполнено
\begin{equation}\label{nech}
f(-x_1,...,-x_m) = - f (x_1,...,x_m).
\end{equation}
 Тогда  
 интеграл
 $I(T,{\bf 0})$  возвращается.
\end{theorem}

 Отметим, что
 $
w_1(m)< w(m)
$,
 и, следовательно,
в теореме С.В.Конягина о нечетной функции условие на гладкость  $f$ 
 является более слабым, чем в теореме
\ref{koz1}.

Заметим, что теорема \ref{koz3}  может быть обобщена следующим образом:
если рассмотреть не одну функцию $f$, удовлетворяющую условию теоремы \ref{koz3}
 а {\it  конечный набор}  функций  $f_1,...,f_r$, каждая из
которых удовлетворяет условию теоремы \ref{koz3}, то
будет выполнено
$$
\liminf_{T\to +\infty} 
\sum_{1\le j \le r}|
I ^{[f_j, \Theta]}(T,{\bf 0})
| = 0.
$$
 Доказательство из 
\cite{KONY}
обобщается непосредственно.
 Сделаем несколько замечаний о том, какой гладкости функции  $f$  недостаточно для возвращаемости.
 Н.Г.Мощевитин (см. \cite{motre}), основываясь на идеях Д.В.Трещева, обобщил пример А.Пуанкаре
 \cite{PUA},\cite{PUA1} и доказал следующее утверждение.

\begin{theorem} \label{koz4}

Пусть  набор частот
$ \Theta= (\theta^1, ... , \theta^m)$ 
является плохо приближаемым, то есть 
с некоторой положительной
$ \Gamma = \Gamma(\Theta) $
  для всех ненулевых целых векторов 
 $ {\bf k} = (k_1,...,k_m) \in \mathbb{ Z}^m$
выполняется неравенство
$$
|
k_1\theta^1 + \dots + k_m\theta^m
|
>
\Gamma\cdot \left(\max_{1\le i\le m}|k_i|\right)^{1-m}.
$$ 
 
 Тогда найдется 1-периодическая по каждой переменной функция
 $f^{[\Theta ]}(x_1,...,x_m)$  с нулевым средним класса
$ C^{m-2} ([0,1]^m)\setminus C^{m-1}([0,1]^m)
$
 такая, что
$$
\lim_{T\to +\infty }
I^{[f^{[\Theta ]},\Theta]} (T,{\bf 0}) 
=+\infty.
$$

\end{theorem}
Отметим, что для плохо приближаемого набора условие линейной независимости над $\mathbb{Z}$
 выполнено автоматически.

С другой сторны,  построив некоторый сингулярный набор частот  $\Theta$,
С.В.Конягин в \cite{KONY}  получил следующий результат.

\begin{theorem} \label{koz5}
   При  $m \ge 4$ положим
$$
w_3= w_3(m) =
\left[
\frac{2^{m-1}(m-2)^{m-2}}{(m-1)^{m-1}}\right]-1
.
$$
Тогда найдется функция  периодическая $f\in C^{w_3}([0,1]^m)$
и набор частот
$ \Theta= (\theta^1, ... , \theta^m)$,
состоящий из линейно независимых над $\mathbb{Z}$ чисел, 
такие, что выполнено (\ref{nech}) и что
$$
\lim_{T\to +\infty }
I^{[f^{[\Theta ]},\Theta]} (T,{\bf 0}) 
=+\infty.
$$
\end{theorem}

Заметим, что при  $m\ge 9$  теорема \ref{koz5}  дает пример более гладкой функции  $f$ 
(чем теорема  \ref{koz4}),  для
которой нет  возвращаемости интеграла $I^{[f^{[\Theta ]},\Theta]} (T,{\bf 0}) $.

Отрицательный результат об отсутствии одновременной возвращаемости в общем случае имеется у Н.Г. Мощевитина
в \cite{momatza}. Он тоже связан с использованием сингулярных наборов  $\Theta$.
 Приведем его формулировку.

\begin{theorem} \label{koz6}
  
Пусть $ \Phi  (t)$ -
сколь угодно быстро убывающая к нулю функция.
Тогда найдутся  
две вещественнозначные функции
$$
f_j (x_1,x_2,x_3) 
=
\sum_{ (k_1,k_2,k_3) \in \mathbb{Z}^3\setminus \{{\bf 0} \}}\textit{}
f_{j;k_1,k_2,k_3} \exp( 2 \pi  i ( k_1x_1 +k_2x_2 + k_3x_3 )),
\,\,\, j = 1,2,$$
такие, что их коэффициенты Фурье допускают оценку
$$
|f_{j;k_1,k_2,k_3}| \le \Phi( \max_{1 \le j \le 3 } |k_j|) \hspace{2mm} \forall
k \in \mathbb{Z}^3,\hspace{2mm}j = 1,2
,$$
но для некоторого набора  $\Theta= (\theta^1,\theta^2,\theta^3)$, состоящего из трех чисел,
линейно независимых над  $\mathbb{Z}$  выполняется
$$
\liminf_{T\to +\infty} 
\sum_{1\le j \le 2}|
I ^{[f_j, \Theta]}(T,{\bf 0})
| = +\infty.
$$

\end{theorem}

Отметим также  работу автора  \cite{mmmz3},
в которой имелось частичное решение проблемы при $m=3$.

\subsection{Равномерная возвращаемость}

Когда мы говорим о {\it  равномерной возвращаемости}, мы имеем в виду что
 возвращаются величины
$$
J(T) = J^{[f, \Theta]}(T) =\sup_{{\bf y}\in \mathbb{R}^m} | I ^{[f, \Theta]}(T,{\bf y})|
,\,\,\,\,
R(Q) = R^{[f, \Theta]}(Q)= \sup_{{\bf y}\in \mathbb{R}^m}|S ^{[F, \Theta]}(Q,{\bf y})|
.
$$

Результат Е.А.Сидорова, упоминавшийся выше, относится именно к равномерной возвращаемости. Сейчас мы его сформулируем.

\begin{theorem} \label{koz6}
 {\rm  (Е.А.Сидоров \cite{SYD})}
Пусть  $m=1$ и функция  $F(x)$  является 1-периодической абсолютно непрерывной функцией вещественной переменной  $x$. Пусть  $\theta \not \in \mathbb{Q}$.
Тогда имеет место равномерная возвращаемость суммы
$S ^{[F, \theta]}(Q,{y}),$
то есть
$$
\liminf_{Q\to +\infty } R(Q) = 0
.$$

\end{theorem}
 Доказательство Е.А.Сидорова использует аппарат  цепных дробей.

 Как обнаружено в \cite{monon}
случай  $m > 1$  существенно отличается от случая  $m=1$.
(В работе  \cite{monon}  автор доказал отсутствие равномерной возвращаемости  уже при  $m=2$ и
даже для гладкой  $F(x_1,x_2)$. Для этого он использовал сингулярный
 набор  $(\theta_1,
\theta_2)$,  передоказав 
  теорему А.Я. Хинчина 1926 года (теорему \ref{T1}),
о которой он в то время не знал. 
  В настоящем обзоре мы приводим  формулировку более общего результата Е.В.Коломейкиной и Н.Г.Мощевитина \cite{kolo}  и основного вспомогательного утверждения (о существовании сингулярных наборов  $\Theta $
специального вида при  $n=1$).

  Рассмотрим разложение в ряд Фурье
$$
F(x_1,...,x_m) =
\sum_{(k_1,...,k_m)
\in \mathbb{Z}^m\setminus\{{\bf 0}\}}
F_{k_1,...,k_m} \exp (2\pi i (k_1x_1+...+k_mx_m)).
$$

 \textit{Спектром функции} $F$ будем называть
множество  $$
{\rm spec} \,F=\{(k_1,...,k_m)\in \mathbb{Z}^m:\,\,\, F_{k_1,...,k_s}\neq 0\}.$$

Следующий результат показывает, что равномерной возвращаемости может не быть даже
в случае очень большой гладкости функции  $F$.

\begin{theorem}
\label{koz100}
 Пусть периодическая  функция $F$ с
периодом 1 принадлежит классу $C^1([0,1]^m)$ и
 выполнено  (\ref{zeromen1}).

Рассмотрим следующие два условия:

{\rm (A)} \,\,  для любой положительнозначной функции $\lambda(t)=o(t),$
$t\rightarrow +\infty,$ найдется такой набор чисел
$\Theta =(\theta^1,...,\theta^m)$ линейно независимых вместе с единицей над
$\mathbb{Z}$ такой, что
$$R^{[F,\Theta]} (Q) 
> \lambda(Q)$$
при всех достаточно больших $Q$;

{\rm (B)}\,\, найдется положительное $R$ и ненулевая целая точка
$(p_1,...,p_m)\in \mathbb{Z}^m$ такие, что ${\rm spec}\, F\subset {\cal B}(R)\bigcup {\cal L}(p),$ где ${\cal B}(R)$ обозначает шар 
в  $\mathbb{R}^m$
с
центром в нуле и радиусом $R$, а ${\cal L}(p)$ обозначает
прямую в  $\mathbb{R}^m$, проходящую через начало координат и точку $p$.

Тогда условие {\rm (A)} эквивалентно отрицанию условия \rm{(B)}.
\end{theorem}

Для доказательства теоремы  \ref{koz100}
используетяся следующая 

\begin{lem}\label{lkolo}
В предположении отрицания условия {\rm (B)} теоремы \ref{koz100}
  для произвольной функции  $\psi (t)$,  монотонно убывающей к нулю при  $t \to +\infty$,
найдется набор  $\Theta = (\theta^1,...,\theta^m)$,
состоящий из чисел,
 линейно независимых вместе с единицей над
$\mathbb{Z}$ и  бесконечная последовательность векторов
 $$
(k_{1,\nu},...,k_{m,\nu}) \in {\rm spec}\, F,\,\,\, \nu = 1,2,3,...
$$
такие, что
для каждого 
$\nu$ выполняется 
\begin{equation}
||\theta^1
k_{1,\nu}+...+\theta^mk_{m,\nu}||\le
\psi \left(
\max_{1\le i\le m}|k_{i,\nu+1}|\right)
.
\label{ssi}
\end{equation}

\end{lem}

Видим, что условие (\ref{ssi})  очень похоже на соотношение (\ref{3}),
фигурирующее в эквивалентном определении $\psi$-сингулярности
(в лемме \ref{lkolo}  речь идет как бы о "наилучших приближениях из множества  ${\rm spec}\, F$").

Отметим, что если набор  $\Theta$  не является  $\psi $-сингулярным (при некоторой функции  $\psi$),
 то   получить утверждение о равномерной возвращаемости можно очень легко: имеет место следующее

{\bf  Предложение 4.}
\,\,\,{\it  

Рассмотрим убывающую к нулю при  $t\to +\infty$ функцию
$\Phi (t)$  и предположим, что ряд
$$
\sum_{k_1,...,k_m =-\infty}^{+\infty}\,\,\,
\Phi \left( \max_{1\le i\le m} |k_i|\right) \cdot  \max_{1\le i\le m} |k_i|
$$
сходится.
 Предположим, что коэффициенты Фурье функции  $F(x_1,...,x_m)$  допускают оценку
$$
|F_{k_1,...,k_m}|
\le \Phi \left( \max_{1\le i\le m} |k_i|\right),
$$
и что выполнено условие нулевого среднего  (\ref{zeromen1}).
Рассмотрим (стремящуюся к нулю) функцию
$$
\Phi_1 (t) =\sum_{(k_1,...,k_m): \max_{1\le i\le m} |k_i|\ge t}
\,\,\,
\Phi \left( \max_{1\le i\le m} |k_i|\right) 
$$
 Пусть $n=1$  и 
набор  $\Theta =(\theta^1,...,\theta^m)$,
состоящий из линейно независимых вместе с единицей над  $\mathbb{Z}$  вещественных чисел,
не является  $\psi$-сингулярным
ни при какой функции  $\psi (t)$  такой, что
$$
\psi (t)\cdot  \left(\Phi_1(t)\right)^{-\frac{1}{m}}\to+\infty,\,\,\,\, t\to +\infty.
 $$
Тогда
 найдется последовательность целых чисел  $q_\mu, \mu =1,2,3,..
$
 такая, что
$$
R^{[F,\Theta]}(q_\mu) \to 0,\,\,\,
\max_{1\le i \le m}||\theta^iq_\mu ||\to 0,\,\,\, \mu \to +\infty.
$$

}

Сделаем пояснение. Например,
условия предложения будут выполнены, если от коэффициентов Фурье  функции $F$ потребовать
выполненоя неравенства
$$
|F_{k_1,...,k_m}|
\le \Gamma (\max_{1\le i\le m}|k_i|)^{-\gamma},\,\,\,\,\
\forall {\bf k} \in \mathbb{Z}^m \setminus\{{\bf 0}\},
$$
с некоторыми
$$
\Gamma >0,\,\,\,\,\,\gamma > (m+1)m,
$$
 а от набора  $\Theta$  потребовать выполнения соотношения
$$
\limsup_{t\to+\infty}\psi_\Theta (t)\cdot t^{\frac{\gamma}{m}-1} =+\infty.
$$

Видим, что предложение 4 утверждает существование "сильной возвращаемости" когда выполнены одновременно
оба соотношения из (\ref{nadstro}),
причем эта возвращаемость носит равномерный характер.
Из доказательства очевидно, что
аналог предложения 4  справедлив и для одновременной возвращаемости интергалов нескольких функций.
В свете предложения 4 становится ясно, что основную трудность в доказательтсве теорем \ref{koz1},\ref{koz2}
составляет  разбор случая, когда  набор $\Theta$  является $\psi$-сингулярным с очень быстро 
убывающей функцией  $\psi$.

Доказательство предложения 4.

Невыполнение условия  $\psi$-сингулярности, согласно предложению 1, означает, что для некоторой бесконечной подпоследовательности  $\nu_\mu$ наилучших приближений выполнено
\begin{equation}\label{ooppv}
\zeta_{\nu_\mu}\cdot (\Phi_1(M_{\nu_\mu+1}))^{-\frac{1}{m}}
\to +\infty,\,\,\, \mu \to +\infty.
\end{equation}

  Ясно, что
$$
\lambda_{\nu_\mu} =
\zeta_{\nu_\mu}\cdot \Phi_1 (M_{\nu_\mu +1}) \to 0,\,\,\,
\mu \to +\infty.
$$
Положим  
$$Q_\mu = 
\left(\zeta_{\nu_\mu}\cdot 
 \Phi_1 (M_{\nu_\mu +1}) \right)^{-\frac{m}{m+1}}
$$
 и определим по теореме Дирихле такое  $q_\mu$, что
$$
1\le q_\mu \le Q_\mu,\,\,\,\,
\max_{1\le i \le m}||\theta_i q_\mu ||\le Q_\mu^{-\frac{1}{m}}.
$$
  Ясно, что
$$
q_\mu \to +\infty,\,\,\,
\max_{1\le i \le m}||\theta_i q_\mu ||
\to 0,\,\,\,\, \mu \to +\infty.
$$
 С другой стороны
$$
R^{[F,\Theta]}(q_\mu) \le
\sum_{{\bf k}\in \mathbb{Z}^m\setminus \{{\bf 0}\}}
\left|F_{k_1,...,k_m}\right|
\cdot \left|\sum_{s=1}^{q_\mu}
\exp (2\pi i (k_1\theta^1+...+k_m\theta^m)s )\right|
\ll
$$
$$
\ll
\sum_{{\bf k}: \max_i |k_i|< M_{\nu_\mu +1}}
\Phi \left(\max_i |k_i|\right)\cdot
\frac{\max_i |k_i|\cdot \max_i ||\theta_i q_\mu||}{||k_1\theta^1+...+k_m\theta^m||} +
q_\mu \cdot \Phi_1 (M_{\nu_\mu+1})\ll
$$
$$
\ll
\frac{1}{\zeta_{\nu_\mu}\cdot Q_\mu^{1/m}}+Q_\mu \cdot \Phi_1 (M_{\nu_\mu+1})
= 2\left(\zeta_{\nu_\mu}^m(\Phi_1(M_{\nu_\mu+1}))^{-1}\right)^{-\frac{1}{m+1}}
\to 0,\,\,\,
\mu \to +\infty.
$$
В последнем неравенстве использовано (\ref{ooppv}).  
 Предложение доказано.

\section{Сингулярные системы специального вида.}\label{specivid}

Некоторые утверждения из предыдущих пунктов нашей статьи представляют из себя теоремы существования сингулярных матриц, обладающих специальными свойствами. К таковым относятся, например, теоремы 
\ref{Ta},\ref{tritri},\ref{t:zeta_convergence_implies_degeneracy},\ref{dime2},\ref{suxx}
  и лемма \ref{lkolo}.

В настоящем пункте мы приведем еще два утверждения о существоавнии сингулярных систем со специальными свойствами. 

Первому утверждению следует предпослать  два замечания.  
Во-первых, если рассмотреть случай  $n=1,m=2$  то (как отмечалось в пунктах  \ref{ppnp},\ref{mravno2}
бесконечно много раз будут встречаться тройни подряд идущих независимых наилучших приближений, и, таким образом,
для бесконечно многих значений $\nu$  будет выполнено
$$
\zeta_\nu \ge \frac{1}{6M_{\nu+1}M_{\nu+2}}.
$$
Во-вторых, хотим отметить, что  для  сколь угодно быстро убывающей к нулю при $t \to +\infty$ 
функции $\psi (t)$
при $n=1$ и любом  $m\ge 2$ автором в \cite{MDAN}  было анонсировано (относительно подробное доказательство есть
в \cite{ME})  существование векторов
$\Theta\in\mathbb{R}^m$ таких, что ${\rm dim}_\mathbb{Z}\Theta = m+1$
и
$$
\zeta_\nu \le \psi(M_{\nu+m-1}).
$$
Последнее утверждение существование допускает существенное усиление.

\begin{theorem}\label{specii}
 Пусть  $m\ge 3$, Пусть задана
 сколь угодно быстро убывающая к нулю при $t \to +\infty$ 
функция   $\psi (t)$.
Пусть задана сколь угодно быстро растущая к бесконечности последовательтность натуральных чисел  $\tau (\nu), \nu =1,2,3,...$.
Тогда существуют 
$\Theta\in\mathbb{R}^m$ такие, что ${\rm dim}_\mathbb{Z}\Theta = m+1$
и
$$
\zeta_\nu \le \psi(M_{\nu+\tau(\nu)}).
$$
\end{theorem}

Мы поясним доказательство этой теоремы в случае  $m=3$.

Для векторов ${\bf w,e}\in\mathbb{R}^3$
(предполагаем, что  ${\bf e}$  есть вектор единичной длины)
и положительных $\eta,\delta$  определим шаровой сектор
\begin{equation}\label{vid}
B_{\eta,\delta}({\bf w,e}) =
\{
\Theta:\,\,\,
|\Theta -{\bf w}|\le \eta,
\,\,\,
|(\Theta-{\bf w})/|\Theta-{\bf w}|-{\bf e}|\le \delta\}
\end{equation}
($|\cdot|$  в данном случае
обозначает евклидову норму, хотя, вообще говоря, специфика нормы не имеет значения).
Для вектора ${\bf w} =(w^1,w^2,w^3) \in \mathbb{R}^3$  будем обозначать
$\overline{\bf w} =( w^1,w^2,w^3,1) \in \mathbb{R}^4$.
Будем отождествлять пространство  $\mathbb{R}^3$  с аффинным подпространством
$\{ ({\bf x},y):\, y =1\}\subset \mathbb{R}^4$.
Для доказательства теоремы \ref{specii}
используется следующая

\begin{lem}\label{specilem}
Пусть имеется трехмерное линейное вполне рациональное подпространство
$\Pi\subset \mathbb{R}^4, \,{\rm dim }\, \Pi = 3$  с нормалью  $\overline{\bf w} $ 
  и два линейных подпространства  $\pi_1,\pi_2 \subset 
\Pi,\,\, {\rm dim}\pi_1 ={\rm dim }\pi_2 = 2.$
 Пусть заданы натуральные числа
 $a<b<c$ и векторы
$$
{\bf z}_\nu = ({\bf x}_\nu,y_\nu) =
(x_{1,\nu},x_{2,\nu},x_{3,\nu},y_\nu)\in \mathbb{Z}^4,\,\,\,
1\le \nu \le c $$
 такие, что

{\rm (i)}
${\bf z}_a,{\bf z}_{a+1},... {\bf z}_{b-1},{\bf z}_b \in \pi_1$;

{\rm (ii)} 
${\bf z}_b,{\bf z}_{b+1},... {\bf z}_{c-1},{\bf z}_c \in \pi_2$;

{\rm(iii)}
${\bf z}_{b-1},{\bf z}_b,{\bf z}_{b+1}$  линейно независимы.

Пусть имеются положительные  $\eta,\delta$  и вектор  ${\bf e}\in\mathbb{R}^3$  такие, что выполнено следующее:
вектор  ${\bf e}$  ортогонален подпространству  $\pi_2$,  и
для любого 
$\Theta \in 
B_{\eta,\delta}({\bf w,e})
$,  с условием  ${\rm dim}_\mathbb{Z}\Theta = 4$,
первые $c$  штук вектров наилучших приближений суть  именно векторы  ${\bf z}_\nu,\, 1\le \nu\le c$.

Пусть задано натуральное  $ d >c$.

Тогда найдется вполне рациональное линейное подпространство
$\Pi_*, {\rm dim }\Pi_* = 3$
 с нормалью  $\overline{\bf w}_*$, подпространство
$\pi_3 \subset \Pi_*,\,{\rm dim}\,\pi_3 =2$ 
и последовательность целочисленных векторов
$$
{\bf z}_{c+1},...,{\bf z}_{d-1}, {\bf z}_d
\in \pi_3
$$
 такие, что
${\bf z}_{c-1}, {\bf z}_c,{\bf z}_{c+1}$  независимы, и 
для некоторых положительных
$\eta_*,\delta_*$  и некоторого вектора  ${\bf e}_*\in\mathbb{R}^3$  
 выполнено следующее:
вектор  ${\bf e}_*$  ортогонален подпространству  $\pi_3$,  и
для любого
$\Theta \in 
B_{\eta_*,\delta_*}({\bf w_*,e_*})\subset B_{\eta,\delta}({\bf w,e})
$, с условием  ${\rm dim}_\mathbb{Z}\Theta = 4$,
первые $d$  штук векторов наилучших приближений суть именно векторы  ${\bf z}_\nu,\, 1\le \nu\le d$.

  \end{lem}

Изложим схему доказательства леммы \ref{specilem}.  
Целочисленная решетка  $\mathbb{Z}^4$  раcпадается на трехмерные слои $\Pi_r$, параллельные подпространству 
 $\Pi$:
$$
\mathbb{Z}^4 = \bigcup_{r\in \mathbb{Z}} \Pi_r
$$
Сначала выбираем 
целую точку ${\bf z}_{c+1}$ из слоя  $\Pi_1$  или  $\Pi_{-1}$,  "соседнего"  к слою  $\Pi_0=\Pi$
(какой именно надо выбирать соседний слой определяется направлением вектора ${\bf e}$).
Можем считать, что абсолютная величина $M_{c+1} =\max_{i=1,2,3}|x_{i,c+1}|
$
 очень велика по сравнению со всеми параметрами, фигурирующими в индуктивной предпосылке.
Далее обозначаем  $\Pi_* = {\rm span} ({\bf z}_{c-1}, {\bf z}_c,{\bf z}_{c+1})$, 
$\pi_3 = {\rm span} ({\bf z}_c,{\bf z}_{c+1})$,
и пусть
$\overline{\bf w}_* $
 -  нормаль к  $\Pi_*$.

Векторы ${\bf z}_c,{\bf z}_{c+1}$  лежат в двумерном вполне рациональном подпространстве  $\pi_3$.
 Мы можем  найти последовательность векторов  $
{\bf z}_{c+2},...,{\bf z}_d\in \pi_3\cap\mathbb{Z}^4$  так, чтобы они составляли последовательность 
{\it  всех }
наилучших приближений
(больших по норме, чем  $|{\bf z}_{c+1}|_\bullet$)
 для последнего вектора  ${\bf z}_d$ 
{\it относительно решетки }
$\pi_3\cap \mathbb{Z}^4$  с индуцированной нормой.
Так как  $\Pi_*\supset \pi_3\ni {\bf z}_c, ...,{\bf z}_d$,
то за счет выбора малых  $\eta_*$
и $\delta_*$
 (в зависимости от  величин $\rho(\Pi_*)$ и $M_d = \max_{i=1,2,3}|x_{i,d}|$)  мы добиваемся выполнения утверждения леммы
\ref{specilem}.  При этом направление вектора  ${\bf e}_*$  выбирается, чтобы именно векторы
${\bf z}_{d-1}$  и  ${\bf z}_d$
(а не ${\bf z}_d-{\bf z}_{d-1}$  и  ${\bf z}_d$ )
 были бы последними наилучшими приближениями из первых $d$  приближений для любого  $\Theta$ из
строящейся окрестности.
(Подобная конструкция имеется у О.Н.Германа в \cite{sado},\cite{germa}.)

 Для доказательства теоремы \ref{specii} лемма \ref{specilem}
применяется следующим образом.  Сначала занумеро-вываются все трехмерные вполне рациональные подпространства в  $\mathbb{R}^4$.
Затем запускается индуктивный процесс построения вложенных окрестностей вида (\ref{vid}).
Определяется последова-тельность индексов
$\nu_k$  так, чтобы
 $\nu_{k+1} = \nu_k+\tau (\nu_k)$.
Каждый раз лемма \ref{specilem}
применяется для  $a= \nu_{k-1}, b=\nu_k, c= \nu_{k+1}$.
При этом дополнительно требуется выбирать величину  $\eta_*$  настолько малой, чтобы
$$
\eta_* < \frac{\psi (M_d)}{{8M_b}}.
$$ 
Тогда 
при
$b\le \nu\le c$ 
имеем  $\nu+\tau (\nu) \le c+\tau (c) = d$ и
будет выполняться
$$
\zeta_\nu
\le \zeta_b < \psi (M_d)\le \psi (M_{\nu+\tau (\nu)}).
$$
Кроме того, на каждом шаге  следует отделяться от вполне рационального трехмерного подпростра-нства, номер которого совпадает с номером шага.
 
Пользуясь случаем, приведем формулировку одного  утверждения, анонсированного 
Н.Г.Мощеви-тиным и О.Н.Германом 
в  \cite{sado}

\begin{theorem} \label{prop:psi_singularity_plus_linear_independence}
Пусть  $ n=1, m\ge  2$. 
Тогда для произвольной 
 сколь угодно быстро убывающей к нулю при $t \to +\infty$ 
функции   $\psi (t)$ 
найдется 
$\psi$-сингулярный    набор $\Theta \in \mathbb{R}^m$ такой, что
${\rm dim}_\mathbb{Z}\Theta = m+1$ и
 для любого  $\nu$  векторы
${\bf z}_\nu,{\bf z}_{\nu +1},\ldots,{\bf z}_{\nu+m}$ будут независимы.
\end{theorem}

\section{Добавление.}
\label{dobavl}
\subsection{О лакунарных последовательностях.}\label{lakey}

Последовательность $\{t_j\},\,\,\ j = 1,2,3,..$ положительных вещественных чисел называется лакунарной, если для некоторого $ M>0$ выполнено
\begin{equation}
\frac{t_{j+1}}{t_j} \ge 1+\frac{1}{M},\,\,\, \forall j \in \mathbb{N}. \label{aaaa} \end{equation}

А.Я.Хинчин в 1924 году
доказал, что для любой последовательности, удовлетворяющей (\ref{aaaa}) найдется вещественное число $\alpha$   и положительное $\gamma$ такие,
что   будет выполнено
$$
|| t_n \alpha ||\ge \gamma \,\,\, \forall n \in \mathbb{N}.
$$
Этот результат опубликован в работе \cite{HINS} (Hilfssatz III),
которую мы цитировали неоднократно. Его можно также найти в недавнем переиздании трудов А.Я.Хинчина \cite{H}
(работа \cite{HINS}, переведенная на русский язык, вошла в книгу \cite{H}).
Лемма \ref{YT2} из пункта \ref{ododo} настоящей работы, естественно, являлясь обобщением этого одномерного результата А.Я.Хинчина.

 Отметим, что конструкция А.Я.Хинчина позволяет установить
существование вещественного $\alpha$ и положительной абсолютной постоянной $\gamma$ таких, что
$$
|| t_n \alpha ||\ge \frac{\gamma}{(M\log M)^2} \,\,\, \forall n \in \mathbb{N}.
$$
50 лет спустя, в 1974 году
 П.Эрдеш
  \cite{E}
высказывал гипотезу, что для произвольной лакунарной последовательности найдется вещественное число  $\alpha$ такое, что дробные доли $ \{
\alpha t_j\} , \,\, j \in\mathbb{N}$ не будут всюду плотны в $[0,1]$.  Ответ на вопрос П.Эрдеша, естественно, следовал из процитированной выше
леммы А.Я.Хинчина. Однако этот результат А.Я.Хинчина был забыт. Решение гипотезы П.Эрдеша было опубликовано А.Поллингтоном \cite{PO} и Б.
де Матаном \cite{MA}. Позже количественные уточнения имелись у И.Кацнельсона \cite{K}, Р.К.Ахунжанова и Н.Г.Мощевитина \cite{MOS}, А.Дубицкаса
\cite{D}. Наилучший известный количественный результат принадлежит Ю.Пересу и В.Шлагу \cite{PS}. Они доказали, что с некоторой абсолютной положителной
постоянной
  $\gamma
> 0$
для каждой рассматриваемой последовательности $\{t_j\}$ найдется вещественное число $\alpha $ такое, что
\begin{equation}
||\alpha t_j  || \ge \frac{\gamma}{M\log M},\,\,\, \forall j \in \mathbb{N}.
\label{ququ}
\end{equation}
Ю.Перес и В.Шлаг использовали весьма оригинальную конструкцию, связанную с применением специального варианта локальной леммы Ловаса.

После работы 
Ю.Переса и В.Шлага появился ряд статей, в которых предложенный метод применялся в различных диофантовых задачах
(см. \cite{bumo},\cite{Db},\cite{MAR0},\cite{MAR1},\cite{MJNT},\cite{Mo4},\cite{Mo5},\cite{RO}).
В следующем пункте мы дадим краткий обзор некоторых из этих работ. Здесь же отметим, что в задаче о лакунарной последовательности наилучшее значение постоянной  $\gamma$  из  (\ref{ququ}), по-видимому, получается из работы И.П.Рочева \cite{RO}.

Предложенный Ю.Пересом и В.Шлагом в работе \cite{PS} подход
 к построению плохо приближаемых (по отношению к лакунарной последовательности) вещественных чисел оказался
применим во  многих задачах теории диофантовых приближений.  Ниже мы перечислим некоторые из них.
Но прежде мы сделаем ряд замечаний метрического характера.

\subsection{О некоторых метрических результатах}\label{memei}

Во-первых, рассмотрим одномерную задачу о приближениях.
Процитируем
результат Дж.В.С.Кас-селса из работы \cite{CAS}.
Напомним, что возрастающая последовательность
$\{t_n\}_{n=1}^{\infty}$ целых чисел называется
$\Sigma$-последовательностью, если, обозначая через $\mu_n$
количество дробей вида $\frac{j}{t_n}, 0{<}j{<}t_n$ не
представимых в виде $\frac{i}{t_q}$ с натуральным $i$ и $q{<}n$,
имеет место соотношение
$$
\liminf_{N\to\infty}\frac{1}{N}\sum_{n\le N}
\frac{\mu_n}{t_n}{>}0. \label{sigma}
$$
Примерами $\Sigma$-последовательностей могут служить
произвольная лакунарная последовательность натуральных чисел,
последовательность чисел  вида $ n^d, n =1,2,3,...$
при фиксированном натуральном $d$,
 или последовательность Фюрстенберга (о которой идет речь ниже).

Теорема Дж.В.С.Касселса \cite{CAS} утверждает, что если
$\{t_n\}_{n=1}^{\infty}$ является $\Sigma $-последовательностью
 и функция  $\psi (n)$   монотонно убывает к нулю при  $t\to+\infty$
то в случае расходимости ряда $\sum\limits_{n=1}^{\infty}\psi(n)$
 для почти всех $\xi $ 
выполняется
$$
\liminf_{tn\to+\infty} ||t_n\xi ||\cdot (\psi(n))^{-1}
=0.
$$
Во-вторых, сформулируем случай расходимости из теоремы А,Я. Хинчина \cite{h1926a}(в задаче о совместных приближениях и в том виде, в котором он присутствует в главе III книги \cite{SCH}). Пусть набор положительнозначных функций $\psi_1(q),...,\psi_n(q)$  натурального аргумента $q$  таков, что функций
$$
\psi(q) =\prod_{j=1}^n\psi_j(q)
$$
не возрастает, а ряд
$
\sum_{q=1}^{+\infty} \psi (q)
$
 расходится, то для почти всех наборов
$(\theta_1,...,\theta_n)\in \mathbb{R}^n$
 имеется бесконечно много натуральных чисел $q$,
таких что
$$
||q\theta_j||<\psi_j(q),\,\,\,\,1\le j\le n.
$$

В-третьих упомянем результат П.Галлагера из \cite{Gal62} (точнее не общее утверждение, а его следствие),
тоже касающийся случая расходимости. Этот результат гласит, что для почти всех пар вещественных чисел
$\xi_1,\xi_2$
выполнено
$$
\liminf_{n\to\infty} n\log^2 n \, ||n\xi_1||\, ||n\xi_2 || =0.
$$

\subsection{О применении метода Переса-Шлага.}

Здесь мы обрисуем круг задач, в которых применение метода Переса-Шлага позволяет
получать нетривиальные результаты о существовании плохо приближаемых чисел.

$\bullet$\, {\bf Приближения с субэкспоненциальными последовательностями.}

Рассмотрим последовательность вещественных чисел  $t_n, n = 1,2,3,..$  такую, что
\begin{equation}
t_n \asymp \exp(\gamma t^\beta ),\,\,\, \gamma > 0,\,\,\, 0<\beta<1.
\label{eexp}
\end{equation}
В \cite{MAR1}  показано, что для любой последовательности вещественных чисел
$\eta_n$
 множество
$$ 
\{ \xi \in \mathbb{R}: \, \inf_{n\in \mathbb{N} }||t_n\xi  -\eta_n||\cdot t^{1-\beta} \log t>0\}
$$
имеет размерность Хаусдорфа равную 1
(точнее, в \cite{MAR1} рассмотрен только однородный случай $\eta_n  = 0 \,\,\,\forall n$,
однако доказательство общего случая дословно повторяет доказательство частного случая).

Например, рассмотрим  последовательность Фюрстенберга
$s_n, n = 1,2,3,...$,
 состоящую из целых чисел вида
$2^k3^m
, k, m = 0,1,2,3,...$, упорядоченных по возрастанию.
 Для нее условие (\ref{eexp})  выполнено с
$\beta = 1/2$.
 Таким образом, получаем, что
 при всякой фиксированной последова-тельности неоднородностей $\eta_n$
 множество
$$ 
\{ \xi \in \mathbb{R}: \, \inf_{n\in \mathbb{N} }||s_n\xi  -\eta_n|| \cdot t_n^{1-\beta} \log t_n>0\}
$$
 имеет
хаусдорфову размерность 1.

Г.Фюрстенберг
\cite{FU} (простое доказательство  дал М.Бошерницан  \cite{BOSH}) доказал, что для иррацио-нального  $\xi$  множество
дробных долей $\{ 2^n3^m\alpha\}$ на отрезке $[0,1]$ всюду плотно и, стало быть,
$$
\liminf_{n\to\infty}||s_n \xi ||{=}0.
$$
Отметим, что о характере стремления величины  $||s_n\alpha ||$ к
нулю замечательный результат получили недавно Ж.Бургейн, Е.Линденштраусс, П.Мишель и А.Венкатеш \cite{JBL}.
Этот результат связан с применением оценок снизу линейных форм от двух логарифмов алгебраических чисел.

$\bullet$ \,{\bf Полиномиальные последовательности. }

Пусть
последовательность  $t_n$  растет полиномиально, то есть
$$
t_n \asymp n^\beta,\,\,\, \beta >0,$$

В работе \cite{MJNT}  доказано, что множество
$$ 
\{ \xi \in \mathbb{R}: \, \inf_{n\in \mathbb{N} }||t_n\xi  -\eta_n|| \cdot t_n \log t_n>0\}
$$
имеет размерность Хаусдорфа  $\ge \frac{1}{1+\beta}$.

 Этим дан ответ на вопрос В.М. Шмидта
из работы  \cite{SCHP}.
 В этой же работе В.М. Шмидт задает более сложный вопрос, на который ответ до сих про неизвестен. Мы приводм его формулировку: 
существуют ли такие    вещественные $\xi$, для которых
$$
 \inf_{n\in \mathbb{N} }||n^2\xi|| \cdot n >0
\,\,
?$$
Приведем здесь также формулировку результата А.Захареску \cite{Z} о том, что для любого иррационального  $\xi$
 и для любого положительного  $\varepsilon$  выполнено
$$
\liminf_{n\to \infty} ||n^2\xi ||\cdot n ^{2/3-\varepsilon } = 0.
$$
Результат А.Захареску является наилучшим известным на настоящее время.

Упомянутая выше теорема Касселса оказывает, что построенные в двух рассмариваемых выше задачах множества могут иметь лишь нулевую меру Лебега. 
Хочется
также отметить, что в работе И.П.Рочева  \cite{RO}  имеется общая теорема о  приближениях линейной формы,
содержащая в себе изложенные выше результаты о существовании плохо приближаемых чисел.

Далее в текущем пункте мы остановимся на двух многомерных задачах.

$\bullet$\,\,{\bf  Проблема Литтлвуда.}

Весьма знаменитая гипотеза Литтлвуда предполагает, что для любых вещественных чисел
$\xi_1,\xi_2$
выполняется
$$
\liminf_{n\to\infty} n\, ||n\xi_1||\, ||n\xi_2 || =0.
$$
Пользуясь случаем, автор хочет сказать,
что формулировку гипотезы Литтлвуда он узнал, будучи студентом второго курса: ему рассказал о ней С.А.Довбыш.

С помощью метода Переса-Шлага можно доказать (см. \cite{Mo5}) существование таких
$\xi_1,\xi_2$
выполняется
$$
\liminf_{n\to\infty} n\log^2 n \cdot \, ||n\xi_1||\, ||n\xi_2 ||>0.
$$
 Из упомянутой в начале пункта теоремы П.Галлагера вытекает, что наборы $\xi_1,\xi_2$,
обладающие таковым свойством, могут образовывать в  $\mathbb{R}^2$ множество не более чем нулевой меры Лебега.

$\bullet$\,\,{\bf BAD-гипотеза.}

 Для чисел $\alpha, \beta \in [0,1] $ при условии $ \alpha +\beta = 1 $ и для $ \delta \in (0,1/2)$
 рассмотрим множество
$$
{\rm BAD}(\alpha, \beta ;\delta ) = \left\{\xi = (\xi_1,\xi_2 ) \in [0,1]^2:\,\,\,\inf_{p\in \mathbb{N}} \,\, \max \{ p^\alpha ||p\xi_1||,
 p^\beta ||p\xi_2||\} \ge \delta \right\}
.$$
В.М.Шмидт \cite{SCHL} 
высказал предположение о том, что для любых двух наборов  $\alpha_1,\alpha_2, \beta_1,\beta_2 \in [0,1],\, \alpha_1 +\beta_1 =\alpha_2+\beta_2 =1$
 пересечение
  $$
{\rm BAD}(\alpha_1, \beta_1  ) \bigcap {\rm BAD}(\alpha_2, \beta_2  )
$$
непусто. Это предположение, двойственное  к гипотезе Литтлвуда, получило название  BAD-гипотезы.
Оно долгое время оставалось открытым.
Недавно Д.Бодягин, А.Поллингтон и С.Велани анонсировали решение BAD-гипотезы.
Полное (или даже схематическое) доказательство пока недоступно. 

Метод Переса-Шлага дает следующий результат.
Для произвольно последовательности неоднородностей $\eta =\{\eta_j\}_{j=1}^\infty $
Рассмотрим
множества
$$
{\rm BAD}_\eta^* (\alpha, \beta ;\delta ) = \left\{\xi = (\xi_1,\xi_2 ) \in [0,1]^2:\,\,\,\inf_{p\in \mathbb{N}} \,\, \max \{ (p\log( p+1))^\alpha
||p\xi_1||,
 (p\log (p+1))^\beta ||p\xi_2-\eta_p||\} \ge \delta \right\}.
$$
Понятно, что
$$
{\rm BAD} (\alpha, \beta ;\delta ) \subseteq {\rm BAD}_0^* (\alpha, \beta ;\delta )$$
($\eta =0$  здесь означает, что все $\eta_j$ равны нулю).
и что, поскольку ряд 
$ \sum_{p=1}^\infty \frac{1}{p\log (p+1)}$ 
расходится, при каждой фиксированной последовательности  $\eta$ объединение
$$
{\rm BAD}_\eta^*(\alpha, \beta  ) = \bigcup_{\delta > 0} {\rm BAD}_\eta^*(\alpha, \beta ;\delta )
$$
есть множеество нулевой меры Лебега

В \cite{Mo4}  доказано, что
 для произвольных 
 $\alpha_1,\alpha_2, \beta_1,\beta_2 \in [0,1],\, \alpha_1 +\beta_1 =\alpha_2+\beta_2 =1$ и для произвольного $ 0<\delta \le 2^{-20} $
 пересечение
\begin{equation}
 {\rm BAD}_0^*(\alpha_1, \beta_1 ;\delta ) \bigcap {\rm BAD}_0^*(\alpha_2, \beta_2 ;\delta )
\label{baad}
\end{equation}
непусто.
Доказательство 
утверждения о непустоте пересечений вида
$$
 {\rm BAD}_{\eta^1}^*(\alpha_1, \beta_1 ;\delta ) \bigcap {\rm BAD}_{\eta^2}^*(\alpha_2, \beta_2 ;\delta )
$$
для любых двух фиксированных последовательностей  $\eta^1,\eta^2$ дословно повторяет доказательстов
соотношения (\ref{baad}),  проведенное в \cite{Mo4}.

Как показал Я.Бюжо в совместной работе \cite{bumo}, во всех перечисленных выше задачах
можно доказать существование множества плохо приближаемых чисел {\it  полной}
размерности Хаусдорфа. 
Для этого следует применить принцип распределения масс (см. \cite{bubu}, глава V или \cite{falco}).
 Кроме этого,  работа \cite{bumo}  содержит ряд других применений метода к задачам
типа проблемы Литтлвуда.

 В настоящей статье мы не предполагаем делать сколь-нибудь подробный обзор работ, связанных с
проблемами типа гипотезы Литтлвуда.  Таковых работ очень много.
Заметим лишь, что имеется большое количество работ посвященным исследованиям обсуждавшихся выше вопросов с помощью теории динамических систем. Мы упомянем обзор \cite{oogoro}.

Отметим, что в доказательствах всех сформулированных ы настоящем путкте результатов о существовании плохо приближаемых чисел так или иначе  присутствует следующий вариант локальной леммы Ловаса (см. \cite{RO}):

\begin{lem}\label{lll}Пусть $\{A_n\}_{n=1}^N$ --- система событий в
вероятностном пространстве $(\Omega,\mathcal F,{\bf P})$,
$\{x_n\}_{n=1}^N$ --- набор чисел из $[0;1]$. Обозначим
$B_0=\Omega$, $B_n=\bigcap\limits_{m=1}^nA_m^c$ ($1\le n\le N$),
где $A_m^c=\Omega\setminus A_m$. Пусть для каждого
$n\in\{1,\ldots,N\}$ найдется такое $m=m(n)\in\{0,1,\ldots,n-1\}$,
что выполнено
$$
{\bf P}(A_n\cap B_m)\le x_n\prod_{m<k<n}(1-x_k)\cdot{\bf  P}(B_m)
$$
(если $m=n-1$, то считаем $\prod\limits_{m<k<n}(1-x_k)=1$). Tогда
при $1\le n\le N$ справедливо неравенство
$$
{
\bf P} (B_n)\ge(1-x_n){\bf P}(B_{n-1}).
$$
\end{lem}

\subsection{Об   $(\alpha,\beta)$-играх. }\label{uinni}

     Напомним определение $(\alpha,\beta)$ игры Банаха-Мазура-Шмидта   в пространстве  $\mathbb{R}^d$.
 Итак, пусть  $\alpha,\beta \in (0,1)$  и пусть  $S\subseteq \mathbb{R}^d$ --  некоторое множество.
Играют белые и черные.
 Сначала черные выбирают
 замкнутый шар (в sup-норме) $B_1$
диаметра $l(B_1){=}2\rho$. Затем белые выбирают замкнутый шар
$W_1{\subset}B_1$ диаметра $l(W_1){=}\alpha l(B_1)$. Потом чёрные
выбирают замкнутый шар $B_2{\subset}W_1$ диаметра $l(B_2){=}\beta
l(W_1)$, и т.д. В результате образуется
  последовательность вложенных замкнутых
шаров
$B_1{\supset}W_1{\supset}B_2{\supset}W_2{\supset}{\cdots}$ с
диаметрами  $l(B_i){=}2(\alpha\beta)^{i{-}1}\rho$ и
$l(W_i){=}2\alpha(\alpha\beta)^{i{-}1}\rho\;(i{=}1,2,\dots)$
Очевидно, что множество
$\cap_{i{=}1}^{\infty}B_i{=}\cap_{i{=}1}^{\infty}W_i$ состоит из
одной единственной точки. Если $\cap_{i=1}^{\infty}B_i{\in}S$, то
говорят, что белые выиграли игру. Множество $S$ называется
$(\alpha,\beta)$-выигрышным, если белые могут выиграть игру
независимо от того, как играют чёрные. Множество $S$ называется
$\alpha$-выигрышным, если оно является $(\alpha,\beta)$-выигрышным
при каждом $\beta \in (0,1)$.
 Выигрышной размерностью ${\rm windim }{\cal A}$ множества  ${\cal A}\in \mathbb{R}$  называется супремум тех  $\alpha$, при которых множество  ${\cal A}$  является $\alpha$-выигрышным.   Если  
${\rm windim }{\cal A}> 1/2$ то, как легко видеть  ${\cal A} =\mathbb{R}^d$. 
При  $0<\alpha <1/2$  существуют нетривиальные $\alpha$-выигрышные множества.
В.М.Шмидт(\cite{SCH1},\cite{SCH}))  доказал о выигрышных множествах ряд общих теорем.  
Ниже мы приводим две из них.

\begin{theorem}\label{T02}
  Если $ \alpha > 0$ , то всякое $\alpha$-выигрышное множество   в $\mathbb{R}^d$ имеет   размерность Хаусдорфа
 равную  $d$.
\end{theorem}

\begin{theorem}\label{T002}
 Если $ \alpha > 0$ , то пересечение любого
счетного числа $\alpha$-выигрышных множеств снова будет
$\alpha$-выигрышным.
\end{theorem}

Далее мы приводим ряд примеров выигрышных множеств.

 $\bullet $ 
\,\,  
 Следующий результат имеется у В.М.Шмидта  \cite{SCH1}.

\begin{theorem}\label{T0002}
Пусть $d = 1$.
 Зафиксируем натуральное  $ q \ge 2$.
  Множество
$$
{\cal N}_q = \{ x\in \mathbb{R}:\,\,\, \exists C(x) > 0\,\,\,
\forall n \in \mathbb{N}\,\,\, ||q^n x|| \ge C(x)\},
$$
состоящее из  вещественных чисел, не являющихся нормальными по натуральному
основанию $q \ge 2$ имеет выигрышную размерность  $1/2$.
\end{theorem}

$\bullet$
\,\,
При
  $d = 1$  рассмотрим лакунарную последовательность положительных чисел  $t_j$
(определение лакунарности см. в предыдущем пункте).  Обобщая рассуждения В.М.Шимдта,
 легко видеть, что множество
$$
{\cal N}= \{ x\in \mathbb{R}:\,\,\, \exists C(x) > 0\,\,\,
\forall  \in \mathbb{N}\,\,\, ||t_n x|| \ge C(x)\},
$$
 имеет выигрышную размерность  $1/2$.

$\bullet$
\,\,
Для последовательностей сублакунарного роста автор \cite{Mos}
 доказал следующий результат.

\begin{theorem}\label{T00002}
 Пусть
последовательность положительных чисел
 $t_n$ такова, что
$$\forall \varepsilon > 0 \,\,\, \exists N_0 \,\,\, \forall n
\ge N_0\,:\,\,\, \frac{t_{n+1}}{t_n}\ge 1+\frac{1}{n^\varepsilon}.
$$
Тогда при каждом положительном $\delta$ множество
$$
{\cal A}_\delta = \{ x \in \mathbb{R} \, :\,\,\, \exists c(x)>
0\,\,\, \forall n \in \mathbb{N}\,\,\,
 ||t_nx||> c(x)/n^\delta \}
$$
будет $\alpha$-выигрышным с любым заданным наперед $\alpha \in
(0,1/2)$, то есть ${\rm windim  }{\cal A}_\delta = 1/2$.
\end{theorem}

$\bullet$ \,\,
 Особо следует поговорить о плохо приближаемых системах линейных форм.
 Рассматриваемый в настоящей статье набор вещественных чисел   
$
\Theta =\{ \theta^i_j,\,\,\, 1\le i \le m,\,\,\, 1\le j\le n\}
$
называется {\it  плохо приближаемым}, если  найдется положительная постоянная  $\gamma$  такая, что для соответствующей системы линейных форм
$
{\bf L}_\Theta ({\bf x}) = \{L_j ({\bf x}),\,\,\, 1\le j\le n\},
$  выполнено
$$
\max_{1\le j\le n} ||L_j({\bf x})||\ge \gamma \cdot \left(\max_{1\le i \le m }|x_i|\right)^{-\frac{m}{n}},\,\,\,\,
\forall {\bf x} \in \mathbb{R}^m \setminus \{{\bf  0}\}.$$
 Хорошо известно, что плохо приближаемые наборы  $\Theta$
существуют, и что при всех   $n,m$  множество всех прохо приближаемых наборов  $\Theta$
 есть множество нулевой меры Лебега в  $\mathbb{R}^{nm}$
(последнеее вытекае, например, из метрической теоремы
А.Я.Хинчина, упомянутой в пункте \ref{meemei}).

С другой стороны, как легко видеть набор $\Theta$  будет плохо приближаем в том и только том случае, когда для функции Ярника  (\ref{fj}) выполняется
$$
\liminf_{t\to+\infty} t^{\frac{m}{n}}
\psi_\Theta (t) >0.
$$
 Ясно, что плохо приближаемый набор  $\Theta$  образует регулярную систему (в теоминологии А.Я.Хинчина).

  Сформулируем теорему В.М.Шмидта из \cite{Schforms}.

\begin{theorem}\label{0000000T2}
 При каждом значении размерностей  $n,m$  множество плохо приближаемых наборов  $\Theta$ в $\mathbb{R}^d, d = nm$
 имеет выигрышную размерность  $1/2$.
\end{theorem}

 Отметим, что плохо приближаемым системам чисел посвящено огромное количетво литературы,
обзор которой мы просто не в состоянии сделать.

$\bullet$  Отметим в заключение, что о выиргышных множествах имеется довольно много других работ
(см., например,  
\cite{ah},\cite{ah1},\cite{cantor},\cite{dre},\cite{tsengA},\cite{FISH},\cite{fb},\cite{fb1},\cite{KLEI},\cite{Ts}).

\vskip+2.0cm

\end{document}